%% file: dual.tex
\documentclass[10pt,a4paper]{article}

\usepackage[latin1]{inputenc}
\usepackage[T1]{fontenc}
\usepackage[english]{babel}
\usepackage{amsfonts}
\usepackage{euscript,amsthm,pstricks}
\usepackage[vcentermath]{youngtab}

\newtheorem{conj}{Conjecture}[section]
\newtheorem{theo}[conj]{Theorem}

\newtheorem{prop}{Proposition}[section]
\newtheorem{coro}[prop]{Corollary}
\newtheorem{lemm}{Lemma}[section]

\newtheorem{clai}[lemm]{Claim}
\newtheorem{fait}[lemm]{Fact}
\newtheorem{defi}{Definition}[section]
\newtheorem{nota}[defi]{Notation}
\newtheorem{exem}[defi]{Example}
\newtheorem{naiv}[defi]{Naive guess}
\newenvironment{proo}{{\flushleft \bf Proof :}}{\hfill $\square$ \vspace{2mm}}
\newenvironment{rema}{{\flushleft \bf Remark :}}{\hfill}
\begin{document}

\input commandes_pe.tex

\title{Dual varieties of subvarieties of homogeneous spaces}
\author{Pierre-Emmanuel Chaput}
\maketitle

{\def\thefootnote{\relax}
\footnote{\hskip-0.6cm
{\it AMS mathematical classification \/}: 14N99, 14L35, 14L40.   }}

\begin{center}
{\bf Abstract}
\end{center}

If $X \subset \p^n_\C$ is an algebraic complex projective variety, one
defines the dual variety $X^* \subset {(\p^n)}^*$ as the set of
tangent hyperplanes to $X$.
The purpose of this paper is to generalise this notion when $\p^n$ is
replaced by a quite general partial flag variety. A similar
biduality theorem is proved, and the dual varieties of Schubert varieties
are described.
\lpara

\begin{center}
{\bf Introduction}
\end{center}

Let $X \subset \p V$ be a complex projective algebraic variety, with
$V$ a $\C$-vector space. If $h \in \p V^*$ is a hyperplane and $x \in X$
is a smooth point, we say that $h$ is tangent to $X$ at $x$ if $h$
contains the embedded tangent space of $X$ at $x$. Equivalently, the
intersection $X \cap h$ is singular at $x$. The closure of
the set of all $h \in \p V^*$
which are tangent at some smooth point of $X$
is denoted $X^*$ and called the
dual variety of $X$; for given $h \in X^*$, the closure
of the set of smooth points $x \in X$ such that $h$ is tangent at $x$
is called the tangency locus of $h$.

This notion of dual varieties is a very classical one, and it is used
plentifully in both classical and modern articles. The very powerfull
biduality theorem, to the effect that ${(X^*)}^* = X$,
and its corollary, which states that the tangency locus
at a smooth point $h \in X^*$ is a linear space, are ubiquitous.
To state only one example, this result is crucial in Zak's classification
of Severi varieties, since it allows proving that the entry locus of
a Severi variety is a smooth quadric \cite[proposition IV.2.1]{zak}.

This biduality
theorem deals with subvarieties of projective space, which have
been studied by so many classical algebraic geometers. More recently, work has
been done in a new direction which consists in
considering subvarieties of other homogeneous
spaces. For example, G. Faltings \cite{faltings}
and O. Debarre \cite{debarre1,debarre2} have shown
some connectivity theorems that hold
in an arbitrary homogeneous space, E. Arrondo has proved a classification
of some subvarieties of Grassmannians similar to Zak's result
\cite{arrondo_severi}, and some
topological results on low-codimensional subvarieties of some homogeneous
spaces emerge in works of E. Arrondo - J. Caravantes 
\cite{arrondo_low} and N. Perrin \cite{perrin_low}.

\lpara

Obviously, to study subvarieties of homogeneous spaces, a similar
notion of dual variety and a biduality theorem are lacking. The aim of
this article is to fill this gap as much as it is possible.

Since homogeneous spaces $G/P$ are by definition projective algebraic
varieties, it is certainly possible to embed them in a projective space, and
therefore a subvariety $X \subset G/P$ is a fortiori a
subvariety of a projective space, so that one can consider
the usual dual variety of $X$.

However I claim that in many cases this is not the best thing to
do. Let us consider an example. Let $V$ be a $\C$-vector space
equipped with a non-degenerate quadratic form. If $Q \subset \p V$ is 
the smooth
quadric it defines, then it is well-known that $Q^* \subset \p V^*$ is also a
smooth quadric, canonically isomorphic with $Q$. 
Now, let $r$ be an integer and let us consider the variety $\G_Q(r,V)$
parametrising $r$-dimensional isotropic subspaces as a subvariety of
a suitable projective space using Plücker embedding. Then clearly we no
longer have $\G_Q(r,V)^* \simeq \G_Q(r,V)$. On the contrary,
let us consider $\G_Q(r,V)$ as a subvariety of the Grassmannian
$\G(r,V)$;
proposition \ref{prop_symmetric}
shows that for my definition of dual varieties, $\G_Q(r,V) \subset \G(r,V)$ 
has a well-defined
dual variety in $\G(r,V^*)$ which is canonically isomorphic with
$\G_Q(r,V)$.

In fact, homogeneous spaces are often minimally
embedded in projective spaces of very big dimension, so that the usual
dual variety of a subvariety of a homogeneous space will happen to be
very large and often untractable. A notion of dual varieties within
homogeneous spaces is probably more suitable if one wants to deal with
low-dimensional or low-codimension subvarieties (of course, the price to
pay is that the ambient space is a bit more complicated than a projective
space).

\lpara

My definition of dual varieties uses a class of birational
transformations called stratified Mukai flops by Namikawa
\cite{namikawa}. These are birational maps
$\mu:T^*G/P \dasharrow T^*G/Q$ defined in terms of nilpotent orbits
for some semi-simple group $G$ and some parabolic subgroups
$P,Q$. For
given $G,P,Q$, if there exists such a map, then we say that
$G/P$ and $G/Q$ allow
duality. For $X \subset G/P$, we consider its conormal bundle 
$N^*X \subset T^*G/P$ and define the dual variety 
$X^Q = \pi_Q \circ \mu (N^*X) \subset G/Q$ ($\pi_Q : T^* G/Q
\rightarrow G/Q$ denotes the projection) if $N^*X$ meets the locus
where $\mu$ is defined (in which case we say that $X$ is suitable). 
For example, if $G = SL(V)$ and
$G/P = \p V$, the only possibility for $Q$
leads to $G/Q = \p V^*$; any
proper subvariety $X \subset \p V$ will be suitable and $X^Q = X^*$.
Another example is the fact that a Grassmannian $\G(r,V)$ and its dual
Grassmannian $\G(r,V^*)$ allow duality, as one could naturally expect.

One advantage of this definition is that it uses the so-called
Springer resolutions of the corresponding nilpotent orbit, which are
symplectic resolutions, and this
article uses heavily informations which come from the study of such
resolutions \cite{namikawa,flop_scorza}.
Another advantage is that it exhibits
the symplectic nature of dual varieties. In fact, 
$T^*G/P$ and $T^*G/Q$
are symplectic varieties and $N^*X$, as a subvariety of
$T^* G/P$, is a Lagrangian subvariety. These properties suffice to
show very easily the biduality theorem \ref{theo_bidualite} in our
setting.

\lpara

However, this definition also has its drawbacks. The most important is
probably that it is not so much intuitive, so that given $x \in X$ and
$h \in X^Q$, it is not obvious at all what the sentence ``$h$ is
tangent to $X$ at $x$'' should mean. However, 
in the case of a Grassmannian,
using the natural rational 
map $Hom(\C^r,V) \dasharrow \G(r,V)$, I show that the
dual variety of $X \subset \G(r,V)$ can be computed in terms of the usual
dual variety of an adequate subvariety of $\p Hom(\C^r,V)$ (see subsection
\ref{subsection_grassmannienne}).
Therefore, this is a way of understanding more easily
dual varieties in the case of Grassmannians.
In general however, there are two fundamental
differences between our setting and usual duality.

First of all, given $G,P$, there may be many different $Q$'s, or none,
such that $G/P$ and $G/Q$ allow duality. Therefore, given suitable 
$X \subset G/P$, we will get a dual variety $X^Q$ for each such $Q$. If one
restricts to maximal parabolic subgroups, thanks to
\cite{namikawa}, this difficulty disappears
because for given $G/P$ there will be at most one parabolic subgroup 
$Q$ such that
$G/P$ and $G/Q$ allow duality. Moreover, section
\ref{section_reduction} shows that one can understand all dual
varieties if they are understood when $P$ and $Q$ are maximal parabolic
subgroups. These cases are therefore called fundamental cases. They
include the duality between the Grassmannian $\G(r,V)$ and its dual
Grassmannian $\G(r,V^*)$, but also a duality
between the two spinor varieties of a quadratic space of dimension
$4p+2$, and between the exceptional homogeneous spaces 
$E_6/P_1 \leftrightarrow E_6/P_6$ and $E_6/P_3 \leftrightarrow E_6/P_5$.

The second difference is that not all proper subvarieties 
$X \subset G/P$ will have a dual variety. Note that $X = \p V$ has no
dual variety in $\p V^*$, because for any $x \in X$,
no non-zero cotangent form can vanish
on $T_xX$. From this point of view, the situation is
quite similar in our setting~: too big subvarieties $X$ of $G/P$ don't have
dual varieties because for any $x \in X$
there is no generic cotangent form in $T^*_x G/P$ which
vanishes on $T_xX$.

\lpara

In the classical setting, a hyperplane $h$ is tangent to $X$ at $x$
\iff the intersection $h \cap X$ is singular.
As I already alluded to, I have not been able to give a similar
geometric notion of ``tangent element''. The only sensible definition
seemed to state that $h \in X^Q$ is tangent to $X$ at $x \in X$ if $h$
belongs to the image of $N^*_xX$ under $\pi_Q \circ \mu$.
Since there is an incidence
variety in $G/P \times G/Q$ (the closed $G$-orbit), any $h \in G/Q$
still defines, exactly as in the classical situation, a subvariety $I_h
\subset G/P$. Lemma \ref{lemm_q} implies that if $h$ is tangent to $x$
at $X$, then the intersection $I_h \cap X$ is not transverse, but the
reciprocal of this fact is false.

Section \ref{section_tangency} deals with this matter. Corollary
\ref{coro_incident} states that if $h$ is tangent to $X$ at $x$, then
$x \in I_h$. For $x \in X$ with $X$
suitable, the tangent cone $\overline{T_xX} \subset G/P$ of $X$ at $x$
is defined in a roundabout manner as the dual variety 
of the variety of $h$'s in $X^Q$ which are
tangent to $X$ at $x$. Theorem \ref{theo_cone} implies that
$\overline{T_xX}$ is a ``cone'' with vertex $x$, where
definition \ref{def_cone}
generalises the classical notion of
cones from subvarieties of projective space to subvarieties of
fundamental homogeneous spaces.

\lpara

Finally, section \ref{section_exemple} studies dual
varieties of Schubert varieties. In the classical setting, the dual
variety of a linear subspace is again a linear subspace. In our
setting, it is a formal consequence of the definitions that the dual
variety of a Schubert variety is again a Schubert
variety (see proposition
\ref{prop_schubert} which relies on the functorial property of dual
varieties given in proposition \ref{foncteur}).

Let $B \subset G$ be a Borel subgroup. It turns out that the combinatorial
involution $X \mapsto X^Q$ between $B$-stable suitable
Schubert subvarieties of $G/P$ and $G/Q$
is no longer decreasing, as it was the case for $G/P = \p V$.
For this reason, the description of this map is quite intricate. In
the case of Grassmannians and spinor varieties, we give explicitly
in terms of
partitions the map $X \mapsto X^Q$,
see propositions \ref{grass_adapte} and \ref{prop_spinoriel}.
This relies on a general
recipy for finding $X^Q$ when $X$ is a Schubert variety
which is given in subsection
\ref{subsection_combinatoire}. For the exceptional cases, this recipy
theoretically defines the involution (there is only a finite number of
calculations to do to compute the dual variety of a Schubert
subvariety), but I will not give a more
explicit description of it. As a first step, I describe a criterion
for a Schubert subvariety to be suitable. Remarkably enough, this
criterion can be stated in a uniform way for all the fundamental
cases, using the combinatorics of some quivers studied in
\cite{carquois}~: see theorem \ref{theo_adapte}.

\lpara
{\bf Further questions :}
Of course this study only gives basic properties of our generalised
dual varieties~: if one compares with usual dual varieties, what
essentially has been proved is the biduality theorem and the
computation of the dual variety of a quadric and a linear subspace. The
power of the classical notion of dual varieties gives hope to me that much more
can be said on this topic, including~:
\begin{itemize}
\item
Is it true that for a smooth subvariety $X \subset G/P$
one has $\dim\ X^Q \geq \dim\ X$ ? This question has been raised by
Laurent Manivel.
\item
Many Fano 3-folds are defined as subvarieties of some homogeneous
spaces. What are the dual varieties of these Fano 3-folds ?
\item
What is the dual variety of a divisor in $G/P$ ? If this is a divisor,
what is the degree of this divisor ? The
answer to this question for $G/P = \G(2,V)$ or $G/P = E_6/P_1$ and a
divisor of degree 1 has been given in \cite{hermitian} (the dual
variety is again a divisor of degree 1).
\item
Classification problems~: for example find all smooth varieties
with dual variety a divisor of low degree.
\end{itemize}

\lpara
\noindent
{\bf Acknowledgements : }I am very grateful to M. Brion for suggesting
that maybe nilpotent orbits could help defining an
interesting equivariant rational map
$T^*G/P \dasharrow G/Q$, as it was finally exactly the case. Thanks
are also due to B. Fu who pointed to me the reference \cite{namikawa}.

\tableofcontents

\sectionplus{Definition of the dual variety}

\label{section_definition}

\subsectionplus{Notations and definition}

In this subsection, I introduce the (abstract) definition of dual 
varieties, which allows easy proofs of general results; in subsection
\ref{subsection_grassmannienne}, an equivalent but more ``down-to-earth'' 
definition will be given in the case of Grassmannians.

\lpara

Before giving this definition, which is not so intuitive, I give some
``naive guesses'' and explain why the corresponding notion of dual
varieties would not be interesting. In this way, I hope to convince
the reader that it is not possible to avoid some technicalities.
Let us try our unsuccessful experiments in the case of Grassmannians.

So assume $G/P = \G(r,V)$ and $G/Q = \G(r,V^*)$ and assume
$2r<\dim V$. Any element
$h \in \G(r,V^*)$, representing a codimension $r$ 
subspace of $V$ denoted $L_h$, defines (at
least) two subvarieties in $\G(r,V)$. The first (resp. the second)
is the subvariety of $x \in \G(r,V)$ such that
$L_x \subset L_h$ (resp. $\dim(L_x \cap L_h)>0$). It will be denoted
$I_h$ (resp. $h^\bot$).
Assume $X \subset \G(r,V)$ is a subvariety and $x \in X$.
In the following, we give some naive definitions of the fact that $h$
is tangent to $X$ at $x$ in terms of the intersection of $X$, $h^\bot$
and $I_h$.

\begin{naiv}
``$h$ is tangent to $X$ at $x$ if $x \in I_h$
and the intersection $h^\bot \cap X$ is
singular at $x$.''
\end{naiv}
\noindent
This is really stupid, because if $x \in I_h$, then
$h^\bot$ is singular at $x$, and so is the intersection
$h^\bot \cap X$. So any $h$
such that $L_x \subset L_h$ will satisfy this condition, regardless to
the tangent space $T_xX$.

\begin{naiv}
``$h$ is tangent to $X$ at $x$ if $x \in h^\bot$
and the intersection $h^\bot \cap X$ is
singular at $x$.''
\end{naiv}
\noindent
For the same reason as above, it suffices that $L_h$ contains $L_x$ in
order that this condition holds. So if we define $X^*$ as the set of
$h$'s satisfying the above condition, we will not have a biduality
theorem. In fact, if for example $X=\{x\}$ is a point, then $X^*$ will
contain $\{h:L_h \supset L_x\}$ and ${(X^*)}^*$ will certainly not be
reduced to $\{x\}$.

\lpara

Therefore, it seems necessary to use the smooth subvariety $I_h$. In
this case, assuming that $I_h \cap X$ is singular is not quite
acurate, because $I_h$ has codimension larger than 1, so this
condition should be replaced by the fact that the intersection is not
transverse~:
\begin{naiv}
``$h$ is tangent to $X$ at $x$ if $x \in I_h$ and the
intersection $I_h \cap X$ is not transverse at $x$.''
\end{naiv}
\noindent
Again, if we take $X = \{x\}$, then 
$X^* = \{h : L_h \supset L_x \}$, and
${(X^*)}^* = \{y : \dim(L_x \cap L_y) > 0 \}$. So we don't have a
biduality theorem.

\lpara

Of course we could multiply such definitions; let us just consider one
more possibility~:
\begin{naiv}
``$X^*$ is the intersection of the usual dual variety of $X$ (in the
Plücker embedding) with $\G(r,V^*)$.''
\end{naiv}
\noindent
Already in case $r=2$ and $\dim V$ even, it is easy to see that
biduality will not hold. Let again $X = \{x\}$. The usual dual variety
of $X$ in the Plücker embedding is a hyperplane; therefore $X^*$ will
be a hyperplane section of $\G(2,V^*)$. As it is well-known,
the dual variety of
$\G(2,V^*) \subset \p \wedge^2 V^*$ is a hypersurface in 
$\p \wedge^2 V$. Therefore it follows that the usual dual variety of
$X^* \subset \p \wedge^2 V^*$ 
will have codimension at most 2 in $\p \wedge^2 V$. Thus
its intersection ${(X^*)}^*$ with $\G(2,V)$ cannot be a point.

\para

I hope that the previous unsuccessfull experiments will convince the
reader to accept a more conceptual definition of generalised dual
varieties.
Let $G$ be a semi-simple simply-connected complex algebraic group with
Lie algebra $\g$, and let $\g^*$ be the dual vector space of $\g$.
We fix $T\subset B \subset G$ a
maximal torus and a Borel subgroup of $G$.
If $P \subset G$ is a
parabolic subgroup, let $G/P$ denote the corresponding flag variety.
If $X$ is a variety, let $T^*X$ denote its cotangent bundle;
we denote $t_P:T^*G/P \rightarrow \g^*$ the natual map.

\begin{defi}
Let $P,Q$ be parabolic subgroups of $G$.
We say that $G/P$ and $G/Q$ allow duality if there is a nilpotent
orbit ${\cal O} \subset \g^*$ such that 
$t_P:T^*G/P \rightarrow \g^*$ and $t_Q:T^*G/Q \rightarrow \g^*$ are birational
isomorphisms between the cotangent bundles and $\cal O$.
\label{defi_dualite}
\end{defi}

Assume that $G/P$ and $G/Q$ allow duality. The birational map
$t_Q^{-1} \circ t_P : T^*G/P \dasharrow T^* G/Q$ will be denoted
$\mu$. Let $\co$ be such that $t_P(T^*G/P) = t_Q(T^*G/Q) = \overline \co$.
Let $X \subset G/P$ be any subvariety. Let $X^{sm}$ denote its smooth
locus and let 
$N^*X \subset T^* G/P$ denote the conormal bundle to $X^{sm}$~: we have
$(x,f) \in N^*X$ \iff $x \in X^{sm}$,
$f \in T^*_xG/P$, and $f_{|T_xX} = 0$.

\begin{defi}\ 
\label{defi_adapte}
\begin{itemize}
\item
A form $f \in T^*G/P$ (resp. $f\in T^*G/Q$) is called generic if it
belongs to $t_P^{-1}(\co)$ (resp. $t_Q^{-1}(\co)$).
\item
A subvariety $X \subset G/P$ 
is suitable if it is irreducible and there are generic forms
in $N^*X$.
\item
A point $x$ of a suitable variety $X$ is itself suitable if there
are generic forms in $N^*_xX$.
Let $X^s$ denotable the locus of suitable
points of $X$.
\end{itemize}
\end{defi}
\rek
One could also consider reducible suitable subvarieties : they
would be subvarieties such that every irreducible component is
suitable; we could then define the dual variety of a reducible suitable
variety as the union of dual varieties of its irreducible
components.

\begin{nota}
Let $\pi_P$ denote the projection $T^* G/P \rightarrow G/P$.
\end{nota}

\begin{defi}
If $X \subset G/P$ is suitable then we define $X^Q \subset G/Q$ as
the image of $N^*X$ by the rational map
$\pi_Q \circ \mu$.
\end{defi}

In the rest of the article, $P$ and $Q$ will denote parabolic
subgroups of a reductive simply-connected group $G$ allowing duality.
Moreover, we denote
$p := \pi_P \circ \mu\ : T^* G/Q \dasharrow G/P$,\ 
$q := \pi_Q \circ \mu^{-1}\ : T^* G/P \dasharrow G/Q$ 
the relevant rational maps. Finally, let
$\co \subset \g^*$ be the $G$-orbit which is dense in 
$t_P(T^* G/P) = t_Q(T^* G/Q)$.

\begin{defi}\ 
\label{i_h}
\begin{itemize}
\item
Let $x \in X \subset G/P$. We say that $h \in G/Q$ is tangent to $X$
at $x$ if $h \in q(N^*_x X)$.
\item
If $h \in G/Q$, let $I_h$ denote the Schubert variety of elements in
$G/P$ which are incident to $h$, in the sense that $x$ is incident to
$h$ if the intersection of the stabilisors of $x$ and $h$ (in $G$)
contain a Borel subgroup.
\item
As a corollary of Borel's conjugacy theorem, $I_h$ is homogeneous
under the stabilisor of $h$.
\end{itemize}
\end{defi}

\noindent
The notion of tangency will be studied in more details in subsection
\ref{section_tangency}. Here we only remark the following~:

\begin{fait}
\label{fait_tangent}
If $h$ is tangent to $X$ at $x$, then the intersection
$I_h \cap X$ is not transverse at $x$.
\end{fait}
\noindent
The proof of this fact is postponed to section
\ref{section_tangency}~: see lemma \ref{lemm_q}.
Note that the converse does not hold in general, contrary
to the case when $G/P = \p V$.

\subsectionplus{Fundamental cases}

\label{subsection_fondamental}

\begin{defi}
Let $P,Q \subset G$ allow duality. We say that $P,Q,G$ is a
fundamental case if one of the following hold :
\begin{itemize}
\item
$G=SL_n$, $P$ and $Q$ are the stabilisors of supplementary subspaces
of $\C^n$.
\item
$G=Spin_{4p+2}$, $P$ and $Q$ are the stabilisors of supplementary
(and so of different families) isotropic subspaces of $\C^{4p+2}$.
\item
$G$ is of type $E_6$, and, with Bourbaki's notations 
\cite[p.261]{bourbaki} $(P,Q)$ correspond either to the roots 
$(\alpha_1,\alpha_6)$ or $(\alpha_3,\alpha_5)$.
\end{itemize}
If this holds, $G/P$ and $G/Q$ are called fundamental homogeneous spaces.
\label{fondamental}
\end{defi}

By \cite[theorem 6.1]{namikawa}, these examples are all the examples of 
maximal parabolic subgroups
allowing duality. Recall that the corresponding rational map
$\mu : T^* G/P \dasharrow T^* G/Q$ is called a stratified Mukai flop.

Moreover, in all the other cases, the rational map
$\mu : T^* G/P \dasharrow T^* G/Q$ (and, as we will see in subsection 
\ref{subsection_sequence}, the dual varieties) may be described using 
only fundamental stratified Mukai flops~: let us recall
this construction \cite[theorem 6.1]{namikawa}. 
Assume $P,Q \subset G$ are parabolic subgroups
included in a common parabolic subgroup
$R$. Then we have fibrations
$$
\begin{array}{ccccc}
G/P & & & & G/Q\\
& f_P\searrow & & \swarrow f_Q\\
&& G/R
\end{array}
$$
with fibers $R/P$ and $R/Q$. Let $U(R)$ denote the unipotent radical of
$R$ and $Z(R)$ its connected center; let $L=R/Z(R)U(R)$; $R/U(R)$
is isomorphic with
a levi factor of $R$ and $L$ is semi-simple.
Moreover, $R/P$ and $R/Q$ are $L$-homogeneous varieties~:
denote $\pi : R \rightarrow L$ the projection, and denote
$P_L:=\pi(P)$ (resp. $Q_L:=\pi(Q)$) we have
$R/P \simeq L/P_L$ and $R/Q \simeq L/Q_L$. Assume now that $P_L,Q_L$
allow duality. Therefore there is a rational map
$\mu_L : T^* L/P_L \dasharrow T^* L/Q_L$ 
which can be used to define the stratified Mukai flop.

In fact, let $z \in G/R$, and denote $\cf_z:=f_P^{-1}(z)$
(resp. $\cg_z:=f_Q^{-1}(z)$), and let $i_z:\cf_z \rightarrow G/P$
(resp. $j_z:\cg_z \rightarrow G/Q$) be the natural inclusions. We have
$\cf_z \simeq L/P_L$ and $\cg_z \simeq L/Q_L$.
Let $L_z = R_z/Z(R_z)U(R_z)$ denote the group
isomorphic with $L$ acting on $\cf_z$ and $\cg_z$.
Because $\mu_L$ is canonical, it defines an algebraic family
of rational maps
$\mu_z : T^* \cf_z \dasharrow T^*\cg_z$ parametrised by $G/R$.
Now, if $\alpha$ is an element of $T^* G/P$,
say $\alpha \in T^*_x G/P$ with $x \in G/P$, then we can restrict this
linear form
to $T_x \cf_{f_P(x)}$; this gives an element in the bundle
$T^*\cf_{f_P(x)}$ which we denote $f_x$.
Finally, recall that $\pi_P:T^* G/P \rightarrow G/P$ and
$\pi_Q:T^* G/Q \rightarrow G/Q$ denote the bundle projections.
With these notations we have the following proposition~:

\begin{prop}
\label{prop_fibration}
If $f \in T_x^* G/P$ belongs to the open $G$-orbit, then
$f_x \in T^*\cf_{f_P(x)}$ belongs to the open $L_{f_P(x)}$-orbit, and
$\pi_Q(\mu(f)) = j_{f_P(x)} \circ \pi_{Q_L} \circ \mu_x(f_x)$.
\end{prop}
\noindent
Then, using \cite[theorem 4.1]{flop_scorza},
one can deduce a description of the flop
$T^* G/P \dasharrow T^* G/Q$.
\begin{proo}
Since both maps of the proposition are equivariant,
we can assume that $x$ corresponds to the base point in $G/P$.
If the restriction of $f$ to $T_x \cf_{f_P(x)}$ would belong to a
closed $L$-stable strict subvariety of $T^* \cf_{f_P(x)}$, then forms
in the
$G$-orbit of $f$ would restrict to non generic forms; therefore this
$G$-orbit could not be dense in $T^* G/P$.

Let $\u({\mathfrak r})$ and $\z({\mathfrak r})$ denote the
nilpotent part and the centraliser of the Lie algebra
${\mathfrak r}$ of $R$.
Let $\plie$ be the Lie algebra of $P$. Under $t_P$,
$f$ is mapped to an element in $\g^*$ which is orthogonal to $\plie$
and therefore to
$\u({\mathfrak r}) \oplus \z({\mathfrak r})$. It thus defines an
element $\overline f$ in ${\mathfrak l}^*$, if ${\mathfrak l}$ denotes
the Lie algebra of $L$. If $y \in L/Q_L$ denotes the element
$\mu_L(\overline f)$, then, by definition of the Mukai flop,
$\overline f$ is orthogonal to
$\overline \q_y$ 
($\overline \q_y$ denotes the Lie algebra of the stabilisor of $y$ in
$L$). Thus it follows that $f$ vanishes on $\q_{j(y)}$, the Lie
algebra of the stabiliser of $j(y)$ in $G/Q$. Therefore, $y$ equals
$\pi_Q \circ \mu(f)$.
\end{proo}

Now, \cite[theorem 6.1]{namikawa} states that for any pair $(P,Q)$
of parabolic
subgroups allowing duality, we can find a chain 
$(P_0=P,P_1,\ldots,P_n=Q)$ of parabolic subgroups such
that all the pairs $(P_i,P_{i+1})$ are as above and the corresponding
pair $P_L,Q_L \subset L$ is a fundamental case. Therefore, the description of
stratified Mukai flops in the fundamental cases is enough to understand
all stratified Mukai flops. As we will see in section \ref{section_reduction},
the same is true as far as dual varieties are concerned.

\subsectionplus{Recallections about fundamental homogeneous spaces}

\label{subsection_notation}

We now introduce some notations and recall some results
for fundamental homogeneous spaces
which will be used throughout the article.
In particular, we give in each case
a simple way of understanding the rational map
$q : T^* G/P \dasharrow G/Q$.

Let $r$ and $n$ be integers with $2r<n$.
The Grassmannian parametrising $r$-linear subspace of a fixed
$n$-dimensional vector space $V$ will be denoted
$\G(r,V)$. The dual Grassmannian, parametrising
codimension $r$ subspaces of $V$, will be denoted $\G(r,V^*)$.
Let $x \in \G(r,V)$; it represents a linear subspace of $V$ which will
be denoted $L_x$. Moreover, we have a natural identification
$T^*_x \G(r,V) \simeq Hom(V/L_x,L_x)$. If 
$\varphi \in Hom(V/L_x,L_x)$ is generic (that is, of rank $r$), then
its kernel is a codimension $r$ subspace of $V$ containing $L_x$. In
fact, we have $q(\varphi) = \ker \varphi$.

\lpara

Let $p$ be an integer. Let $V$ be a vector space of dimension $4p+2$,
equipped with a quadratic form. In case we need a basis for $V$, we
will take a hyperbolic one, of the form 
$(e_1^+,\ldots,e_{2p+1}^+,e_1^-,\ldots,e_{2p+1}^-)$, such that the
quadratic form is given by
$Q(\sum x_i^+e_i^+ + \sum x_i^-e_i^-) = \sum x_i^+x_i^-$.
Recall that the variety parametrising isotropic subspaces of $V$ of
dimension $2p+1$ has two connected components, which will be denoted
$G/P = \G_Q^+(2p+1,4p+2)$ and $G/Q = \G_Q^-(2p+1,4p+2)$. As in the
case of Grassmannians, for $x \in \G_Q^+(2p+1,4p+2)$ and
$h \in \G_Q^-(2p+1,4p+2)$, we denote $L_x,L_h$ the corresponding
isotropic subspaces. The relation $x \in I_h$ amounts to
$\dim(L_x \cap L_h) = 2p$. Given $x \in \G_Q^+(2p+1,4p+2)$ and
$L \subset L_x$ of dimension $2p$,
there is exactly one $h \in \G_Q^-(2p+1,4p+2)$ such that
$L_x \cap L_h = L$~: this yields a natural isomorphism between
$I_x$ and $\p L_x^*$.

The map $q$ may be defined as follows. Let $x \in \G_Q^+(2p+1,4p+2)$;
the cotangent space $T_x^* \G_Q^+(2p+1,4p+2)$ identifies with
$\wedge^2 L_x$. If $\omega \in \wedge^2 L_x$ is a skew form of rank
$2p$, let $L_\omega$ be its image. It is a hyperplane in $L_x$;
therefore it defines a unique element $h \in \G_Q^-(2p+1,4p+2)$
such that $L_x \cap L_h = L_\omega$. We have
$q(\omega) = h$.

\lpara

As far as the exceptional group $E_6$ is concerned, we denote $V_i$
the $i$-th fundamental representation of $E_6$, so that
$E_6/P_i \subset \p V_i$. We have $V_6 = V_1^*$ and $V_5 = V_3^*$. In
terms of this embedding, an element $h \in \p V_1^*$ belongs to
$E_6/P_6$ \iff it contains the linear span of two tangent spaces
$T_xE_6/P_1,T_yE_6/P_1$, for some distinct $x,y \in E_6/P_1$.

We refer to \cite{flop_scorza} for the
proofs of the following results. Let $x \in E_6/P_1$. The cotangent
space $T_x^* E_6/P_1$ identifies with $\oc \oplus \oc$, if $\oc$
denotes the algebra of complexified octonions, an 8-dimensional 
non-associative and non-commutative algebra over $\C$. This algebra is a
normed algebra~: there is a quadratic form $N:\oc \rightarrow \C$ such
that $N(z_1z_2) = N(z_1)N(z_2)$ for all $z_1,z_2 \in \oc$. The variety
$I_x$ is an 8-dimensional smooth quadric. It is convenient to denote
$Z = \C \oplus \oc \oplus \C$ a 10-dimensional space, equipped with
the quadratic form $Q(t,z,u) = tu - N(z)$. Then $I_x$ is the smooth
quadric defined by $Q$ and $q$ is defined by
$q((z_1,z_2)) = [N(z_1):z_1\overline z_2:N(z_2)] \in \p Z$
\cite[theorem 3.3 and corollary 3.2]{flop_scorza}.

\lpara

To visualise the homogneous space $E_6/P_3$ (resp. $E_6/P_5$), 
we use the fact
that its points parametrise
projective lines included in $E_6/P_1$ (resp. $E_6/P_6$)
\cite[theorem 4.3 p.82]{landsberg}.
To avoid confusions between points in $E_6/P_1$ and
$E_6/P_3$, we will denote the latters with greek letters.
Since the marked Dynkin diagrams of $E_6/P_3$ and $E_6/P_5$ are
respectively
$\dynkinep{9}{4}{2}$ and $\dynkinep{9}{4}{6}$, we see that for
$\kappa \in E_6/P_5$, $I_\kappa \simeq \G(2,5)$. Let us describe this
isomorphism $I_\kappa \simeq \G(2,5)$ more explicitly, since this will be
needed to describe the rational map $q$.
If $\alpha \in E_6/P_3$, we will denote $l_\alpha \subset E_6/P_1$ the
corresponding
line and $L_\alpha$ the linear subspace it represents. 
By \cite[proposition 3.6]{flop_scorza}, the span of 
the affine tangent spaces
$\widehat{T_x E_6/P_1}$ in $V_1$ for $x$ in
$l_\alpha$ is a 22-dimensional linear subspace in $V_1$ denoted 
$S_\alpha$.
Therefore, any 25-dimensional space which contains $S_\alpha$ defines a
pencil of hyperplanes which belong to $E_6/P_6 \subset \p V_1^*$, that
is, a point in $E_6/P_5$. Denoting $Q_\alpha = V_1/S_\alpha$, 
this shows that
$I_\alpha = \G(3,Q_\alpha) \simeq \G(2,5)$. Dually, for 
$\beta \in E_6/P_5$, $I_\beta \simeq \G(2,W_\beta)$, 
where $W_\beta \subset V_1$ is a
5-dimensional linear subspace such that 
$\p W_\beta \subset E_6/P_1$.

A Levi factor of $P$
contains $L' \simeq SL_2 \times SL_5$, and $L_\alpha$ (resp. $Q_\alpha$)
is the natural representation of $SL_2$ (resp. $SL_5$).
These representations are usefull describing $T^* E_6/P_3$~:
let $[e] \in E_6/P_3$ denote the base point; according to
\cite[propositions 3.6 and 3.7]{flop_scorza},
$T^*_{[e]} E_6/P_3$ is no longer an irreducible
$L'$-module, but there are exact sequences of $L'$-representations
\begin{equation}
\label{suite_exacte}
0 \rightarrow L_\alpha^* \otimes \wedge^2 Q_\alpha \rightarrow
T_{[e]} E_6/P_3 \rightarrow Q_\alpha^* \rightarrow 0
\end{equation}
$$
\label{suite_exacte_cotangent}
0 \rightarrow Q_\alpha \rightarrow
T^*_{[e]} E_6/P_3 \stackrel{\pi}{\rightarrow}
L_\alpha \otimes \wedge^2 Q_\alpha^*
\rightarrow 0.
$$

We now describe the rational map $q$.
Choose a base
$e_1^*,e_2^*$ (resp. $f_1,\ldots,f_5$) of $L_\alpha^*$ (resp. 
$Q_\alpha$). The
rational map $q:T_{[e]}^* E_6/P_3 \dasharrow I_{[e]}$ factors through
$L_\alpha \otimes \wedge^2 Q_\alpha^*$, and the induced rational map
$\overline q:L_\alpha \otimes \wedge^2 Q_\alpha^* 
\dasharrow I_{[e]} = \G(2,Q_\alpha^*)$ is described as
follows : let $\varphi \in L_\alpha \otimes \wedge^2 Q_\alpha^* \simeq
Hom(L_\alpha \otimes \wedge^2 Q_\alpha^*)$ be generic. Its image in
$\G(2,W_5^*)$ under $\overline q$
represents the linear subspace generated by
\begin{itemize}
\item
the
orthogonal for the alternate form $\varphi(e_2^*)$ of the kernel of
$\varphi(e_1^*)$,\hspace{.4cm} and
\item the orthogonal for $\varphi(e_1^*)$ of the kernel
of $\varphi(e_2^*)$.
\end{itemize}
This is well-defined \iff $\varphi(e_1^*)$ and
$\varphi(e_2^*)$ have maximal rank 4 and the corresponding orthogonals
are different lines in $V_2^*$. This is proved in
\cite[theorem 4.3]{flop_scorza}.

\subsectionplus{Dual schemes}

For some purposes (for example \cite{hermitian}),
it may be usefull to extend the above definition of
dual varieties to more general schemes. The goal of this subsection is 
to explain how this is possible.

Let us first define the cotangent scheme of a subscheme. So let $S$ be
an arbitrary scheme and $f:X \rightarrow Y$ a morphism above $S$.
The cotangent scheme $T^*X$ of $X$ is
${\bf Spec}\ S\point {\cal H}om(\Omega_{X/S},\co_X)$; it is a scheme
over $X$, equipped with a natural section, the zero section.
Now $f$ induces a natural morphism of sheaves
$f^*\Omega_{Y/S} \rightarrow \Omega_{X/S}$, and so a morphism
$f^* T^*Y \rightarrow T^*X$. We finally define the cotangent
scheme $N^*_{X,Y}$ as the fiber above the zero section of this map.

Let $G$ be a semi-simple Chevalley group scheme over $\Z$, $P$ and $Q$
parabolic subgroups.

\begin{defi}
$P$ and $Q$ allow duality if the complex groups $P(\C),Q(\C)$ do.
\end{defi}

If $P$ and $Q$ allow duality, although the moment map 
$T^*G/P \rightarrow \g^*$
may fail to be birational in positive caracteristic, there is still a
well-defined birational map $T^*G/P \dasharrow T^*G/Q$, defined over
$\Z$~:

\begin{prop}
There is a $G$-equivariant birational map
$\mu : T^*G/P \dasharrow T^*G/Q$ defined over $\Z$.
\end{prop}
\pr
By \cite[theorem 6.1]{namikawa} and proposition \ref{prop_fibration},
any pair of parabolic subgroups allowing duality is
related by a chain of pairs $(P,Q)$ of parabolic subgroups for which the
birational map $T^* G/P \dasharrow T^* G/Q$ is locally isomorphic with
a family of birational maps given by a
fundamental stratified Mukai flop. It is therefore
enough to check the proposition
for fundamental cases. In these cases it is a consequence of the explicit
description of this flop recalled in \ref{subsection_notation}.
\qed

If $S$ is a scheme and $G,P,Q$ are as above, let $G_S,P_S,Q_S$ the
groups obtained by base change $S \rightarrow Spec\ \Z$.

\begin{defi}
Let $S$ be a reduced irreducible scheme, and
let $f: X \rightarrow G_S/P_S$ be an irreducible
closed $S$-subscheme. We say that
$X$ is suitable if $\mu$ is defined at the generic point
of $N^*_{X,G_S/P_S}$.
In this case, the dual scheme of $X$ is the
scheme-theoretic image of $N^*_{X,G_S/P_S}$ under $\pi_Q \circ \mu$.
\end{defi}

\subsectionplus{Functorial property of dual varieties}

\label{subsection_fonctoriel}

We come back to our setting of complex geometry. In the usual setting,
if $X_1,X_2 \subset \p V$ are subvarieties, with $X_1 \subset X_2$,
there is in general no relation of inclusion between the dual
varieties of $X_1$ and $X_2$. Thus dual varieties have bad functorial
properties. The only thing one can say is the following obvious result.

\begin{prop}
Let $P,Q \subset G$ allow duality.
Let $X \subset G/P$ be suitable, $g \in G$, and $Y=g(X)$. Then $Y$ is
suitable and $g(X^Q)=Y^Q$.
\label{foncteur}
\end{prop}
\pr
Let $x \in X^s$ and $f \in N_x^*X \subset \g^*$ an element in the open
$G$-orbit. Then
$\tr g^{-1}.f \in N_{g(x)}^*Y$ is also in the open $G$-orbit. Therefore,
$Y$ is suitable. Moreover, since $q$ is equivariant, 
$q(\tr g^{-1}.f)=g.q(f)$.
Therefore, $g(X^Q) \subset Y^Q$. By symmetry, we have
also $g^{-1}(Y^Q) \subset X^Q$, so $g(X^Q) = Y^Q$.
\qed

\subsectionplus{Dual varieties in type A}

\label{subsection_grassmannienne}

In this section, I give a description of the dual variety of a subvariety
$X \subset \G(r,V)$
using an analog of the quotient map $V \dasharrow \p V$ for
Grassmannians.

\lpara

If $A$ and $B$ are vector spaces, $Inj(A,B)$ will
denote the sets of linear (resp. linear and injective) maps from $A$ to $B$.
Let $\varpi : Hom(\C^r,V) \dasharrow \p Hom(\C^r,V)$ denote the
natural rational map, and let
$\pi : \p Hom(\C^r,V) \dasharrow \G(r,V)$ map $\varphi$ of rank $r$ on
its image.
Dually, consider $\varpi' : Hom(V,\C^r) \dasharrow \p Hom(V,\C^r)$
and $\pi' : \p Hom(V,\C^r) \dasharrow \G(r,V^*)$ mapping $\varphi'$
of rank $r$ on its kernel.
If $X \subset \G(r,V)$ is a subvariety, let $\overline X^o$ denote
the set $\pi^{-1}(X)$ and $\overline X$ its closure in $\p Hom(\C^r,V)$.

\begin{prop}
Let $X \subset \G(r,V)$ be a suitable variety. Then
$X^Q = \pi' [{(\overline X)}^*]$, where ${(\overline X)}^*$ is the usual dual
variety of the subvariety $\overline X \subset \p Hom(\C^r,V)$
of a projective space.
\label{type_A}
\end{prop}
\pr
We fix a smooth point $x \in X$ and $\overline f \in \p Hom(\C^r,V)$ 
such that $\pi(\overline f)=x$,
and start with two easy lemmas.
\begin{lemm}
$\overline X$ is smooth at $\overline f$.
\end{lemm}
\begin{proo}
In a neighbourhood of $\overline f$ we have $\overline X = \overline X^o$.
Moreover,
the map $\pi:\overline X^o \rightarrow X$ is locally a trivial fibration
with fiber at $x$ the smooth variety $Inj(\C^r,L_x)$.
\end{proo}

\noindent
We denote $f \in Hom(\C^r,V)$ such that 
$\varpi(f) = \overline f$.
\begin{lemm}
The affine tangent space $\widehat{T_{\overline f}\overline X}$ 
is the linear space of
maps $g:\C^r \rightarrow V$ such that the composition
$L_x \stackrel{f^{-1}}{\rightarrow} \C^r 
\stackrel{g}{\rightarrow}V \rightarrow V/L_x$
belongs to $T_xX$.
\label{lemme_TxX}
\end{lemm}
\noindent
Recall that for $Z \subset \p W$ a projective variety and $z \in Z$,
the affine tangent space $\widehat{T_zZ} \subset W$ is the tangent
space of the affine cone over $Z$ at a lift of $z$ in $W$.
\begin{proo}
Let $\widehat X \subset Hom(\C^r,V)$ be the affine cone over
$\overline X$.
Since $\widehat X$ is smooth at $f$, any tangent
vector is the direction of a curve included in $\widehat X$.
Let $\gamma : (C,0) \rightarrow (\widehat X,f)$ be a curve in
$\widehat X$ and let $g=\gamma'(0) \in Hom(\C^r,V)$. 
Under the well-known identification of
$T_x \G(r,V)$ with $Hom(L_x,V/L_x)$, the composition of the lemma
equals 
$(\pi \circ \varpi \circ \gamma)'(0)$. Therefore it belongs to
$T_xX$. By dimension count, the lemma follows.
\end{proo}

\noindent
\underline{\bf Proof of proposition \ref{type_A} :}
The linear subspace $(T_f\widehat X)^\bot \subset Hom(V,\C^r)$ is
the set of $g$'s such that for all $h \in T_f\widehat X$, the
composition 
$\C^r \stackrel{h}{\rightarrow} V \stackrel{g}{\rightarrow} \C^r$ is
traceless. Since $T_f \widehat X$ contains $Hom(\C^r,L_x)$, this means
that $g$ is induced by a morphism $\overline g: V/L_x \rightarrow \C^r$ such
that $f \circ \overline g$ is orthogonal to 
$T_xX \subset Hom(L_x,V/L_x)$, by lemma \ref{lemme_TxX}. Therefore,
for $h \in \G(r,V^*)$, we have 
$h \in \pi' \circ \varpi((T_f \widehat X)^\bot)$
\iff
$h \in q(N_x^*X)$.
\qed




\sectionplus{Reduction to fundamental examples}

\label{section_reduction}

From section \ref{section_definition}, we see that there are a lot of pairs of
parabolic subgroups which allow duality. In this section, I will show
that to understand all the dual varieties, it is enough to
understand dual varieties for fundamental cases. 

For example, the
varieties corresponding to the marked diagrams
$$
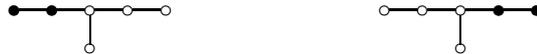
\dynkineDeuxMarque{11.5}{4}{0}{2}
\hspace{2cm}
\dynkineDeuxMarque{11.5}{4}{6}{8}
$$
both have dimension 26. Using tables in \cite[p.202]{mcgovern}, we see
that there is a unique nilpotent
orbit of dimension 52 in ${\mathfrak e}_6$ and that the disconnected
centralizer of an element of this orbit is trivial. Therefore, the two
corresponding parabolic subgroups $P,Q \subset E_6$
allow duality. It may seem at first that the corresponding
duality
$X \subset G/P \mapsto X^Q \subset G/Q$ has to do with the
exceptional geometry of $E_6$. However, we will see that it is not
the case; indeed, $X^Q$ can be described using dual varieties in four
classical homogeneous spaces. Indeed, \cite[theorem 6.1]{namikawa}
is verified
in this case thanks to the sequence of parabolic subgroups
$$
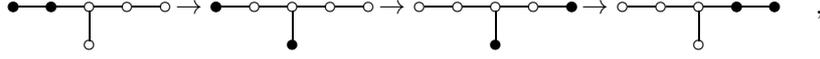
\dynkineDeuxMarque{8.5}{4}{0}{2} \rightarrow
\dynkineUneMarqueP{8.5}{4}{0} \rightarrow
\dynkineUneMarqueP{8.5}{4}{8} \rightarrow
\dynkineDeuxMarque{8.5}{4}{6}{8} \hspace{.5cm},
$$
and we will see in
this section how to compute accordingly dual varieties.
We will show that the computation of the dual variety
$X^Q$ for $X \subset G/P$ can be
done in three steps, the first and the last in a family of spinor varieties
$\G_Q^+(5,10)$, and the second in a family $\p^5$'s.

\subsectionplus{Biduality theorem}

Let $G$ be as above, $P,Q,R \subset G$ be subgroups 
such that $P$ and $Q$ allow
duality, and $Q$ and $R$ allow duality. 
Then, by definition $P$ and $R$ also allow duality.

\begin{theo}[Biduality theorem]
Let $X \subset G/P$ be an suitable variety. Then $X^Q$ is suitable
and $\mu(N^*X) = N^*X^Q$. In particular, 
${(X^Q)}^R = X^R$.
\label{theo_bidualite}
\end{theo}
\noindent
If $G=SL_n$, $P=R$ is the stabilisor of a line and $Q$ is the
stabilisor of a hyperplane, we recover the
usual biduality theorem.
\begin{proo}
We follow the argument of \cite[pp.27 to 30]{gkz}.\\
Let $N=\mu(N^*X) \subset T^* G/Q$.
Recall that $T^* G/Q$ is a symplectic variety. Moreover, it is proved
in \cite{gkz} that $N^*X \subset T^* G/P$ is a lagrangien
subvariety of $T^* G/P$. Let $\co \subset \g^*$ denote the 
nilpotent orbit which
closure is the image of $T^* G/P$. Since the birational morphisms
$T^* G/P \stackrel{\sim}{\dasharrow} {\cal O}
\stackrel{\sim}{\dasharrow} T^* G/Q$ are symplectic, it follows that
$N$ is also lagrangien.

Moreover, it has the property that if 
$(x,f) \in N$ and $\lambda \in \C$, then $(x,\lambda f) \in N$. 
This follows from the fact that the image of $N^*X$ in
$\overline {\cal O}$ is stable under multiplication by scalars.
From \cite[proposition 3.1]{gkz}, we know that
$N=N^*Z$ for $Z=\pi_Q(N)=X^Q$. Therefore, $\mu(N^*X)=N^*X^Q$
and $X^Q$ is suitable.

Since ${(X^Q)}^R$ (resp. $X^R$) is the image of $N^*X_Q$
(resp. $\mu(N^*X)$) under the
rational map $T^* G/Q \dasharrow G/R$, these varieties are equal.
\end{proo}

The following corollary shows that the name of biduality theorem for
the above result is justified~:

\begin{coro}
\label{coro_bidualite}
Let $P,Q \subset G$ allow duality. If $X \subset G/P$ is suitable,
then $X^Q$ is suitable and ${(X^Q)}^P = X$. Moreover,
if $x \in X$ and $h \in X^Q$, then $h$ is
tangent to $X$ at $x$ \iff $x$ is tangent to $X^Q$ at $h$.
\end{coro}
\begin{proo}
To prove that ${(X^Q)}^P = X$, 
it is enough to take $R=P$ in theorem \ref{theo_bidualite}, after
observing that for suitable $X \subset G/P$, $X^P = X$.
The second result, that $h$ is tangent to $X$
at $x$ \iff $x$ is tangent to $X^Q$ at $h$ follows from the fact the
first (resp. the second) affirmation means that $(x,h)$ lies in the
image by $(p,\pi_Q)$ of an element in $\mu(N^*X)$ (resp. $N^*X^Q$).
\end{proo}

\subsectionplus{Families of dual varieties}

\label{subsection_sequence}

We consider the following situation :
let $R \subset G$ be a parabolic subgroup.
Let $P,Q \subset R \subset G$ be parabolic subgroups and recall
notations of subsection \ref{subsection_fondamental}.
If $X \subset G/P$ and $z \in G/R$, denote $X_z:=X \cap \cf_z$.
Assume $P_L,Q_L \subset L$ allow duality.
For $z \in G/R$ and suitable $Y \subset \cf_z \simeq L/P_L$, let 
$Y^{Q_L} \subset \cg_z \simeq L/Q_L$ denote its generalised dual
variety.

\begin{theo}
With the previous notations, assume that $P,Q \subset G$ allow
duality, and also $P_L,Q_L \subset L$.
If $X \subset G/P$ is suitable, then
for generic $x \in X$, $X_{f_P(x)} \subset \cf_{f_P(x)}$ is 
suitable. Moreover, $X^Q$ is the closure of the union of
the $X_{f_P(x)}^{Q_L}$, for such $x$ in $X$.
\label{sous-groupe}
\end{theo}
\pr
Let $f \in N^*_xX$ an element which $G$-orbit in $T^* G/P$ is dense
and set $z = f_P(x)$.
We have seen in the proof of proposition \ref{prop_fibration}
that the restriction $f_x$ of $f$ to 
$T_x \cf_z$ is a generic element in $T^* \cf_z \simeq T^* L/P_L$.
Moreover, this restriction belongs to $N^*_x X_z$, so that $X_z$ is
suitable.

Let $q_z : T^* \cf_z \dasharrow \cg_z$ be the composition of
$\mu_z : T^* \cf_z \dasharrow T^* \cg_z$ and the projection
$T^* \cg_z \rightarrow \cg_z$. Proposition \ref{prop_fibration}
states that $q(f) = j_z \circ q_z(f_x) \in G/Q$. Therefore it follows
that $q(N^*X_{|X_z}) = j_z(X_z^{Q_L})$.
The description of $X^Q$ in the theorem follows.
\qed

\lpara

As a consequence of theorems \ref{theo_bidualite} and \ref{sous-groupe},
if $P=P_1 \times P_2$ and $Q=Q_1 \times Q_2$ are parabolic subgroups
of $G=G_1 \times G_2$, 
and if $X=X_1 \times X_2$, then we have
$X^Q = X_1^{Q_1} \times X_2^{Q_2}$.




\sectionplus{Tangency for fundamental examples}

\label{section_tangency}

In this section, if $x \in X \subset G/P$,
I introduce a definition of the embedded tangent cone at $x$,
$\overline{T_xX}$, which is a subvariety of $G/P$ and a cone
at $x$ (in a suitable
sense). I also introduce the cotangent
variety at $x$, $\overline{N_xX}$, which is
a subvariety of $G/Q$. Moreover a notion of ``linear varieties'' is defined
and linear varieties are classified.

From now on, $P,Q \subset G$ are fundamental subgroups of $G$ allowing
duality.

\subsectionplus{A tangent element is incident}

In this subsection, we prove that if $x \in X \subset G/P$ and 
$h \in G/Q$ is tangent to $X$ at $x$ 
(see definition \ref{i_h}), then $h$ is incident to $x$
(in the sense that the stabilisors of $x$ and $h$ contain a common
Borel subgroup). This only holds in fundamental cases.

\begin{nota}
Let $x \in \g$ nilpotent. Then there exists $y,h \in \g$ such that $(x,y,h)$
is a $\mathfrak {sl}_2$-triple. For $i \in \Z$, let $\g_i$ denote
$\{X \in \g : [h,X]=iX\}$. The parabolic subalgebra 
$\plie_x := \oplus _{i \geq 0} \g_i$ does not depend on $y$ and $h$
\cite[theorem 3.8]{mcgovern}, and is
called the canonical parabolic subalgebra of $x$.
\end{nota}

In the following lemma, I say that $\plie \subset \g$ is a maximal
parabolic subalgebra of $\g$ of fundamental type if the pair
$(\g,\plie)$ is the pair of Lie algebras of groups $(G,P)$ as in
definition \ref{fondamental}. 

\begin{lemm}
Let $x \in \g$ and $\plie$ be a polarisation of $x$. Assume that
$\plie$ is a maximal parabolic subalgebra of fundamental type. 
Then $\plie_x \subset \plie$.
\end{lemm}
\begin{proo}
Let $\plie$ be a maximal parabolic subalgebra which is a polarisation of $x$.
Let $(x,y,h)$ be a $\mathfrak {sl}_2$-triplet, $\h$ a Cartan
subalgebra containing $h$ and $\Delta=\{\alpha_1,\ldots,\alpha_r\}$ 
a basis of the root system such
that $\forall \alpha \in \Delta,\alpha(h) \geq 0$. 

We denote $\plie_1$ the following maximal parabolic subalgebra~:
$$
\plie_1 := \h \oplus \bigoplus_{
\begin{array}{c}
\alpha = \sum_j k_j \alpha_j\\
k_i \geq 0
\end{array}
} \g_{\alpha}\ \ ,
$$
where $i$ is chosen such that $\plie$ is conjugated to $\plie_1$
(such $i$ exists because $\plie$ is a maximal parabolic subalgebra).

Let us prove that $x \in \u(\plie_1)$.
According to the
decomposition $\g = \h \oplus \bigoplus_{\alpha} \g_{\alpha}$, we can
write $x = h_x + \sum_\alpha x_\alpha$, with $h_x \in \h$ and
$x_\alpha \in \g_\alpha$. Now, since $[h,x]=2x$, we deduce that 
$h_x = 0$ and that for any root $\alpha$, either $x_\alpha = 0$ or
$\alpha(h) = 2$. 

\begin{clai}
If $\alpha = \sum k_j \alpha_j$ is a root, then
$\alpha(h)=2 \Longrightarrow k_i > 0$.
\end{clai}
\begin{proo}
This is proved by ad hoc arguments in all cases. Assume first that 
$\g = \liesl_n$ and that $\plie$ is the stabilisor of an $r$-dimensional
subspace. Thus $i=r$.
Recall that the weighted diagram of $x$ is by definition the list
of the values $\alpha_j(h)$. The weighted diagrams of nilpotent
elements in $\liesl_n$ are well-known; in our case, since
$x$ is a generic element of $\u(\plie)$ with $\plie$ conjugated to
$\plie_1$, we have 
$\alpha_i(h)=\alpha_{n-i}(h)=1$ and the other values $\alpha_j(h)$ equal 0.
The equality $\alpha(h)=2$ with $\alpha = \sum k_j \alpha_j$ amounts
to $k_i + k_{n-i} = 2$, which implies $k_i = k_{n-i} = 1$.

Assume now that $\g = \spin_{4p+2}$. In this case, there is only one 
possibility for the $G$-orbit in $\spin_{4p+2}$
of $x$, and $\alpha_j(h)=1$ \iff $\alpha_j$
is a spin root (ie $j \in \{2p,2p+1\}$); otherwise $\alpha_j(h)=0$. Therefore
$\alpha(h) = 2$ implies that $\alpha$ is not less than the root
$\alpha_{2p-1} + \alpha_{2p} + \alpha_{2p+1}$, which implies the claim.

If $\g$ is of type $\e_6$ and $\plie$ corresponds to the first root, then
the weighted diagram of $x$ is $\poidsesix 100010$
(see \cite[table p.202]{mcgovern}). Since for all roots
$\sum k_j \alpha_j$ we have $-1\leq k_1,k_6\leq 1$, we again have
$\sum k_j \alpha_j(h) = 2 \Rightarrow k_1=k_6=1$. In case $\plie$
corresponds to the second root, the weighted diagram is
$\poidsesix 010100$. The equality $\alpha(h)=2$ for
$\alpha = \sum k_j \alpha_j$ implies that $k_3 + k_5 = 2$. If
$k_3 = 2$, then necessarily $k_5 \geq 1$ (see the list of roots in
\cite{bourbaki}), so we get a contradiction. Similarly $k_5 \leq 1$.
So $k_3 = k_5 = 1$, and again the claim is proved.
\end{proo}

This claim therefore proves that if $x_\alpha \not = 0$, with
$\alpha = \sum k_j \alpha_j$, then $k_i > 0$. 
This proves that $x$ belongs
to 
$$
\bigoplus_{
\begin{array}{c}
\alpha = \sum_j k_j \alpha_j\\
k_i \geq 1
\end{array}
} \g_{\alpha}\ ,
$$
which is readily seen to be $\plie_1^\bot = \u(\plie_1)$.
Thus $x \in \u(\plie_1)$
and $\plie_1$ is a polarisation of $x$.
Now, since the map $T^* G/P \rightarrow \g$ is birational on its
image, there is a unique polarisation of $x$ in the conjugacy class of
$\plie$. Therefore $\plie = \plie_1$.

\lpara

Let us now show that $\plie_x \subset \plie_1$. Since $\plie_x \supset \h$,
it is the sum of $\h$ and some root spaces. Now, assume 
$\g_\alpha \subset \plie_x$, with $\alpha = \sum k_j \alpha_j$.
This means that $\sum k_j \alpha_j(h) \geq 0$.
I claim that $k_i \geq 0$. In fact, if $k_i < 0$, then $\alpha$ is a negative
root, so $k_j \leq 0$ for all $j$. We therefore have
$\sum k_j \alpha_j(h) \leq k_i \alpha_i(h)$.
In the proof of the above claim, we have seen that we
allways have $\alpha_i(h)=1$. So we get a contradiction.

Therefore, we have $k_i \geq 0$, and so $\g_\alpha \subset \plie_1$.
Since $\plie_1 = \plie$, we have proved that $\plie_x \subset \plie$, as
claimed.
\end{proo}

\begin{coro}
\label{coro_incident}
If $\plie$ and $\q$ are polarisations of the same nilpotent element
$x$, and are maximal parabolic subalgebras of
fundamental type, then they contain a
common Borel subalgebra.
\end{coro}
\pr
They both contain the canonical parabolic subalgebra $\plie_x$.
\qed

\lpara

We now show, with an example, that the above corollary is wrong if one
considers non maximal parabolic subalgebras.

\begin{exem}
Let $\g = {\mathfrak sl}_n$. Let $x \in \g$ be an element of rank 2,
such that $x^3 = 0$ but $x^2 \not = 0$. Let $\plie$ (resp. $\q$) be the
parabolic subalgebra preserving the image of $x^2$ and the image of
$x$ (resp. the kernel of $x$ and the kernel of $x^2$). Then we have 
$x \in  \u(\plie)$ and $x \in \u(\q)$. However, since
$\im\ x \not \subset \ker x$, $\plie$ and $\q$ are
not incident.
\end{exem}
\begin{proo}
If $y \in \plie$ (resp. $y \in \q$), then the commutator
$[x,y]$ is strictly upper triangular for the filtration
$\im\ x^2 \subset \im\ x \subset \C^n$ (resp.
$\ker x \subset \ker x^2 \subset \C^n$). Therefore, $[x,y]$ is
traceless and so $x \in \plie^\bot$ (resp. $x \in \q^\bot$).
\end{proo}

\lpara

The Schubert varieties $I_h$ (recall definition \ref{i_h})
give a geometric understanding of the
rational map $q : T^* G/P \dasharrow G/Q$ : 
\begin{lemm}
\label{lemm_q}
Assume $P$ and $Q$ are maximal parabolic subgroups.
Let $x \in G/P$ and $h \in G/Q$, and
let $f$ be a generic element in 
$T_x^* G/P$. Then $q(f) = h$ \iff $x \in I_h$ and the cotangent
form $f$ vanishes on $T_x I_h$.
\end{lemm}
\noindent As a consequence of the lemma, there is a unique $h$ such that
$x \in I_h$ and $f$ vanishes on $T_x I_h$. By definition, if $h$ is
tangent to $X$ at $x$, then there exists $f \in N^*_xX$ such that $q$
is defined at $f$ and $q(f)=h$. Thus the lemma implies that the
intersection $I_h \cap X$ is not transverse at $x$, as was stated in 
fact \ref{fait_tangent}.
\begin{proo}
Let $x \in G/P$; $t_P$ restricts to an isomorphism between
$T^*_x G/P$ and $(\g/\plie_x)^* \subset \g^*$ if $\plie_x$ denotes the Lie
algebra of the stabilisor of $x$. Conversely, given $\eta \in {\cal O}$,
$\pi_P(t_P^{-1}(\eta))$ is the unique $x \in G/P$ such that the corresponding
parabolic subalgebra $\plie_x$ is orthogonal to $\eta$.

Let $x \in G/P$, $f \in T_x^* G/P$ generic and 
$\eta = t_P(f) \in \plie_x^\bot$, and let $h=q(f)$. 
The previous argument shows that $h$ is the unique element in $G/Q$
such that $\eta$ vanishes on $\q_h$. Moreover, we know by
corollary \ref{coro_incident} that $x \in I_h$.
Note that $T_xG/P = \g / \plie$ and
$T_x I_h \simeq \q_h/(\plie_x \cap \q_h)$.
Since $\eta$ vanishes on $\plie_x$, it will
vanish on $\q_h$ \iff it vanishes on $\q_h/(\plie_x \cap \q_h)$, namely,
\iff $f$ vanishes on $T_xI_h$.
\end{proo}

\begin{exem}
Let $h \in G/Q$ and let $X = I_h \subset G/P$. Then $X$ is suitable
and $X^Q = \{h\}$. Moreover $p(T^*_hG/Q) = I_h$.
\end{exem}
\begin{proo}
First, let $x \in X$, let $f \in T^*_xX$ be generic and let $h=q(f)$.
By corollary \ref{coro_incident}, $x$ and $h$ are incident, and by
lemma \ref{lemm_q}, $f$ vanishes on $T_x I_h$. Thus, $I_h$ is
suitable. Since $G/Q$ is homogeneous, $I_h$ is suitable for all 
$h \in G/Q$.

Let $x \in X$ and $f \in N^*_xX$ generic. Then by the above $q(f)=:h'$
is well-defined, and by lemma \ref{lemm_q} again, $h'$ is the unique
element in $G/Q$ such that $x \in I_{h'}$ and such that $f$ vanishes on
$T_x I_{h'}$. Since $h$ satisfies these conditions, $h'=h$. Therefore,
$X^Q = \{h\}$.

For the last point, we note that $p(T^*_h G/Q) = \{h\}^P = I_h$, by
biduality theorem \ref{theo_bidualite}, 
since we have proved that $I_h^Q = \{h\}$.
\end{proo}

\subsectionplus{Dual varieties and cones}

If $X \subset \p V$ is included in a hyperplane represented by
$h \in \p V^*$, then the dual variety of $X$, which is a subvariety of
$\p V^*$, is a cone over $h$. The aim of this subsection is to prove
an analogous result for our generalised dual varieties. Our first goal
is to define cones.

\begin{defi}Let $x_1,x_2 \in G/P$ 
\begin{itemize}
\item
$x_1,x_2$ are linked if there exists $h \in G/Q$ such that 
$x_1,x_2 \in I_h$.
\item
If $E \subset G/P$, let 
$\displaystyle I_E := \bigcap_{x \in E} I_x \subset G/Q$.
\item
If $x_1,x_2$ are linked, denote 
$\displaystyle L(x_1,x_2) = \bigcap_{h \in I_{\{x_1,x_2\}}} I_h$.
\end{itemize}
\end{defi}
\noindent
In $\p V$, all points are linked, and $L(x_1,x_2)$ is the line through
$x_1$ and $x_2$. The difference between $\p V$ and our general
situation is that in general $G$ does not 
act transitively on pairs of distinct points $x,y \in G/P$, so that $L(x,y)$
may depend, up to isomorphism, on $x$ and $y$.
However, cones are defined in perfect analogy~:

\begin{defi}
\label{def_cone}
Let $X \subset G/P$ and $x \in X$. Then $X$ is a cone over $x$ if
for all $y \in X$, $x$ and $y$ are linked and $L(x,y) \subset X$.
\end{defi}
\noindent
An equivalent definition is that for generic $y \in X$ the same
condition holds, as will be clear from the following description of
$L(x,y)$~:

\begin{prop}
\label{prop_l}
Let $x \not = y \in G/P$. We have~:
\begin{itemize}
\item
If $G/P = \G(r,V)$, then $(x,y)$ are linked \iff 
$\codim_V(L_x+L_y) \geq r$, in which case $L(x,y) = \G(r,L_x+L_y)$.
\item
If $G/P = \G_Q^+(2p+1,4p+2)$, then $(x,y)$ are linked \iff
we have $\dim(L_x \cap L_y)=2p-1$, in which case
$L(x,y) = \{ z : L_z \supset L_x \cap L_y \} \simeq \p^1$.
\item
If $G/P = E_6/P_1$, then $(x,y)$ are allways linked. In case a line
passes through $x$ and $y$ in $E_6/P_1$, then $L(x,y)$ is this line;
otherwise, there is a unique smooth 8-dimensional quadric through $x$
and $y$, and $L(x,y)$ is this quadric.
\item
If $G/P = E_6/P_3$, then $(x,y)$ are linked \iff there is a $\G(2,5)$
through them. If $\dim(L_x \cap L_y) = 1$ then 
$L(x,y)$ is equal to $\G(2,L_x + L_y) \simeq \p^2$, 
otherwise $L(x,y) \simeq \G(2,5)$.
\end{itemize}
\end{prop}
\noindent
In this proposition, for the two exceptional cases, I use the minimal
projective homogeneous embedding $E_6/P_i \subset \p V_i$. 
For example, in the case of
$E_6/P_3$, the condition that there is a $\G(2,5)$ through $x$ and $y$
means that there is a linear 10-dimensional subspace $W \subset V_3$ 
containing $x$ and $y$ and such that $\p W \cap E_6/P_3$ is projectively
isomorphic with a Grassmanian $\G(2,5)$ in its Plücker
embedding. Recall also that $E_6/P_3$ parametrises projective lines in
$V_1$ which are included in $E_6/P_1$. For $x \in E_6/P_3$, the
corresponding 2-dimensional subspace of $V_1$ has been denoted $L_x$.
\begin{proo}
The first case follows directly from the definition. In the second
case, one only has to note that if there exists $h \in G/Q$ such that
$(x,h),(y,h)$ are incident, then 
$\dim(L_x\cap L_h) = \dim(L_y \cap L_h) = 2p$, so
$\dim(L_x\cap L_h\cap L_y) = 2p-1$.

\lpara

For the exceptional cases one oviously has to use the geometry of the
involved homogeneous spaces. Let us first consider $E_6/P_1$. For all
$x \in E_6/P_1$, $I_x$ is a smooth 8-dimensional quadric. Moreover, for any
$x\not = y \in E_6/P_1$, the intersection of the two quadrics
$I_x$ and $I_y$ is either a point or a $\p^4$. In
fact, this was proved in 
\cite[propositions IV.3.2 and IV.3.3]{zak} in the context of Severi
varieties, but also follows easily from the fact that there are three
$E_6$-orbits in $E_6/P_1 \times E_6/P_1$ \cite[proposition 18]{quantique}.
Given $x,y \in E_6/P_1$, we can have $x=y$, $x \not =
y$ and there is a line through $x$ and $y$, or there is no line
through $x$ and $y$. This describes the three orbits in 
$E_6/P_1 \times E_6/P_1$.
In the degenerate case when a line passes through $x$ and $y$, 
$I_x \cap I_y$ is thus isomorphic with $\p^4$. Dually,
the intersection of all the
$I_h$ for $h$ in this $\p^4$ is a linear space (indeed, $x \in I_h$
\iff $x \in E_6/P_1 \subset \p V_1$ is orthogonal to
$\widehat{T_h E_6/P_6} \subset V_6=V_1^*$) and contains $x$ and $y$;
a direct computation of dimension shows that it is exactly
the line through $x$ and $y$.
In the generic case, $I_x \cap I_y = \{h\}$; therefore 
$L(x,y) = I_h$ is the unique 8-dimensional quadric through $x$ and $y$.

\lpara

Let $\alpha,\beta \in E_6/P_3$ be linked, and
denote $\kappa \in E_6/P_5$ an element such that 
$\alpha,\beta \in I_\kappa$. According to subsection
\ref{subsection_fondamental}, $\alpha$ and 
$\beta$ represent 2-dimensional
subspaces of a 5-dimensional subspace of $V_1$ denoted
$W_\kappa$; we have denoted $L_\alpha,L_\beta$ these spaces.

Assume first that $\dim(L_\alpha \cap L_\beta) = 1$.
It is proved in
\cite[proposition 3.6]{flop_scorza} that the linear
span of all the affine tangent spaces at the points of the projective plane
generated by $L_\alpha$ and $L_\beta$
is 24-dimensional and equal to the span of affine tangent spaces at
points in $l_\alpha \cup l_\beta$. 
Thus $(\alpha,\beta)$ defines a projective plane in $E_6/P_6$ and also in
$E_6/P_5$. Moreover
$I_{\alpha,\beta} = I_{\G(2,L_\alpha+L_\beta)} \simeq \p^2$ 
and $L(\alpha,\beta) = \G(2,L_\alpha+L_\beta) \simeq \p^2$.

Assume finally that $L_\alpha$ and $L_\beta$ don't meet. 
Let $L \subset L_\alpha \oplus L_\beta$ be any 3-dimensional subspace;
the linear span $S_L$ of
the affine tangent spaces at points of $\p L$ is again 24-dimensional,
and any element in $I_{\alpha,\beta}$ must contain it.
Assume that $I_{\alpha,\beta}$
contains two points $\kappa,\lambda \in G/Q$. These points would correspond to
codimension 2 subspaces $L_\kappa,L_\lambda$ of 
$V_1$ containing $S_L$; therefore
$L_\kappa$ and $L_\lambda$ would be contained 
in a common hyperplane of $V_1$.
Since $\alpha,\beta \in I_{\kappa,\lambda}$, by
the case considered above, this would in turn imply that $L_\alpha$ and
$L_\beta$ meet in dimension 1, which we have excluded.
Therefore we have proved that $I_{\alpha,\beta} = \{ \kappa \}$, so
$L(\alpha,\beta) = I_\kappa$ is isomorphic with $\G(2,5)$.
\end{proo}

\begin{theo}
\label{theo_cone}
Let $h \in G/Q$ and let $X \subset G/P$ such that $X \subset I_h$. 
Then $X$ is suitable and $X^Q$ is a cone over $h$.
\end{theo}
\begin{rema}
In fact, as the proof will show, in all the cases but in type $A$,
a stronger result holds~: for any $k \in X^Q$, there is a certain
homogeneous subvariety
$q(\C.f + N^*_x I_h) \subset G/Q$, of type given by lemmas
\ref{lemm_cone_gr}, \ref{lemm_cone_spin}, \ref{lemm_cone_e6_1},
and \ref{lemm_cone_e6_2}, containing (eventually stricly) $L(h,k)$,
and included in $X^Q$. Although the idea of proof of this theorem
is uniform, this proof unfortunately ends up with a case by case
analysis.
\end{rema}
\begin{proo}
If $x \in X$, then $N^*_xX$ contains $N^*_x I_h$ on which $q$ is
well-defined generically, so $X$ is suitable.
Assume $X \subset I_h$ and let $k$ be a generic element in $X^Q$.
By definition of $X^Q$ there is an element $x \in X$ and 
$f \in N^*_xX$ such that $k = q(f)$. Since $x \in X \subset I_h$, we
have $h \in I_x$. By corollary \ref{coro_incident}, $k \in I_x$;
therefore $h$ and $k$ are linked. Moreover, we have 
$f \not \in N_x^* I_h$ (otherwise we would have $q(f)=h$).
Therefore it follows from
the inclusion $q(\C.f + N^*_x I_h) \subset X^Q$ and the following
lemmas \ref{lemm_cone_gr}, \ref{lemm_cone_spin}, \ref{lemm_cone_e6_1}
and \ref{lemm_cone_e6_2} that $L(h,k) \subset X^Q$.
\end{proo}

\begin{lemm}
\label{lemm_cone_gr}
Let $x \in \G(r,V),\ h\not = k \in \G(r,V^*)$ such that $h,k \in I_x$. Let
$f \in N_x^* I_h$ such that $q$ is defined at $f$. Then
$q(\C.f + N_x^* I_k) = L(h,k)$.
\end{lemm}
\begin{proo}
Let $(e_i)$ be a base of $V$ and $(e_i^*)$ the dual base.
Up to the action of $SL(V)$, we may assume that
$L_x$ is the span of $e_1,\ldots,e_r$,
$L_h$ is the span of $e_{n-r+1}^*,\ldots,e_n^*$, $L_k$ that
of $e_{n-r+1-l}^*,\ldots,e_{n-r}^*,e_{n-r+1}^*,\ldots,e_{n-l}^*$, and
finally that $f \in N_x^* I_h \simeq Hom(L_x^*,L_h)$ is defined by
$f(e_j^*) = e_{n-r+j}^*$. Since
$N_x^* I_k = Hom(L_x^*,L_k)$, a straightforward computation proves
the lemma.
\end{proo}

\begin{lemm}
\label{lemm_cone_spin}
Let $x \in \G_Q^+(2p+1,4p+2),\ h\not = k \in \G_Q^-(2p+1,4p+2)$ such that
$h,k \in I_x$. Let
$f \in N_x^* I_h$ such that $q$ is defined at $f$.
Then $q(\C.f + N_x^* I_k) \simeq \p^{2p-1}$.
\end{lemm}
\noindent
This lemma implies theorem \ref{theo_cone} in this case since
$q(\C.f + N_x^* I_k)$ is a linear space containing $h$ and $k$, 
and will therefore contain the line through $h$ and $k$.
\begin{proo}
We may assume that $x$ represents the isotropic subspace
$L_x$ generated by $e_1^+,\ldots,e_{2p+1}^+$. Since $L_k$ meets $L_x$ along 
a hyperplane, we may further assume that this hyperplane is generated by
$e_2^+,\ldots,e_{2p+1}^+$. We therefore have
$N^*_x I_k = \wedge^2 \scal{e_2^+,\ldots,e_{2p+1}^+} \subset
\wedge^2 L_x = T^*_x \G_Q^+(2p+1,4p+2)$. Let 
$f \in T^*_x \G_Q^+(2p+1,4p+2)$; since 
$f \not \in N^*_x I_k$ (otherwise we would have $h=k$),
the class of $f$ modulo $N^*_x I_h$ is the same as that of some form
$e_1^+ \wedge e$, with $e \in \scal{e_2^+,\ldots,e_{2p+1}^+}$, and we may
assume that $e=e_2^+$.

Recall that $I_x \simeq \p L_x^*$~:
I claim that $q(\C.f + N_x^* I_k)$ is the orthogonal of $e_2^+$ in
$\p L_x^*$. In fact, let $\wedge^p (\C.f + N_x^* I_k) \subset 
\wedge^{2p} L_x \simeq L_x^*$ be the linear span of all the forms 
in $\wedge^{2p} L_x$ which can
be written as a wedge product of $p$ forms in 
$\C.f \oplus N_x^* I_k$. We have
$\wedge^p (\C.f + N_x^* I_k) \subset {(e_2^+)}^\bot$; therefore
$q(\C.f + N_x^* I_k) \subset \p {(e_2^+)}^\bot$.

On the other hand, let 
$\displaystyle f_0=\Sigma_{i=1}^p e_{2i-1}^+ \wedge e_{2i}^+$; we have
$\displaystyle f_0^{\wedge(p-1)}=\Sigma_{i=1}^p {e_{2i-1}^+}^* 
\wedge {e_{2i}^+}^* \wedge e_{2p+1}^*$, from which is follows that the 
rational map $\C.f + N_x^* I_k \dasharrow \p {(e_2^+)}^\bot,
g \mapsto [g^{\wedge p}]$
is submersive at $f_0$, which implies the claim and the lemma. 
\end{proo}

\begin{lemm}
\label{lemm_cone_e6_1}
Let $x \in E_6/P_1,\ h \not = k \in E_6/P_6$ such that
$h,k \in I_x$. Let
$f \in N_x^* I_h$ such that $q$ is defined at $f$.
If there passes a line through $h$ and $k$ in $E_6/P_6$
then $q(\C.f + N_x^* I_k) \simeq \p^4$, otherwise
$q(\C.f + N_x^* I_k) = I_x$
\end{lemm}
\begin{proo}
We adopt the same strategy of proof as for lemma \ref{lemm_cone_spin}.
Let $x \in E_6/P_1$ be fixed.
In subsection \ref{subsection_fondamental}, we saw that
$T^*_xX$ identifies with $\oc \oplus \oc$,
$I_x$ with the projective quadric in $\p (\C \oplus \oc \oplus \C)$
defined by $tu-N(z)=0$.
We can assume that $k \in I_x$ is the class of $(0,0,1)$.
Therefore, $N_x^* I_k = \overline{q^{-1}(k)} = \{(0,z) : z \in \oc\}$.

Write $f=(z_0,z_1)$. Since 
$q((z_0,z_1)) = [N(z_0) : z_0 \overline z_1 : N(z_1)]$, 
there will be a line through
$q(f)$ and $k$ in the quadric $I_x$ \iff $N(z_0) = 0$. If this occurs, then
$$q(\C.f + N_x^* I_k) = \{[(0,u,t)]:t\in \C,u\in L(z_0) \},$$
where $L(z_0)$ denotes the set of right multiples of $z_0$~:
$L(z_0)=\{z_0z:z\in \oc\}$.
It is a linear subspace of $\oc$ of dimension 4, so
$q(\C.f + N_x^* I_k)$ is isomorphic with $\p^4$, as desired. If
$N(z_0) \not = 0$, then left multiplication by $z_0$ is invertible, so that
$q : \C.f + N_x^* I_k \dasharrow I_x$ is dominant, and the lemma again
holds.
\end{proo}

\begin{lemm}
\label{lemm_cone_e6_2}
Let $\alpha \in E_6/P_2,\ \kappa,\lambda \in E_6/P_5$ such that
$\kappa,\lambda \in I_\alpha$. Let
$f \in N_\alpha^* I_\kappa$ such that $q$ is defined at $f$.
Then $q(\C.f + N_\alpha^* I_\kappa) = I_\alpha$.
\end{lemm}
\begin{proo}
We fix $\alpha \in E_6/P_2$. Let $f_1^*,\ldots,f_5^*$ be a 
base of $Q_\alpha^*$ and assume 
that $\kappa$ corresponds to the linear subspace
generated by $f_4^*,f_5^*$. Recall that there is a natural
surjective map $\pi:T_\alpha^* E_6/P_2 
\rightarrow Hom(L_\alpha^*,\wedge^2 Q_\alpha^*)$. Moreover,
$\pi(N^*_\alpha I_\kappa) = Hom(L_\alpha^*,L) 
\subset Hom(L_\alpha^*,\wedge^2 Q_\alpha^*)$, where 
$L \subset \wedge^2 Q_\alpha^*$ is
generated by $f_1^* \wedge f_4^*,f_1^* \wedge f_5^*,
f_2^* \wedge f_4^*,f_2^* \wedge f_5^*,f_3^* \wedge f_4^*,f_3^* \wedge f_5^*,
f_4^* \wedge f_5^*$ (for example, this follows from the fact that for any
$\varphi \in Hom(L_\alpha^*,L)$,
$\overline q(\varphi)$, if defined, equals $\kappa$).

Let $M \subset \wedge^2 Q_\alpha^*$ be generated by
$f_1^* \wedge f_2^*,f_1^* \wedge f_3^*,f_2^* \wedge f_3^*$, so that
$L \oplus M = \wedge^2 Q_\alpha^*$; the class of $\pi(f)$ modulo
$\pi(N_\alpha^* I_\kappa)$ is the class of a unique
$\overline f \in Hom(L_\alpha^*,M)$.
Assume first that the rank of $\overline f$ is 1. We can therefore
assume that $\overline f(e_1^*) = f_1^* \wedge f_2^*$ and
$\overline f(e_2^*) = 0$, where $e_1^*,e_2^*$
is a suitable basis of $L_\alpha^*$.

In the array below we give, for $\omega \in \wedge^2 Q_\alpha^*$, the value
of the derivative $d\overline q_{\overline f}(\varphi)$, for 
$\varphi:L_\alpha^* \rightarrow \wedge^2 Q_\alpha^*$ given by
$\varphi(e_1^*) = \omega$ and $\varphi(e_2^*) = 0$~:
$$
\begin{array}{llll}
f_1^* \wedge f_2^* \mapsto 0 & f_1^* \wedge f_3^* \mapsto 0 &
f_1^* \wedge f_4^* \mapsto 0 & \hspace{-1cm}f_1^* \wedge f_5^* \mapsto 0 \\
f_2^* \wedge f_3^* \mapsto 0 & f_2^* \wedge f_4^* \mapsto f_3^* \wedge f_4^* &
f_2^* \wedge f_5^* \mapsto f_1^* \wedge f_4^* + f_3^* \wedge f_5^*\\
f_3^* \wedge f_4^* \mapsto 0 & f_3^* \wedge f_5^* \mapsto 0 &
f_4^* \wedge f_5^* \mapsto f_1^* \wedge f_5^*.
\end{array}
$$
The following gives similar values for $\varphi$ defined by
$\varphi(e_1^*) = 0$ and $\varphi(e_2^*) = \omega$~:
$$
\begin{array}{llll}
f_1^* \wedge f_2^* \mapsto 0 & f_1^* \wedge f_3^* \mapsto 0 &
f_1^* \wedge f_4^* \mapsto 0 & \hspace{-1cm}f_1^* \wedge f_5^* \mapsto 0 \\
f_2^* \wedge f_3^* \mapsto 0 & f_2^* \wedge f_4^* \mapsto f_2^* \wedge f_3^* &
f_2^* \wedge f_5^* \mapsto f_3^* \wedge f_4^* + f_1^* \wedge f_2^*\\
f_3^* \wedge f_4^* \mapsto 0 & f_3^* \wedge f_5^* \mapsto 0 &
f_4^* \wedge f_5^* \mapsto f_1^* \wedge f_4^*.
\end{array}
$$
It follows from these computations that 
$\overline q(\C.\overline f + \pi(N_\alpha^* I_\kappa))$ has
dimension at least 6, so 
$\overline q(\C.\overline f + \pi(N_\alpha^* I_\kappa)) = I_\alpha$ 
in this case.

In case $\overline f$ has rank 2, the dimension of
$\overline q(\C.\overline f + \pi(N_\alpha^* I_\kappa))$ 
will not vary if $\overline f$ is replaced by
$g.\overline f$, where $g \in SL(L_\alpha) \times SL(Q_\alpha)$ 
preserves $\kappa$.
Using a $\C^*$-action we can degenerate 
$\overline f \in Hom(L_\alpha^*,M)$ 
to some element $\overline f_0$ of
rank one, for which we have already seen that 
$\dim \overline q(\C.f_0 + \pi(N_\alpha^* I_\kappa)) = 6$.
Since this dimension is
lower semi-continuous, we have 
$\dim \overline q(\C.\overline f + \pi(N_\alpha^* I_\kappa)) = 6$ and
the lemma is proved.
\end{proo}

\subsectionplus{The cotangent space and the tangent cone of a variety}

\label{subsection_cotangent}

\begin{defi} Let $x \in X$ be suitable.
\begin{itemize}
\label{cotangent}
\item
The embedded cotangent space of $X$ at $x$ is 
$\overline{N_xX}:=q(N^*_xX)\subset G/Q$.
\item
The embedded tangent space of $X$ at $x$ is
$\overline{T_xX}={\overline{N_xX}}^P$.
\item
$X \subset G/P$ is a linear subvariety if $\overline{T_xX}$ does not
depend on $x$ suitable in $X$.
\end{itemize}
\end{defi}

\rks
\begin{itemize}
\item
The notion of (co)-tangent space (and therefore of linear varieties)
of $X \subset G/P$ could be defined for non maximal parabolic $P$, but
then it would depend on the choice of a parabolic subgroup $Q$.
\item
An equivalent definition of linear subvarieties
is that $\overline{N_xX}$ does not depend on
suitable $x$ in $X$, since $\overline{N_xX} = {\overline{T_xX}}^Q$.
\item
By definition, $X^Q = \overline {\cup_{x \in X^s} \overline{N_xX}}$.
\item
In projective spaces, the tangent cone is the usual embedded tangent
space and linear varieties are linear subspaces. Linear varieties will
be classified in the next subsection.
\end{itemize}

\begin{exem}
Let $x \in G/P$ and $X=\{x\}$. Then $\overline{T_xX} = \{x\}$.
\end{exem}
\begin{proo}
In fact, $\overline{N_xX} = q(T^*_xG/P) = X^Q$, so this follows from
theorem \ref{theo_bidualite}.
\end{proo}

\begin{lemm}
For $x \in X^s$,
$\overline{T_xX}$ is a cone over $x$ and therefore $x \in \overline{T_xX}$.
\end{lemm}
\pr
In fact, for $x \in X^s$, we have $\overline{N_xX} \subset I_x$, so 
$\overline{T_xX} = {\overline{N_xX}}^P$ is a cone
over $x$ by theorem
\ref{theo_cone}.
\qed

\subsectionplus{Linear subvarieties}

In this subsection, we classify linear subvarieties.

\begin{prop}
\label{prop_lineaire}
The following array gives the list of all linear subvarieties~:
$$
\begin{array}{|c|c|c|}
\hline
\mbox{G/P} & \mbox{Linear varieties} \\
\hline
\G(r,n) & \G(r,p),\ r\leq p \leq n  \\
\hline
\G^+_Q(2p+1,4p+2) & \{pt\};I_h,h \in \G_Q^-(2p+1,4p+2) \\
\hline
E_6/P_1 & \{pt\};I_h,h\in E_6/P_6\\
\hline
E_6/P_2 & \{pt\};I_\kappa,\kappa\in E_6/P_5\\
\hline
\end{array}
$$
\end{prop}
\begin{proo}
Let $X \subset G/P$ be linear.
First, we prove that
$\forall x \in X,\overline{T_xX} = X$, and that $X^Q$
is linear.
Let $x \in X^s$. Then
$X^Q = \overline {\cup_{y \in X^s} \overline{N_yX}} = \overline {N_xX}$,
since for all $y \in X^s$, $\overline{N_yX} = \overline{N_xX}$.
Therefore, $X = {\overline{N_xX}}^P = \overline{T_xX}$ by
corollary \ref{coro_bidualite}.
Let $h \in X^Q$ and $x \in X$. Then, by biduality theorem again,
$x \in \overline{N_hX^Q}$ \iff $h \in \overline {N_xX} = X^Q$.
Therefore, $\overline{N_hX^Q}  = X$ and $X^Q$ is linear and the claim
is proved.
Since $X = \overline{T_xX}$ for all $x \in X$, $X$ is a cone over all
of its points by theorem \ref{theo_cone}.
\lpara

We finish the proof case by case. 
In the case of Grassmannians, if we denote 
$W = \sum_{x \in X} L_x$, since $X$ is a cone over all of its points,
we have $\G(r,W) \subset X$, and so $X = \G(r,W)$.

In the case of spinor varieties, any $x,y \in X$ must be linked, which
implies that the line through $x$ and $y$ is in $X$, so $X$ is a
linear subspace. As a consequence of the following proposition
\ref{prop_spinoriel}, the only linear subspaces which dual
variety is again a linear
subspace are the point and maximal linear subspaces. Since we have seen
that $X^Q$ must be a linear variety, the proposition follows in this
case.

Let $X \subset E_6/P_1$ be linear. Let $h \in X^Q$. If there are two
points $x,y \in X$ such that there is no line through $x$ and $y$,
then by lemma \ref{lemm_cone_e6_1} $L(x,y) = I_h$. Since $X \subset I_h$
and $L(x,y) \subset X$, we have $X = I_h$ (and $X^Q$ is a point).

Otherwise, by theorem \ref{theo_cone}, $X$ is a linear
subspace. If $X^Q$ is not a linear subspace, by the argument above, 
$X = {(X^Q)}^P$ is a point. Assume now that both $X$ and $X^Q$ are
linear subspaces, not reduced to a point. By lemma 
\ref{lemm_cone_e6_1}, $X^Q$ and $X={(X^Q)}^P$ contain a $\p^4$. But
this implies that $X^Q \subset I_X$ is at most 1-dimensional
(see the proof of theorem \ref{theo_cone}), and we get
a contradiction.

\lpara

Let finally $X \subset E_6/P_3$ be linear. Assume $X$ is not reduced to
a point. Let $h \in X^Q$; we have $X \subset I_h$.
On the other hand, since
$X$ is not a point, by lemma \ref{lemm_cone_e6_2}, it must
contain $I_h$. Therefore, $X = I_h$.
\end{proo}




\sectionplus{Examples of dual varieties}

\label{section_exemple}

\subsectionplus{Dual varieties of isotropic Grassmannians}

Let $V$ be a vector space, $B:V \rightarrow V^*$ a bilinear form.
If $\epsilon = \pm 1$ and
$\tr B = \epsilon B$, we say that $B$ is $\epsilon$-symmetric.
Assume that this is the case. Let $r$ be an integer; we 
consider the variety $\G_B(r,V)$ of isotropic subspaces of $V$ of dimension
$r$. The aim of this subsection is to describe the dual of
$\G_B(r,V)$ in $G(r,V^*)$ in case $2r<\dim\ V$ (the other cases would be
similar). Note that we don't assume that $B$ is an isomorphism.

We have a rational map $\G_B(r,V) \dasharrow \G(r,V^*)$ which maps a linear
subspace to its orthogonal, and which is well-defined at the point $\alpha$
\iff $L_\alpha$ does not meet the kernel of $B$. Assuming there are
such points, we call co-isotropic Grassmannian the image of this
rational map.

\begin{prop}
\label{prop_symmetric}
Assume $\epsilon=1$. Then $\G_B(r,V)$ is suitable \iff and
$r \leq \rank(B)$. In this case, the dual variety of the isotropic Grassmannian
$G_B(r,V)$ is the co-isotropic Grassmannian.
\end{prop}

\begin{prop}
\label{prop_skew}
Assume $\epsilon=-1$. Then $\G_B(r,V)$ is suitable \iff $r$ is even
and $r \leq \rank(B)$. In
this case, the dual variety of the isotropic Grassmannian
$G_B(r,V)$ is the co-isotropic Grassmannian.
\end{prop}

\begin{proo}
We prove propositions \ref{prop_symmetric} and \ref{prop_skew}
simultaneously. Let $x \in \G_B(r,V)$ be generic. Under the natural
isomorphism $T_x \G(r,V) \simeq Hom(L_x,V/L_x)$, we have
the inclusion
$T_x \G_B(r,V) \supset Hom(L_x,L_x^\bot/L_x)$,
where $L_x^\bot$ denotes the orthogonal of $L_x$ with respect to
$B$. It follows that if $\codim\ L_x^\bot < r$, then 
$N_x^* \G_B(r,V)$ does not meet the open orbit in $T_x^* \G_B(r,V)$.
If $r > \rank(B)$, this occurs for all $x \in \G_B(r,V)$, hence
$\G_B(r,V)$ is not suitable.

Assume $r \leq \rank(B)$.
Now, let $x \in \G_B(r,V)$ such that $\codim\ L_x^\bot = r$.
Denote $Q_x = V/L_x$; we have a morphism 
$Q_x \rightarrow L_x^*$,
induced by $B$. Clearly, $T_x \G_B(r,V) \subset Hom(T_x,Q_x)$ 
is the
subspace of linear maps such that the composition 
$L_x \rightarrow Q_x \rightarrow L_x^*$
is $(-\epsilon)$-symmetric. Therefore,
the normal space of $\G_B(r,V)$ at $x$ identifies with
$\epsilon$-symmetric maps $L_x^* \rightarrow L_x$. 
Since $\G_B(r,V)$ will be suitable \iff there are such maps of rank $r$,
this occurs in all cases if $\epsilon = 1$ and exactly when $r$ is even
when $\epsilon = -1$.

Now, the computation of the dual variety is straightforward : since
we have already remarked that 
$T_x \G_B(r,V) \supset Hom(L_x,L_x^\bot/L_x)$,
the image of a generic conormal form at $x$ under the rational map
$q:T^*_x \G(r,V) \dasharrow G(r,V^*)$ is the element in
$\G(r,V^*)$ corresponding to $L_x^\bot$.
\end{proo}

\subsectionplus{Schubert varieties and quivers in the fundamental case}

\label{subsection_carquois}

In this subsection, I recall that to a cominuscule homogeneous space one
can naturally associate a quiver, such that Schubert cells are parametrised
by some subquivers. I also recall the Hasse diagram of a representation,
and show how the quiver of a cominuscule homogeneous space can be
identified with the Hasse diagram of a tangent space.
This identification is
due to Nicolas Perrin and
Laurent Manivel.
Then, I show that this identification behaves well as far as Schubert 
subvarieties are concerned. Finally, I extend these results to
$E_6/P_3$, which is not a cominuscule homogeneous space.

\lpara

The quiver of a cominuscule homogeneous space has been first introduced by
N. Perrin \cite[definition 3.2]{carquois};
here we use the slightly different definition
\cite[definition 2.1]{quantique}. Recall that
$\G(r,V),\G_Q^+(2p+1,4p+2)$ and $E_6/P_1$ are cominuscule spaces (in
fact even minuscule). The quiver is defined using a reduced expression
of $w_{G/P}$, the shortest element in the class of $w_0$ in $W/W_P$
($w_0$ is the longest element in $W$).
Choose a reduced expression
$w_{G/P} = s_{\beta_1} \cdots s_{\beta_N}$, with $N = \dim G/P$; the
vertices of the quiver 
$Q_{G/P}$ are in bijection with $[1,\ldots,N]$, and we refer
to \cite[definition 2.1]{quantique} for the definition of the arrows.
The quivers may be illustrated by relevant
examples as follows~: 
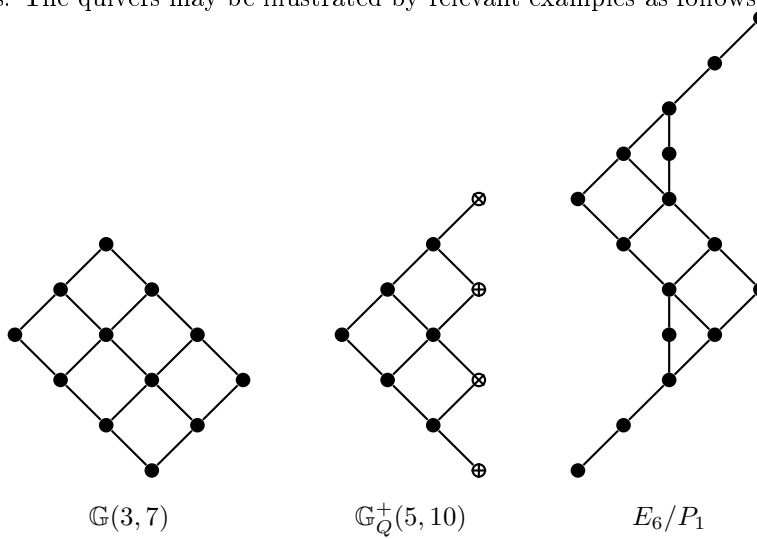
\vspace{-.8cm}

$$
\begin{array}{ccccc}
\input{g37}  &&  \input{gq5}  &&  \input{e6}  \\
\G(3,7) && \G_Q^+(5,10) && E_6/P_1
\end{array}
$$
In these pictures, all arrows are going down. Moreover, we will use
the definition of height of a vertex of such a quiver. More or less by
definition (see \cite[definition 4.7]{carquois}), it is the height of
the vertex in the above drawing, where by convention the lowest vertex
has height 1 (so the highest vertex has height respectively 6,7,11 for
$\G(3,7),\G_Q^+(5,10),E_6/P_1$).
\lpara

Later we will have to identify this quiver with a Hasse diagram. 
Let $V$ be a representation of a semi-simple group $\Lambda$.
Let us recall that the Hasse diagram
of $V$ is a quiver defined as follows.
The vertices of this quiver are the weights of $V$, and
there is an arrow from $\lambda_1$ to $\lambda_2$ if and only if
$\lambda_2 - \lambda_1$ is a simple root. For example, the Hasse diagram
of the $8$-dimensional representation of $Spin_8$ is given on the left~:

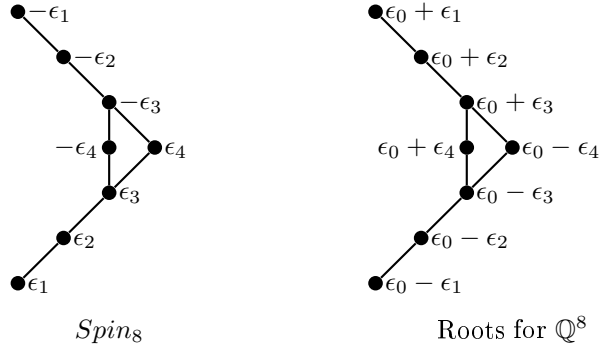
$$
\begin{array}{ccc}
\input{hasseq8} && \input{racineq10} \\
Spin_8 & \hspace{1cm} & \mbox{Roots for }\Q^8
\end{array}
$$

\begin{prop}
\label{hasse}
Let $G/P$ be cominuscule and let $x \in G/P$ be the base point. 
Let $\Lambda$ be a Levi factor
of the stabilisor of $x$. Then the quiver $Q_{G/P}$
of $G/P$ is isomorphic with the
Hasse diagram $H_{G/P}$ of the $\Lambda$-module 
$\widehat{T_xX}/L_x$.
\end{prop}
\noindent
If $G/P \subset V$, recall that $\widehat{T_xG/P} \subset V$ 
is the affine tangent space at $x$; it contains the
line $L_x \subset V$ represented by $x \in \p V$,
so that it makes sense to consider the quotient
$\widehat{T_xG/P}/L_x$.
We have stated this result without proof in
\cite[proposition 7]{quantique}. In this
article
I need the explicit isomorphism, this is why I sketch the proof,
leaving details to the reader.
\begin{proo}
It is known that to each vertex of the quiver one can associate a root of
$G$. In fact, choose a reduced expression
$w_{G/P} = s_{\beta_1} \ldots s_{\beta_N}$ and set
$$ \alpha_i = s_{\beta_N} \circ \ldots \circ s_{\beta_{i+1}} (\beta_i). $$
Since two different reduced expressions for
$w_{G/P}$ only differ by commutation relations,
it is easy to check that the induced map from the set of vertices of
the quiver
to the set of roots is well-defined (it does not depend on the reduced
expression). In the following, we consider that a reduced expression is
chosen, thus identifying this set of vertices with $[1,N]$.

For example, if $G/P$ is a smooth 8-dimensional quadric,
then its quiver, and the corresponding roots, are given above
(here we have shitfted the indices,
denoting $(\epsilon_0,\ldots,\epsilon_4)$ a basis of the weight lattice of
$Spin_{10}$). Note that the highest weight of the corresponding
$Spin_{10}$-representation is $\epsilon_0$, 
and that we recover the Hasse diagram
of $Spin_8$ by considering the weights $\epsilon_0 - \alpha_i$.

By \cite[proposition 4.9]{carquois}, we may reduce the proof of our
proposition to the particular case of a quadric
of any dimension, as above, because if
there is an arrow $i \rightarrow j$ in the quiver of $G/P$, then $i$
and $j$ belong to a subquiver of $Q_{G/P}$ isomorphic with the quiver
of a quadric. It is also possible (and probably shorter)
to check directly in each case that if $\omega$ denotes the highest
weight of $\Gamma(G/P,\co(1))$, then the set 
$\{\omega - \alpha_i : 1\leq i \leq N\}$ is exactly the set of weights
of the tangent space at the base point of $G/P$, and that the
bijection $i \mapsto \omega - \alpha_i$ is an
isomorphism of quivers $Q_{G/P} \rightarrow H_{G/P}$.
\end{proo}

\lpara

Given $[w] \in W/W_P$, we associate the Schubert subvariety 
$C_{[w]} \subset G/P$ which is the $B$-orbit closure of 
$[w] \in G/P$. Assuming that $w$ is the minimal length representative
of its class, we choose a reduced decomposition of $w$, and this
defines a subquiver $Q_w$ of the quiver $Q_{G/P}$
which is an order ideal
(this means that if $i \rightarrow j$ is an arrow in $Q_{G/P}$ and $i
\in Q_w$, then $j \in Q_w$~: 
see \cite[proposition 4.5]{carquois}). We can also consider the
subset $H_w$ of $H_{G/P}$ which elements are the weights of
$w^{-1}.T_{[w]}C_{[w]} \subset T_{[e]} G/P$. The following proposition
will be useful to compute the dual variety of $C_{[w]}$, because it
describes the tangent bundle of $C_{[w]}$~:

\begin{prop}
\label{prop_tangent_cw}
Under the isomorphism $Q_{G/P} \simeq H_{G/P}$
of proposition \ref{hasse}, we have
$Q_w = H_w$.
\end{prop}
\begin{proo}
Recall that $\omega$ denotes the highest weight of 
$\Gamma(G/P,\co(1))$. All the weights of 
$\widehat{T_x G/P}/L_x$ are of the
form $\omega + \alpha$, where $\alpha$ are all the roots not in
$\plie = Lie(P)$ (therefore $\alpha$ is a negative root). 
All the weights of $T_{[w]} C_{[w]}$ are of the form
$w.\omega + \beta$, with $\beta$ a positive root. Therefore, if
$\omega + \alpha$ is a weight of $w^{-1}.T_{[w]}C_{[w]}$,
$w.\alpha$ must be a positive root. So $\alpha$ must be a negative
root sent by $w$ to a positive root. Denote $l(w)$ the length of $w$;
there are $l(w)$ such roots,
namely $\{-\alpha_i:1\leq i \leq l(w)\}$. Since $l(w)$ is also the
dimension of $C_{[w]}$, it follows that
$H_w$ is exactly the set of weights of the form 
$\omega - \alpha_i,1\leq i \leq l(w)$, so the proposition follows.
\end{proo}

\lpara

We now consider the case of $E_6/P_3$. Let 
$[w] \in W/W_3$; we want to define a quiver $Q_{E_6/P_3}$ and a 
subquiver $Q_w$ which pictures the tangent bundle of $C_{[w]}$.
Since $E_6/P_3$ is not cominuscule, the quiver defined as in
\cite[definition 3.2]{carquois} is not well-defined
(it depends on a reduced expression of $w_{E_6/P_3}$), and
as we have already seen, the cotangent bundle $T^* E_6/P_3$ is no
longer irreducible, so its Hasse diagram is not suitable neither.

But our luck is that for $f \in T^*_\alpha E_6/P_3$, $q(f)$ only
depends on $\pi(f) \in L_\alpha \otimes \wedge^2 Q_\alpha^*$;
therefore, what we care for is not really the conormal bundle of
$C_{[w]}$, but rather its projection to the bundle
$L \otimes \wedge^2 Q^*$. This is why we consider the following~:

\begin{defi}\ 
\begin{itemize}
\item
Let $[e] \in E_6/P_3$ be the base point, and let $\Lambda$ denote a Levi
factor of $P_3$.
\item
Let $Q_{E_6/P_3}$ denote the Hasse diagram of the $\Lambda$-module
$L_{[e]}^* \otimes \wedge^2 Q_{[e]} \subset T_{[e]} E_6/P_3$.
\item
For $[w] \in W/W_3$ with $w$ the minimal length representative,
let $Q_w \subset Q_{E_6/P_3}$ denote the set of
weights of 
$w^{-1}.T_{[w]} C_{[w]}\ \cap\ 
(L_{[e]} \otimes \wedge^2 Q_{[e]}^*)$.
\end{itemize}
\end{defi}

\begin{prop}
\label{prop_order_ideal}
For $[w] \in W/W_3$, $Q_w \subset Q_{E_6/P_3}$ is an order ideal.
\end{prop}
\begin{proo}
For $a \in \{-2,-1,0,1,2\}$, let $\g_k \subset \e_6$ denote
$
\displaystyle
\bigoplus_{
\begin{array}{c}
\alpha = \sum_j k_j \alpha_j\\
k_3 = a
\end{array}
} \g_{\alpha}
$
(by this I mean that the Cartan subalgebra is included in $\g_0$).
We have $\g = \g_{-2} \oplus \g_{-1} \oplus \g_0 \oplus \g_1
\oplus \g_2$ and $Lie(P_3) = \g_0 \oplus \g_1
\oplus \g_2$. The tangent space $T_{[e]} E_6/P_3$ decomposes as
$\g_{-2} \oplus \g_{-1}$; let
${\cal P}$ denote the weights of 
$\widehat{T_{[e]} E_6/P_3}/L_{[e]}$ which are of the
form $\omega + \alpha$, with $\alpha$ a root of $\g_{-1}$ (${\cal P}$
is also the set of weights of $L_{[e]} \otimes \wedge^2 Q_{[e]}^*$).
I claim that $w$ induces an increasing bijection between
${\cal P}$ and its image.
The proposition follows from this claim because 
$Q_{[w]}$ is the set of weights of 
$w^{-1}.T_{[w]} C_{[w]}$ which are in ${\cal P}$; arguing as in the
proof of proposition \ref{prop_tangent_cw}, this is the set of roots
of $\g_{-1}$ which are mapped to a positive root by $w$, and this is
obviously an order ideal since $w$ is increasing.

To prove the claim, we note that
$L_{[e]} \otimes \wedge^2 Q_{[e]}^*$ is a minuscule 
$\Lambda$-representation, since $\Lambda$ contains 
$SL_2 \times SL_5$. Therefore $W_P$ permutes transitively the roots in 
$\g_{-1}$. Let $\alpha_0$ be the highest root of $\g_{-1}$, let
$\alpha_1,\alpha_2$ be roots of $\g_{-1}$ and assume
$\alpha_1 \leq \alpha_2$. We can find $w_1,w_2 \in W_P$ such that
$\alpha_i = w_i . \alpha_0$, assume moreover that $w_1,w_2$ are
minimal such elements. Since $\alpha_1 \leq \alpha_2$, we
have $w_1 \geq w_2$ for the Bruhat order, so that we may assume that a
$w_2$ is a product of reflexions appearing in a reduced expression of
$w_1$. Since $w$ is a minimal length representative in $W/W_3$, the
product $w.w_1$ is still a reduced expression, so 
$w.w_1 \geq w.w_2$, and so $w.\alpha_1 \leq w.\alpha_2$, as claimed.
\end{proo}

\begin{rema}
The same proof works for any $G/P$, as soon as $\g_{-1}$ is a
minuscule $\Lambda$-representation, with the notations of the proof.
\end{rema}


\subsectionplus{Schubert varieties and dual varieties}

\label{subsection_combinatoire}

The usual dual variety of a linear subspace is again a linear 
subspace. The
goal of this section is to generalise this result for Schubert varieties.

\begin{prop}
\label{prop_schubert}
Let $X \subset G/P$ be a suitable
Schubert variety. Then $X^Q \subset G/Q$
is a Schubert variety.
\end{prop}
\pr
In fact, $X^Q$ is a $B$-stable (proposition \ref{foncteur})
irreducible closed subvariety of $G/Q$.
\qed

\lpara

Recall that $B$-stable
Schubert varieties in $G/P$ are parametrised by the quotient set
$W/W_P$. For $[w] \in W/W_P$ (resp. $[x] \in W/W_Q$),
we denote as in the previous subsection
$C_{[w]} = \overline {B.[w]} \subset G/P$ the corresponding Schubert
subvariety (resp. $D_{[x]} = \overline {B.[x]} \subset G/Q$).
In the rest of this article, I give a description of the $[w]$'s 
such that the $C_{[w]}$ is suitable, and of the
element in $W/W_Q$ corresponding to the dual Schubert variety, in the
fundamental cases. According to
section \ref{section_reduction}, this is enough to describe all dual varieties
of Schubert varieties. The strategy for this description is first to use
a $T$-fixed point argument, to reduce the task to a purely
combinatorial one. In the types $A$ and $D$, I then give an
explicit solution of this combinatorial problem. For the exceptional
cases, my description of the dual Schubert varieties is not really
explicit, but to compute them there is in principal only a finite
number of computations to make.

\lpara

So we fix the minimal $G$-representation $V$ such that $G/P \subset \p V$. We
denote $V = \oplus V_\lambda$ (resp. $V^* = \oplus V^*_\mu$)
the weight decomposition of $V$ (resp. $V^*$).
Let $C_{[w]} = \overline {B.[w]}$ be a Schubert
variety. Recall that $\overline{N_{[w]} C_{[w]}}$
denotes the variety of $y$'s in
$G/Q$ which are tangent to $C_{[w]}$ at $[w]$, see definition
\ref{cotangent}.

\begin{lemm}
\label{t_fixe}
Let $[x] \in W/W_Q$ such that 
$C_{[w]}^Q = \overline {B.[x]}$. We have 
$[x] \in \overline{N_{[w]} C_{[w]}}$.
\end{lemm}
\begin{proo}
First, notice that $C_{[w]}^Q = \overline
{B.\overline{N_{[w]}C_{[w]}}}$. In fact,
$B.\overline{N_{[w]}C_{[w]}}$ contains the set of $y$'s in $G/Q$
which are tangent at a point in $B.[w]$, therefore at a generic point
of $C_{[w]}$.

Let $\mu_0$ be the highest weight of $V^*$ and denote $\mu = x.\mu_0$. 
Let $y \in V^*$ such that $[y] \in C_{[w]}^Q \subset \p V^*$.
Use the weight decomposition
of $V^*$ to write $y = \sum_{\mu'} y_{\mu'}$.
Since $[y] \in C_{[w]}^Q$ which is the closure of the $B$-orbit of the weight
line of weight $\mu$,
$y_{\mu'} = 0$ if $\mu' \not \geq \mu$. Assume that 
$\forall [y] \in \overline{N_{[w]} C_{[w]}},y_{\mu}=0$.
It would then follow from
the first point that $\forall y \in C_{[w]}^Q, y_{\mu}=0$, contradicting
$[x] \in C_{[w]}^Q$.

Therefore, there exists $[y]$ in $\overline{N_{[w]}C_{[w]}}$ such that 
$y_\mu \not = 0$. Since $\overline{N_{[w]}C_{[w]}}$
is $T$-stable and $V_\mu$ is 
one-dimensional,
we have $[x] \in \overline{N_{[w]} C_{[w]}}$.
\end{proo}

I now explain how to compute
the element $[x]$ of the previous lemma. 
We want to take into account the case of $E_6/P_3$, which cotangent
bundle is not irreducible. Recall that in this case there is a natural
bundle morphism $\pi : T^* E_6/P_3 \rightarrow L \otimes \wedge^2 Q^*$. 
To have uniform notations, in the other cases we denote
$\pi : T^* G/P \rightarrow T^* G/P$ the identity and $\overline q = q$.

Decompose
$\pi(T_{[e]}^* G/P)$ as a sum of weight spaces for the action of a Levi
subgroup of $P$ : $\pi(T_{[e]}^* G/P) = \oplus T^*_\tau$, and write similarly
$\scal{I_{[e]}} = \oplus_\nu N_\nu$, where if $[w'] \in W/W_P$,
$\scal{I_{[w']}} \subset V^*$ denotes the linear span of the Schubert variety
$I_{[w']} \subset G/Q \subset \p V^*$ (see notation \ref{i_h}).
It can be easily checked directly on the
examples that all the weight spaces $T^*_\tau$ and $N_\nu$ have dimension 1.
The rational map 
$\overline q : \pi(T_{[e]}^* G/P) \dasharrow G/Q \subset \p V^*$
is given by a list
of polynomial functions $\pi(T_{[e]}^* G/P) \rightarrow N_\nu$
of the same degree $d$ and with values in the complex line $N_\nu$.
The polarisations of these polynomials yield $d$-linear maps
$T_{\tau_1}^* \times \cdots \times T_{\tau_d}^* \rightarrow N_\nu$, which
will be denoted $P_{\tau_1,\ldots,\tau_d;\nu}$; remark that the space of
such $d$-linear maps has dimension 1. Given $w \in W$, we denote
${\cal P}_w$ the set of weights $\nu$ such that there exist weights
$\tau_1, \ldots , \tau_d$ of 
$\pi(w^{-1}.N^*_{[w]} C_{[w]}) \subset \pi(T_{[e]}^* G/P)$
such that $P_{\tau_1,\ldots,\tau_d;\nu}$ does not vanish.
\lpara

\begin{prop}
\label{prop_combinatoire}
Let $[w] \in W/W_P$, with $w$ its minimal length representative.
The variety $C_{[w]}$ is suitable \iff ${\cal P}_w$ is not empty. In this case,
if we denote $[x] \in W/W_Q$ such that $C_{[w]}^Q = C_{[x]}$,
then $x.\mu_0$ equals $w.\mu_1$, where $\mu_1$ is the
lowest weight in ${\cal P}_w$.
\end{prop}
\begin{proo}
In fact, by lemma \ref{t_fixe} and its proof, $x.\mu_0$ is the lowest weight
$\mu$, if any,
such that the $\mu$-component of the restriction of the rational
map 
$\pi(T^*_{[w]} G/P) \dasharrow G/Q$ to $N^* C_{[w]}$
does not vanish identically. 

The weights of $\scal{I_{[w']}}$ for 
$[w'] \in W/W_P$ are some weights of $V^*$,
a set on which $W$ acts; therefore it makes sense to talk of
$w''.\mu'$, for $w'' \in W$ and $\mu'$ a weight of $\scal{I_{[w']}}$.  
I claim that $w$ induces an increasing bijection between the weights of
$\scal{I_{[e]}}$ and those of $\scal{I_{[w]}}$. 
The argument is similar to that of proposition \ref{prop_order_ideal}.
In fact, a weight of $\scal{I_{[e]}}$ can be
written as $v.\mu_0$, with $v \in W_P$ and $\mu_0$ the highest weight
of $V^*$. Given two such weights $v_1.\mu_0 \geq v_2.\mu_0$, we can assume
that $v_2$ is the minimal length representative of its class in 
$W_P / W_{P \cap Q}$.
Thus $v_1$ can be written as a product of some reflections which occur
in a reduced expression of $v_2$. Since $v_2 \in W_P$ and $w$ is a minimal
length representative, $l(wv_2) = l(w) + l(v_2)$. Therefore
$wv_1 \leq wv_2$ for the Bruhat order, and thus
$w.(v_1.\mu_0) \geq w.(v_2.\mu_0)$. This proves the claim.

Since the
rational map $\overline q : \pi(T^* G/P) \dasharrow G/Q$
is $G$-equivariant, the
weight $\mu$
of the proof of lemma \ref{t_fixe}
is also $w.\mu_1$, where $\mu_1$ is the lowest weight such
that the $\mu_1$-component of the restriction of the rational
map $T^*_{[e]} G/P \dasharrow G/Q$ to 
$w^{-1}.N^*_{[w]} C_{[w]}$ does not vanish identically.
Obviously $\mu_1$ is the lowest weight $\nu$ such that some
$P_{\tau_1,\ldots,\tau_d;\nu}$ 
with $\tau_1,\ldots,\tau_d$ weights of $w^{-1}.N^*_{[w]} C_{[w]}$,
does not vanish, so the proposition is proved.
\end{proo}

We illustrate our method with the easy example $G/P = \p V$.
Let $(e_1,\ldots,e_n)$ be a basis of $V$, let $k$ be an integer and let
$L_k = Vect(e_1,\ldots,e_k),\ M_k = Vect(e_{k+1},\ldots,e_n)$.
We consider the Schubert variety $X = \p L_k \subset \p V$ and compute
its dual variety. The corresponding element of the Weyl group is the
transposition $w=(1k)$. We have
$T_{[w]} X \simeq Hom(e_k,L_k/e_k)$, so
$w^{-1}. T_{[w]} X \simeq Hom(e_1,L_k/e_1)$ and so
$w^{-1}. N^*_{[w]} X \simeq Hom(M_k,e_1)$. Since $q$ is defined taking the
kernel, the lowest weight in $q(w^{-1}. N^*_{[w]})$ is 
$\mu_1 = -\epsilon_{k+1}$. We have $w.\mu_1 = \mu_1$, so that the dual
variety of $X$ is the $B$-orbit closure of $e_{k+1}^*$, as
expected. Note that in this example it would have been easier to compute
directly the lowest weight in
$q(N^*_{[w]}C_{[w]})$, instead of applying first $w^{-1}$ and then $w$. In fact,
this is what we will do to compute dual varieties of Schubert varieties
in the cases $G/P = \G(r,V)$ and $G/P = \G_Q^+(2p+1,4p+2)$.
\lpara

Recall from subsection \ref{subsection_carquois}
the definition of height of a vertex of the quiver
$Q_{G/P} = H_{G/P}$. We denote $h_0$ the maximal $h$ such that there exist
$\tau_1,\ldots,\tau_d \in H_{G/P}$, $\nu$ a weight of $I_{[e]}$ such
that $h(\tau_i) \geq h$ and $P_{\tau_1,\ldots,\tau_d;\nu} \not = 0$.
We have the following values for $h_0$ (I have also indicated the
height $h_{\max}$ of the heighest element of $Q_{G/P}$)~:
$$
\begin{array}{ccccc}
G/P & \G(r,n) & \G_Q^+(2p+1,4p+2) & E_6/P_1 & E_6/P_3\\
h_{\max} & n-1 & 4p-1 & 11 & 8\\
h_0 & \max(r,n-r) & 2p+1 & 8 & 5
\end{array}
$$

\begin{theo}
\label{theo_adapte}
The Schubert subvariety $C_{[w]}$ is suitable \iff all the vertices of
$Q_w$ have height at most $h_0-1$.
\end{theo}
\begin{proo}
Unfortunately, I don't know how to prove in a uniform way this
theorem. It will follow from propositions \ref{grass_adapte},
\ref{prop_spinoriel}, \ref{prop_e6_1_adapte} and
\ref{prop_e6_2_adapte}. The proof of these propositions also imply the
above given values of $h_0$.
\end{proo}


\subsectionplus{Case of Grassmannians}

\label{subsection_schubert_grassmannienne}

Recall that $V$ is an $n$-dimensional vector space.
We will
parametrise Schubert varieties in $\G(r,V)$ by increasing lists of $r$ 
integers,
instead of partitions, because duality will appear easier to formulate
in this way.
The list $(l_i)$ will correspond to the Schubert variety 
$C_l \subset \G(r,V)$ (resp. in $D_l \subset \G(r,V^*)$) which is the
$B$-orbit closure of the linear space spanned by the
$l_i$'s $T$-eigenvectors in $V$ (resp. in $V^*$).
For $x \in \{1,\ldots,n\}$, we will write $x \in l$ to mean that there
exists $i$ such that $x = l_i$.
The $T$-fixed points in $V$ (resp. $V^*$) will be denoted $e_i$
(resp. $e_i^*$).
The $T$-fixed point whose $B$-orbit is dense in $C_l$ 
(resp. $D_l$) will be denoted $x_l$
(resp. $y_l$).
The Bruhat order on Schubert cells is given by
$l \leq m$ \iff $\forall i,l_i \leq m_i$. If $x_1,\ldots,x_r$ are 
distinct integers
not necessarily increasing, we denote the list obtained reordering the
$x_i$ as $[x_1,\ldots,x_r]$.

Let $T_l$ denote $Vect(e_i : i \in l)$ and let
$Q_l$ denote $Vect(e_i : i \not \in l)$.
The tangent
space at $x_l$ identifies with
$Hom(Q_l,T_l)$.
A weight in this space is given by a couple 
$(x,y):x \in l,y \not \in l$.
Recall that the rational map
$$T_{x_l}^* \G(m,n) \simeq 
Hom(Q_l,T_l) \dasharrow \G(r,V^*)$$
is given by $\varphi \mapsto \ker \varphi$.
Thus the degree of $q$ is $r$ and
with the notations before proposition \ref{prop_combinatoire}, 
we have~:

\begin{fait}
\label{fait_P_grassmannienne}
The multilinear form
$P_{(x_1,y_1),\ldots,(x_r,y_r);l'}$ does not vanish \iff the $x_i$'s 
and the $y_i$'s are all distinct, and $l'$ is the set of the $y_i$'s.
\end{fait}

\noindent
Given a list $l$, we consider the list $l^*$ defined inductively by
$$
l^*_i = \min \{ y : y>y_{i-1}, y>x_i, \forall j, y\not = x_j \}.
$$

\begin{lemm}
We have $\forall i,l_i^* \leq n$ \iff $\forall i \in \{1,\ldots,r\},
l_i<n+2i-2r$.
\end{lemm}
\noindent
In terms of partitions, this means that the $i$-th part
must be at least $r+1-i$.
\begin{proo}
Let $i$ be an integer. The integers $l_j$ for $j>i$ and $l_j^*$ for
$j \geq i$ are strictly greater than $l_i$
and distinct, so the lemma follows.
\end{proo}

\noindent
As the following proposition
shows, $l \mapsto l^*$ is the combinatorial model for the duality of
Schubert varieties in Grassmannians~:

\begin{prop}
\label{grass_adapte}
$C_l$ is suitable \iff $\forall i \in \{1,\ldots,r\},l_i < n+2i-2r$.
If $C_l$ is suitable then $C_l^Q = D_{l^*}$.
\end{prop}
\begin{proo}
With the previous notations, the weights of the conormal
space $N^*_{x_l} C_l$ are the couples
$(x,y)$ with $x \in l,y \not \in l$, and $y>l_x$.

By proposition \ref{prop_combinatoire} and the comment after it,
$C_l$ is suitable \iff there are lists $[y_1,\ldots,y_r]$
and $x_1,\ldots,x_r$ with
$(x_i,y_i)$ a weight of $N^*_{x_l} C_l$ and
$P_{(x_1,y_1),\ldots,(x_r,y_r);l'} \not = 0$.
In this case, if we denote $l'$ the list
such that $C_l^Q = D_{l'}$, then $l'$ is the minimal possible such list.
Moreover, in order that 
$P_{(x_1,y_1),\ldots,(x_r,y_r);l'} \not = 0$,
all $x_i$ must be distinct and we must have
$\{x_i\} = l$, so we
may assume by symmetry that $x_i = l_i$.
It is easy to check that the set of such $l'$ is not
empty \iff $\forall i \in \{1,\ldots,r\},l_i<n+2i-2r$. In fact, if
$l_i \geq n+2i-2r$, then the values $y_j$ and $x_j$ for 
$i \leq j \leq r$ must be distinct
and between $n+2i-2r$ and $n$, a contradiction. Conversely, if
$\forall i,l_i < n+2i-2r$, one may choose $y_i = l_i^*$.

We now show that $l^*$ is indeed the minimal list. Let
$[y_1,\ldots,y_r]$ be any list with 
$\forall i, y_i>l_i$ and $\forall i,j,y_i \not = l_j$.
Let $(z_1,\ldots,z_r)$ be the corresponding ordered list
(ie $\{y_1,\ldots,y_r\} = \{z_1,\ldots,z_r\}$ and
$z_1 < z_2 < \cdots < z_r$). Then we have $z_1 > l_1$
and $\forall j, z_1 \not = l_j$,
so $z_1 \geq l_1^*$. Say $z_i = y_{\sigma(i)}$. If $z_1 < l_2$ then
$\sigma(1)=1$. Thus in any case $z_2 > l_2$, so $z_2 \geq l^*_2$. By induction
it follows that $\forall i,l^*_i \leq z_i$, so $l^*$ is the minimal possible
list, and proposition \ref{prop_combinatoire} finishes the proof.
\end{proo}

We illustrate this proposition with two examples. The array below
computes two dual varieties in $\G(3,8)$. It pictures the fact that
for $l=(2,4,5)$ we have $l^*=\lambda=(3,6,7)$, and that for $m=(2,4,6)$ we have
$m^*=\mu=(3,5,7)$~:
$$
\begin{array}{ccc}
\young(\hfil l\lambda\hfil\hfil\hfil\hfil\hfil,\hfil\hfil\hfil l\hfil\lambda\hfil\hfil,\hfil\hfil\hfil\hfill l\hfil\lambda\hfil)
&
\hspace{1cm}
&
\young(\hfil m\mu\hfil\hfil\hfil\hfil\hfil,\hfil\hfil\hfil m\mu\hfil\hfil\hfil,\hfil\hfil\hfil\hfil\hfil m\mu\hfil)
\end{array}
$$

\noindent
Note that we have $C_l \subset C_m$ but we don't have 
$D_{l^*} \supset D_{m^*}$~: contrary to the case $G/P = \p V$, duality
of Schubert cells is no longer decreasing.


\subsectionplus{Case of spinor varieties}

\label{subsection_spinoriel}

Schubert cells in $\G_Q^+(2p+1,4p+2)$
(resp. $\G_Q^-(2p+1,4p+2)$) are parametrised by lists of
$+$ and $-$ signs, with an odd number of $+$ (resp. $-$) signs.
The generic $T$-fixed point corresponding to the list
$(\eta_i)$ is the subspace generated by $e_i^{\eta_i}$, and will be
denoted $x_\eta$.
Schubert cells are also parametrised by strict partitions of size $2p$
(or subsets of $\{1,\ldots,2p\}$), the correspondance
being that we set $x \in \lambda$ ($1 \leq x \leq 2p$)
if $\eta_{2p+1-x} = -\ $.

\begin{defi}\
\begin{itemize}
\item
If $(\eta_i),1\leq i \leq 2p+1,$ is a sequence of signs and $j$
is an integer, we denote
$\varphi(\eta,j)$ the sequence $\eta'$ of signs such that
$\eta'_i = \eta_i$ for exactly all $i$'s but $j$.
\item
A sequence $(\eta_i),1\leq i \leq 2p+1,$ of signs is admissible if
$$
\forall i\in \{1,\ldots,p\},
\#\{j : 1 \leq j \leq 2i, \eta_j = + \} \geq i.
$$
Assume that $\eta$ is admissible~:
\item
If there exists $i\leq p+1$ such that
$\# \{j\ :\ 1\leq j\leq 2i-1,\eta_j=+\} = i-1$, then let $i_0$ be the
minimal such $i$, and set $\eta^* = \varphi(\eta,2i_0-1)$.
\item
Otherwise there exists $i$ such that
$$
\forall k\geq i,\ \#\{j:j\leq k,\eta_j=+\} > \# \{j:j\leq
k,\eta_j=-\}.
$$
Let $i_0$ be the minimal such $i$ and set $\eta^* =
\varphi(\eta,i_0)$.
\end{itemize}
\end{defi}
\noindent
If there does not exist $i\leq p+1$ such that
$\# \{j\ :\ 1\leq j\leq 2i-1,\eta_j=+\} = i-1$,
then $\# \{ j\ :\ 1 \leq j \leq 2p+1 \} \geq p+1$, so $i=2p+1$
satisfies the condition of the last point of this definition.

\lpara

Let $\eta$ be fixed.
Since the positive roots are $\epsilon_i \pm \epsilon_j$ with $i<j$,
the restriction of the Bruhat order on the set of $\varphi(\eta,j)$
is given by~:

\begin{fait}
\label{fait_ordre_spin}
We have $\varphi(\eta,i) \leq \varphi(\eta,j)$ for the Bruhat order
\iff
$$
\mbox{or }\left \{
\begin{array}{l}
\eta_i = \eta_j = + \mbox{ and } i\leq j\\
\eta_i = + \mbox{ and } \eta_j = -\\
\eta_i = \eta_j = - \mbox{ and } i \geq j.
\end{array}
\right .
$$
\end{fait}

\begin{prop}
\label{prop_spinoriel}
$C_\eta$ is suitable \iff $\eta$ is admissible, and in this case
$C_\eta^Q = D_{\eta^*}$.
\end{prop}
\begin{proo}
Recall that $x_\eta \in G/P$ denotes the linear space spanned by
$e_i^{\eta_i}$. It is well-known that $T_{x_\eta} G/P$ identifies
with $\wedge^2 {Vect(e_i^{\epsilon_i})}^*$. Moreover,
$N^*_{x_\eta} C_\eta \subset \wedge^2 {Vect(e_i^{\epsilon_i})}$ 
is generated by $e_i^{\eta_i} \wedge e_j^{\eta_j}$ for $(i,j)$ such
that $\eta_i = +$ and $i<j$. In fact, with the notations of 
\cite[PLANCHE IV]{bourbaki}, the weight of $x_\eta$ is
$\rho_\eta = \frac 12 \sum \eta_i \epsilon_i$, and the weights of
$T_{x_\eta} C_\eta$ are the weights of the form
$\frac 12 \sum \eta'_i \epsilon_i$ which can be expressed as
$\rho_\eta + \alpha$, where $\alpha$ is a positive root. Therefore the
claim follows from the fact that the positive roots are 
$\epsilon_i \pm \epsilon_j$ with $i<j$.

The weights of $T^*_{x_\eta} G/P \simeq \wedge^2 {Vect(e_i^{\epsilon_i})}$
are parametrised by couples $(x,y)$ of integers, with $x<y$.
Now let $x_k,y_k,1 \leq k \leq p,$ be integers with $x_k < y_k$. With
the notations of subsection \ref{subsection_combinatoire}, $\mu$ is given by a
set of polynomials of degree $p$ and the $p$-multilinear map
$P_{(x_k,y_k);\eta'}$ does not vanish \iff the $x_k$'s and the $y_k$'s
are all distinct and $\eta'_i = \eta_i$ for exactly all $i$'s which
belong to the set $U := \{x_k\} \cup \{ y_k \}$.

Given the previous description of 
$N^*_{x_\eta} C_\eta$, the Schubert
variety $C_\eta$ will be suitable \iff we can find
$(x_k,y_k)$ such that
\begin{equation}
\label{condition_spin}
x_k < y_k, \mbox{ the }x_k,y_k \mbox{ are all distinct, and }
\eta_{x_k} = +.
\end{equation}
Therefore, for all
$i$'s with $1\leq i \leq p$, we have the inequality
$$
2i-1 \leq \#(U \cap [1,2i]) \leq 2\# \{j\ :\ 1\leq j\leq 2i,\eta_j = +
\}.
$$
This implies that $\eta$ should be admissible.

Conversely, assuming that $\eta$ is admissible, let us consider the
following algorithm which produces a list of distinct elements
$(x_k,y_k)$ with $\eta_{x_k} = +$ and $x_k < y_k$.
If $\forall i>1, \eta_i = +$, set $x_k = 2p$ and $y_k = 2p+1$.
Otherwise, let $i_0$ be the minimal $i>1$ such that $\eta_i = -$;
set $x_1 = i_0 - 1$ (the fact that $\eta$ is admissible garanties that
even in the case $i_0=2$, we have $\eta_{i_0-1} = +$) and $y_1 =
i_0$. Remove $\eta_{x_1}$ and $\eta_{y_1}$ from the list $\eta$ : this
new list is again admissible, as one checks readily. Therefore, it is
possible to define $(x_k,y_k)$ for $k \geq 2$ inductively.

\lpara

We therefore have proved that $C_\eta$ is suitable \iff $\eta$ is
admissible. Let us now compute the dual variety. Assume first that
there exists $i \in \{1,p+1\}$ such that
\begin{equation}
\label{condition}
\# \{j : 1 \leq j \leq 2i-1, \eta_j = + \} = i-1.
\end{equation}
Let $i_0$ be the minimal such $i$. Admissibility of $\eta$ implies
that $\# \{j\ :\ 1 \leq j \leq 2i_0-2, \eta_j = + \} = i_0 - 1$ and so
$\eta_{2i_0-1} = -$. Therefore, if $(x_k,y_k)$ is any sequence
satisfying (\ref{condition_spin})
and $U = \{x_k\} \cup \{y_k\}$, there exists $j \leq 2i_0-1$
such that $\eta_j = -$ and $j \not \in U$. Thus if there exists
$(x_k,y_k)$ such that $P_{(x_k,y_k);\eta'} \not = 0$ for some $\eta'$,
this implies $\eta' \geq \varphi(\eta,2i_0-1)$
(recall fact \ref{fait_ordre_spin}).
Conversely, the previous algorithm produces a sequence $(x_k,y_k)$ for
which it is easy to see that
$P_{(x_k,y_k);\varphi(\eta,2i_0-1)} \not = 0$. Thus
$\eta^* = \varphi(\eta,2i_0-1)$ is the lowest
list one can obtain in this way,
so that $C_\eta^Q = D_{\eta^*}$ as claimed in this case.

Assume finaly that (\ref{condition}) holds for no $i \in \{1,p+1\}$.
Therefore, as we have seen, there exists
$i$ (for example $i=2p+1$) such that
$$
\forall k \geq i, \# \{j\ :\ j \leq k,\eta_j = +\}
>
\# \{j\ :\ j\leq k,\eta_j = -\}.
$$
Let $i_0$ be the minimal such $i$. Obviously, if $i$ is any integer,
and $(x_k,y_k)$ satisfies (\ref{condition_spin})
and $x_k,y_k \not = i$, then
$$\forall k \geq i, \# \{j\ :\ j \leq k,\eta_j = +,j\not = i\}
\geq
\# \{j\ :\ j\leq k,\eta_j = -,j \not = i\},$$ so that
$i \geq i_0$. So $P_{(x_k,y_k);\eta'} \not = 0$ implies
$\eta' \geq \varphi(\eta,i_0)$. Again, the explicit algorithm provides
a sequence $(x_y,y_k)$ such that
$P_{(x_k,y_k);\varphi(\eta,i_0)} \not = 0$, and therefore
$D_{\varphi(\eta,i_0)}$ is the dual variety of $C_\eta$.
\end{proo}

\subsectionplus{Case of $E_{6,I}$}

\label{subsection_e6_1}

We now consider the exceptional cases when $G$ is of type
$E_6$. Recall that there are two possibilities for $(P,Q)$~: either
they correspond to the roots $(\alpha_1,\alpha_6)$ or
$(\alpha_3,\alpha_5)$. In each case I explain in which case
$P_{\tau_1,\ldots,\tau_d;\nu}$ does not vanish. Using subsection
\ref{subsection_combinatoire}, this describes in principle all dual
varieties to Schubert varieties, although I will not give a simple
combinatorial recipy for this correspondance (note however that to
give such a description there is ``only'' a finite number of computations to
do). My description of which $P_{\tau_1,\ldots,\tau_d;\nu}$ don't
vanish will however yield a simple caracterisation of the suitable Schubert
varieties.

As we have seen in subsection \ref{subsection_fondamental}, a Levi
factor $L$ of $P_1$ is isomorphic with $\C^* \times Spin_{10}$, and
$T_{[e]} E_6/P_1$ identifies with a
$16$-dimensional spinor representation
of $L$. Moreover, the closed $L$-orbit in $\p T^*_{[e]} G/P$
identifies with a $L$-homogeneous spinor variety~: it is a connected
component of the variety
parametrising isotropic linear spaces of dimension 5 in a certain
quadratic vector space of dimension 10 that we will denote $M$.
It is proved in \cite[corollary 3.2]{flop_scorza}
that $I_{[e]} \subset \p M$ is the
corresponding 8-dimensional quadric acted upon by $L$ 
and that the rational map 
$T^*_{[e]} E_6/P_1 \dasharrow I_{[e]}$
is induced by the unique $L$-equivariant
quadratic map $T^*_{[e]} G/P \rightarrow M$. The polarisation 
${\cal P} : T^*_{[e]} E_6/P_1 \times T^*_{[e]} E_6/P_1 \rightarrow M$ of this
equivariant map has the following geometric interpretation~: for 
$x,y \in T^*_{[e]} E_6/P_1$ representing points of the spinor variety
corresponding to the isotropic linear spaces $L_x,L_y$, the class
of ${\cal P}(x,y)$ in $\p M$ is the intersection of $L_x$ and $L_y$ if
this intersection has dimension 1, and ${\cal P}(x,y)=0$ otherwise.

Denote as in subsection \ref{subsection_spinoriel}
$(e_1^+,\ldots,e_5^+,e_1^-,\ldots,e_5^-)$ a base of $M$ such that the
quadratic form $Q$ satisfies
$Q(\sum x_i^+ e_i^+ + \sum x_i^- e_i^-) = \sum x_i^+x_i^-$. An $L$-weight
of $M$ can therefore be denoted 
$\nu \in \{1^+,\ldots,5^+,1^-,\ldots,5^-\}$, and a weight 
$\eta = (\eta_1,\ldots,\eta_5)$ in
$T^*_{[e]} G/P$ corresponds to a list of plus or minus signs, with an
odd number of plus signs.
The condition for $P_{\eta,\eta';\nu}$ not to vanish is thus that
$\eta$ and $\eta'$ have exactly one sign in common which is $\nu$.

From this description we can describe the suitable Schubert varieties.
In the array below, I recall the quiver of $E_6/P_1$ and 
I define an
element $[w_{max}] \in W/W_P$ by its subquiver~:

$$
\begin{array}{ccc}
\input{e6} && \input{w_max_1}\\
\mbox{Quiver of }E_6/P_1  &  \hspace{2cm}  &  \mbox{Subquiver of }[w_{max}]
\end{array}
$$

\noindent
We have the following proposition~:

\begin{prop}
\label{prop_e6_1_adapte}
Let $[w] \in W/W_P$. Then the Schubert variety
$C_{[w]}$ is suitable \iff $[w] \leq [w_{max}]$.
\end{prop}
\begin{proo}
Below we give the Hasse diagram $H$ (resp. $H^*$) of the $L$-module
$T_{[e]} E_6/P_1$ (resp.  $T^*_{[e]} E_6/P_1$), which, by proposition
\ref{hasse}, is isomorphic with the quiver of $E_6/P_1$~:

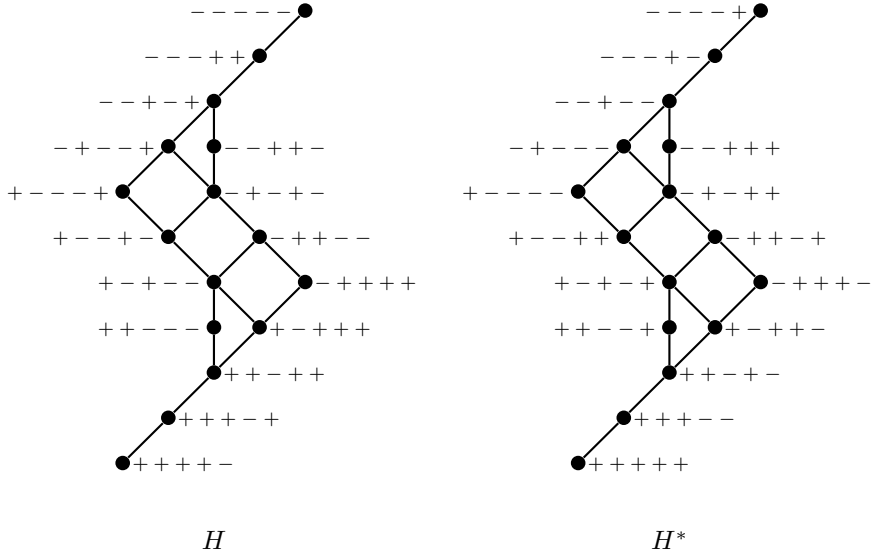
$$
\begin{array}{cc}
\input{h} & \input {hd}\\
\\
H & H^*
\end{array}
$$

Let $\iota : H \rightarrow H^*$ be
induced by the map $\eta \mapsto -\eta$ (in terms of quivers, this
corresponds to the obvious symmetry).
Let $[w] \in G/P$ and $Q_{[w]} \subset H$ the corresponding subquiver, marking
the weights of $w^{-1}.T_{[w]} C_{[w]}$.
Thanks to proposition \ref{prop_combinatoire},
the proposition amounts to the fact that
we can find two weights $\eta,\eta' \in H^* - \iota(Q_{[w]})$
and which have only one sign in common \iff 
$Q_{[w]} \subset Q_{[w_{max}]}$. This may be seen as follows~:

\begin{itemize}
\item
If $Q_{[w]} \subset Q_{[w_{max}]}$, we can set
$\eta = (++--+)$ and $\eta' = (+-++-)$ to check that the corresponding
Schubert variety is suitable (below the subset $\iota(Q_{[w_{max}]})$ 
is drawn).
\item
If $Q_{[w]}$ contains the vertex corresponding to the
weight $(-+--+)\ $, $\iota(Q_{[w]})$ contains the subset drawn below.
Thus $\eta$ and $\eta'$ are weights which begin with $++$, so
they have two common signs. The corresponding Schubert variety
is not suitable.
\item
The last case is that the subquiver contains the vertex corresponding
to the weight $(--++-)$. 
Thus $\eta$ and $\eta'$ have at least 3 plus signs among the 
4 first signs, and
therefore have at least 2 common signs.
\end{itemize}

$$
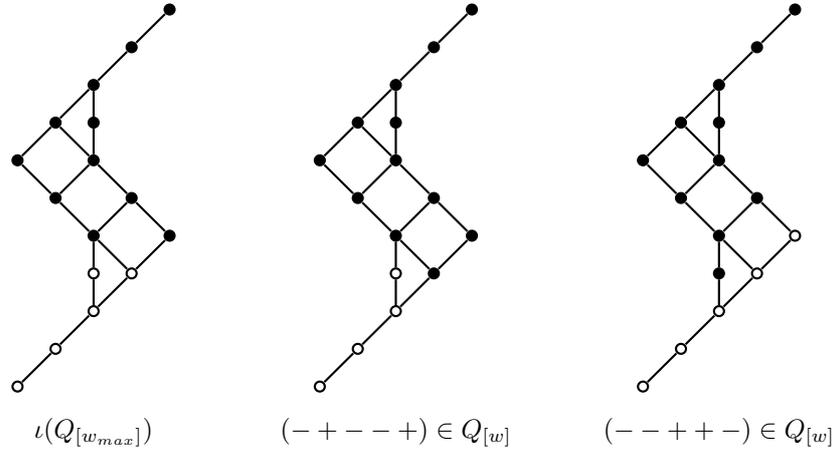
\begin{array}{ccccc}
\input{cas0} && \input{cas1} && \input{cas2} \\
\iota(Q_{[w_{max}]}) & \hspace{.5cm} & 
(-+--+) \in Q_{[w]} & \hspace{.5cm} & (--++-) \in Q_{[w]}
\end{array}
$$

\end{proo}

\begin{exem}
\label{exem_p5e6}
Let $X \subset E_6/P_1$ be a linear subspace of maximal dimension 5.
Then $X^Q \subset E_6/P_6$ is also a linear subspace of dimension 5.
\end{exem}
\noindent
Therefore, this provides, in our setting, a new example of a variety
which is isomorphic to its dual. Similar examples in the usual setting
$X \subset \p^n$
are projective subspaces, quadrics, 
$\G(2,2p+1)$, and the spinor variety
$\G_Q^+(5,10)$.
\begin{proo}
Let $X \subset E_6/P_1$ be a linear subspace of dimension 5.
The variety parametrising linear subspaces of maximal dimension 5 is given by
Tits shadows, according to \cite[theorem 4.3]{landsberg}.
In particular, it is a
homogeneous variety, and so we can assume that $X$ is the Schubert
variety corresponding to the Weyl group element
$w=s_6s_5s_4s_3s_1$. The corresponding quiver $Q_w$ and $\iota(Q_w)$
follow; we have also drawn the quiver $Q_{w^*}$ of the dual variety.
$$
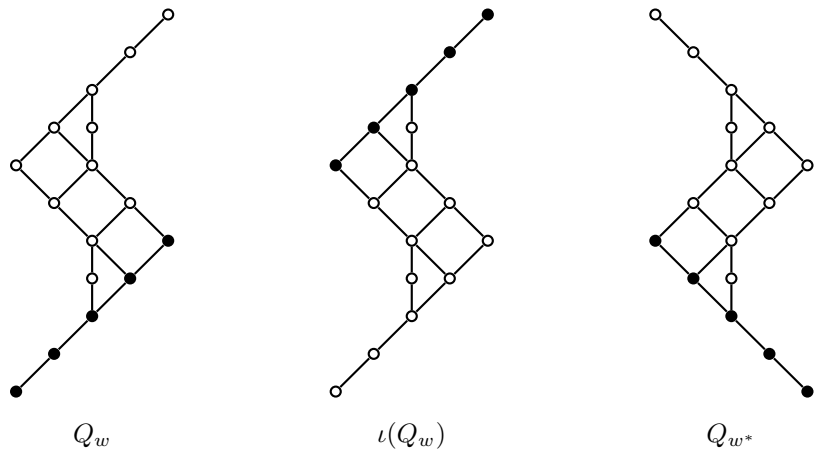
\begin{array}{ccccc}
\input{p5}  &&  \input{ip5}  &&  \input{dp5}\\
Q_w & \hspace{1cm} & \iota(Q_w) & \hspace{1cm} & Q_{w^*}
\end{array}
$$

According to proposition \ref{prop_combinatoire}, we must look for two
weights of $H^*$ not in $\iota(Q_w)$ which have only one common
sign. If this common sign is a minus sign, then among these two
weights there are 6 minus signs. Therefore one of the weight has 4
minus signs which is impossible given $\iota(Q_w)$ and the Hasse
diagram $H^*$. Since $(--+++)$ and $(++--+)$ are weights not in
$\iota(Q_w)$, the lowest weight is $5^+$. Note that this weight is
obtained applying $w_6 = s_3s_4s_5s_6$ to the highest weight.
Therefore, proposition
\ref{prop_combinatoire} shows that the dual variety to $X$ corresponds
to the class of $w.w_6 = s_6s_5s_4s_3s_1s_3s_4s_5s_6$ modulo $W_6$.
Modulo $W_6$, we have
$$
\begin{array}{cccccc}
& s_6s_5s_4s_3s_1s_3s_4s_5s_6 & = & s_6s_5s_4s_1s_3s_1s_4s_5s_6 & = &
s_6s_5s_4s_1s_3s_4s_5s_6s_1\\
\equiv & s_6s_5s_1s_4s_3s_4s_5s_6 & = & s_6s_5s_1s_3s_4s_3s_5s_6 &
\equiv & s_6s_5s_1s_3s_4s_5s_6\\
= & s_6s_1s_3s_5s_4s_5s_6 & \equiv & s_6s_1s_3s_4s_5s_6 &
\equiv & s_1s_3s_4s_5s_6.
\end{array}
$$

This proves the claim.
\end{proo}

\subsectionplus{Case of $E_{6,II}$}

\label{subsection_e6_2}

Let $\alpha \in E_6/P_3$ be the base point.
Recall that there is a surjection
$\pi : T^*_\alpha E_6/P_3 \rightarrow 
L_\alpha \otimes \wedge^2 Q_\alpha^*$ and that the
rational map $q:T^*_\alpha E_6/P_3 \dasharrow I_\alpha$ is induced by a
rational map $\overline q : L_\alpha \otimes \wedge^2 Q^*_\alpha 
\dasharrow I_\alpha = \G(2,Q_\alpha^*)$.

The description of $\overline q$ given in subsection
\ref{subsection_fondamental} implies that
$\overline q$ has degree 6. Consider as in subsection
\ref{subsection_combinatoire} its polarisation, with
coordinates denoted $P$. In order to give a non-vanishing
criterium for $P$, let us introduce some notation. 
Let $e_1,e_2$ be a basis of $L_\alpha$ and 
$f_1^*,\ldots,f_5^*$ a basis of $Q_\alpha^*$.
The weight of the
vector $e_1 \otimes (f_i^* \wedge f_j^*)$ with $i<j$ will be denoted $ij$,
and the weight of the vector $e_2 \otimes (f_i^* \wedge f_j^*)$ will be
denoted $\underline{ij}$.
Finally, the weight of 
$f_i^* \wedge f_j^*$
will be denoted $ij^*$.

Let $\tau_1,\ldots,\tau_6$ be weights of
$L_\alpha \otimes \wedge^2 Q_\alpha^*$ and $\nu$ a weight of 
$\wedge^2 Q_\alpha^*$,
if $P_{\tau_1,\ldots,\tau_6,\nu} \not = 0$, then
$\# \{k:\tau_k \in \{ij\}\ \} = 
\# \{k:\tau_k \in \{\underline{ij}\}\ \} = 3$.
So we assume that this is the case and that
$\tau_1,\tau_2,\tau_3$ (resp. $\tau_4,\tau_5,\tau_6$) are of the form
$ij$ (resp. $\underline{ij}$).

With this setting,
$P_{i_1j_1,i_2j_2,i_3j_3,
\underline{k_1l_1},\underline{k_2l_2},\underline{k_3l_3}
;mn^*}$
will not vanish if and only if, up to permuting the three first weights and the
three last, we have that $i_1,j_1,i_2,j_2$ (resp. $k_1,l_1,k_2,l_2$) 
are all dinstinct; say they take all values in $\{1,\ldots,5\}$ except
$u$ (resp. $v$). Moreover we must have $u \in \{k_3,l_3\}$ (resp.
$v \in \{i_3,j_3\}$), say $\{k_3,l_3\} = \{u,u'\}$ (resp.
$\{i_3,j_3\} = \{v,v'\}$). Finally, we must have $u' \not = v'$ and
$\{u',v'\} = \{m,n\}$.

In principle, this combinatorial rule describes dual varieties of
Schubert cells in this case. However, as in the case of $E_{6,I}$, we
can be more precise as far as suitability is concerned. The Hasse
diagram $H$ of
$L_\alpha^* \otimes \wedge^2 Q_\alpha$ is given below, as well as a subquiver
denoted $Q_{max}$. I have also indicated the Hasse diagram $H^*$ of
$L_\alpha \otimes \wedge^2 Q_\alpha^*$~:

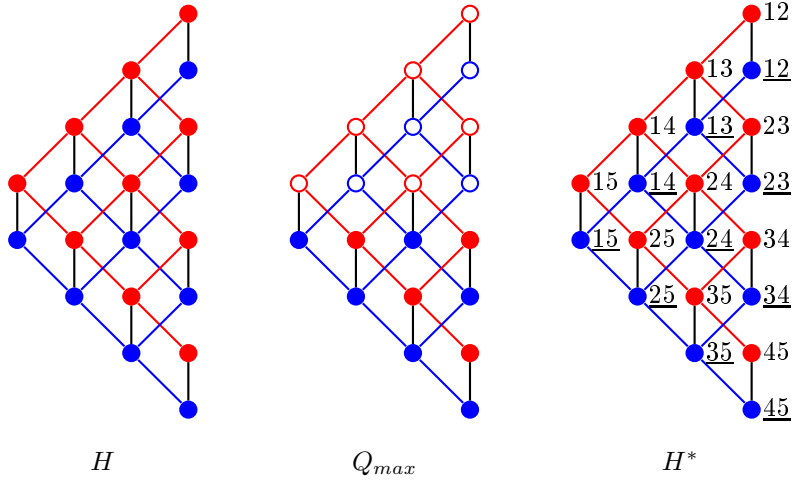
$$
\begin{array}{ccccc}
\input{hasse_e62} && \input{qmax} && \input{e6_2}\\
H && Q_{max} && H^*
\end{array}
$$

In these pictures and the following, the weights $ij$ are drawn in red
and the weights $\underline{ij}$ are drawn in blue.
Given a Schubert cell $C_{[w]}$, with $w$ the minimal length
representative of $[w] \in W/W_P$, recall that we associated
(not injectively) to this Schubert cell
a subquiver $Q_{[w]}$ of $H$ in subsection \ref{subsection_carquois}.

\begin{prop}
\label{prop_e6_2_adapte}
The Schubert cell $C_{[w]}$ is suitable \iff we
have $Q_{[w]} \subset Q_{max}$.
\end{prop}
\begin{proo}
As in the preceeding subsection, we define
$\iota:H \rightarrow H^*$ given by
$\eta \mapsto -\eta$.
The weights of $H^*$ have been given above.

Since $q$ is induced by $\overline q$,
a Schubert variety $C_{[w]}$ will be suitable \iff
the rational map 
$\overline q:L_\alpha \otimes \wedge^2 Q_\alpha^* \dasharrow I_\alpha$
is defined generically on $\pi(w^{-1}.N^*_{[w]} C_{[w]})$; equivalently,
$\overline q$ should be defined on the orthogonal of
$w^{-1}.T_{[w]} C_{[w]} \cap L_\alpha^* \otimes \wedge^2 Q_\alpha$ in
$L_\alpha \otimes \wedge^2 Q_\alpha^*$.
Equivalently again, we should be able to
find 6 weights $\tau_k$ not in $\iota(Q_{[w]})$ and some integers $i,j$
such that $P_{\tau_k;ij^*}$ does not vanish.

In case $Q_{[w]}$ is
included in $Q_{max}$, we can consider the weights
$34, 25, 34$, $\underline{15}$, $\underline{24}$, $\underline{15}$
(the corresponding subset $\iota(Q_{max})$ is drawn below), which satisfy
the relation
$P_{34,25,34,\underline{15},\underline{24},\underline{15};45^*} \not
= 0$ and do not belong to $\iota(Q)$. Otherwise, there are four cases~:
\begin{itemize}
\item
If $Q_{[w]}$ contains the weight $15$, by proposition
\ref{prop_order_ideal}, $\iota(Q_{[w]})$ contains the
corresponding subset in the array below. The remaining weights are of
the form $ij$ or $\underline{ij}$ with $1<i<j$, so the Schubert variety cannot
be suitable.
\item
If $Q_{[w]}$ contains the weight $\underline{14}$, 
the remaining weights are of the form
$ij$ or $\underline{kl}$ with $2<i<j$ (see the array below), so again
the Schubert variety cannot be suitable.
\item
If $Q_{[w]}$ contains the weight $24$, let 
$i_1j_1,i_2j_2,i_3j_3,\underline{i_4j_4},\underline{i_5j_5},
\underline{i_6j_6}$\vspace{1mm} be a list of weights not in $\iota(Q_{[w]})$.
Note that we have
$i_kj_k \in \{12,13,14,15,34\}$ for all $k$. Assume that there exists
$kl^*$ such that 
$P_{i_1j_1,i_2j_2,i_3j_3,\underline{i_4j_4},\underline{i_5j_5},
\underline{i_6j_6};kl^*} \not = 0$. This implies that the
only integer $x$ (resp. $y$)
which does not belong to the set $\{i_1,j_1,i_2,j_2\}$
(resp. $\{i_4,j_4,i_5,j_5\}$) is
either 2 or 5. This integer must
therefore belong to the set $\{i_3,j_3\}$ (resp.
$\{i_3,j_3\}$), so we must have $\{i_3,j_3\} = \{1,x\}$
(resp. $\{i_6,j_6\} = \{1,y\}$). This implies that $k=l=1$, a
contradiction.
\item
Assume finaly that $Q_{[w]}$ contains the weight $\underline{23}$.
In this case all the
weights which are not in $\iota(Q_{[w]})$
and of the form $ij$ satisfy $j=5$. Again,
the Schubert variety is not suitable.
\end{itemize}

$$
\begin{array}{ccccc}
\input{Cas0} & \hspace{.9cm} & \input{Cas1} & \hspace{.9cm} & \input{Cas2}\\
\mbox{Case } \iota(Q_{max}) && \mbox{Case } {15} \in Q_{[w]} &&
\mbox{Case } \underline{14} \in Q_{[w]}
\end{array}
$$
$$
\begin{array}{ccc}
\input{Cas3} && \input{Cas4}\\
\mbox{Case } {24} \in Q_{[w]} & \hspace{1cm} & 
\mbox{Case } \underline{23} \in Q_{[w]}
\end{array}
$$

\end{proo}




\vskip .2cm

\noindent
Author's address :

\noindent
Laboratoire de Math\'ematiques Jean Leray UMR 6629

\noindent
2 rue de la Houssini\`ere - BP 92208 - 44322 Nantes Cedex 3

\vskip .1cm
\noindent
Pierre-Emmanuel.Chaput@math.univ-nantes.fr

\end{document}

%% file: commandes_pe.tex


\newcommand{\N}{\mathbb{N}}
\newcommand{\Z}{\mathbb Z}
\newcommand{\R}{\mathbb{R}}
\newcommand{\Q}{\mathbb{Q}}
\newcommand{\C}{\mathbb{C}}
\newcommand{\G}{\mathbb{G}}
\renewcommand{\H}{\mathbb{H}}
\renewcommand{\O}{\mathbb{O}}
\newcommand{\F}{\mathbb{F}}
\renewcommand{\S}{\mathbb{S}}

\renewcommand{\a}{{\cal A}}
\newcommand{\az}{\a_\Z}
\newcommand{\ak}{\a_k}

\newcommand{\rc}{\R_\C}
\newcommand{\cc}{\C_\C}
\newcommand{\hc}{\H_\C}
\newcommand{\oc}{\O_\C}

\newcommand{\rk}{\R_k}
\newcommand{\ck}{\C_k}
\newcommand{\hk}{\H_k}
\newcommand{\ok}{\O_k}

\newcommand{\rz}{\R_\Z}
\newcommand{\cz}{\C_\Z}
\newcommand{\hz}{\H_\Z}
\newcommand{\oz}{\O_\Z}

\newcommand{\RR}{\R_\R}
\newcommand{\CR}{\C_\R}
\newcommand{\HR}{\H_\R}
\newcommand{\OR}{\O_\R}

\newcommand{\re}{\mathtt{Re}}
\newcommand{\str}{\mathfrak{str}}
\newcommand{\rank}{\mathtt{rk}}


\newcommand{\dem}{\noindent \underline {\bf D\'{e}monstration :} }
\newcommand{\pr}{\noindent {\bf Proof :} }
\newcommand{\indic}{\noindent \underline {\bf Indication :} }
\newcommand{\rem}{\noindent \underline {\bf Remarque :} }
\newcommand{\rek}{\noindent \underline {\bf Remark :} }
\newcommand{\rks}{\noindent \underline {\bf Remarks :} }
\newcommand{\fin}{\begin{flushright} \vspace{-16pt}
 $\bullet$ \end{flushright}}
\newcommand{\lpara}{
\ \vspace{3pt}

\noindent}
\newcommand{\para}{
\

\

\noindent}
\newcommand{\Para}{
\

\

\

\noindent}

\newcommand{\sectionplus}[1]{\section{#1} \vspace{-5mm} \indent}
\newcommand{\subsectionplus}[1]{\subsection{#1} \vspace{-5mm} \indent}


\newcommand{\dual}{{\bf v}}
\newcommand{\com}{\mathtt{Com}}
\newcommand{\rg}{\mathtt{rg}}

\newcommand{\g}{\mathfrak g}
\newcommand{\h}{\mathfrak h}
\renewcommand{\u}{\mathfrak u}
\newcommand{\n}{\mathfrak n}
\newcommand{\z}{\mathfrak z}
\newcommand{\e}{\mathfrak e}
\newcommand{\plie}{\mathfrak p}
\newcommand{\q}{\mathfrak q}
\newcommand{\liesl}{\mathfrak {sl}}
\newcommand{\so}{\mathfrak {so}}
\newcommand{\spin}{\mathfrak {spin}}




\newcommand{\dynkina}[2]{
\setlength{\unitlength}{2.5mm}
\begin{picture}(#1,3)
\put(0,0){$\circ$}
\multiput(2,0)(2,0){#2}{$\circ$}
\multiput(.73,.4)(2,0){#2}{\line(1,0){1.34}}
\end{picture}}

\newcommand{\dynkinap}[3]{
\setlength{\unitlength}{2.5mm}
\begin{picture}(#1,3)
\put(0,0){$\circ$}
\multiput(2,0)(2,0){#2}{$\circ$}
\multiput(.73,.4)(2,0){#2}{\line(1,0){1.34}}
\put(#3,0){$\bullet$}
\end{picture}}

\newcommand{\dynkinapp}[4]{
\setlength{\unitlength}{2.5mm}
\begin{picture}(#1,3)
\put(0,0){$\circ$}
\multiput(2,0)(2,0){#2}{$\circ$}
\multiput(.73,.4)(2,0){#2}{\line(1,0){1.34}}
\put(#3,0){$\bullet$}
\put(#4,0){$*$}
\end{picture}}

\newcommand{\dynkindpspinp}[2]{
\setlength{\unitlength}{2.5mm}
\begin{picture}(#1,1.5)(.15,0)
\put(2,0){$\circ$}
\multiput(4,0)(2,0){#2}{$\circ$}
\multiput(2.73,.4)(2,0){#2}{\line(1,0){1.34}}
\put(-.15,1.17){$\bullet$}
\put(-.15,-1.17){$\circ$}
\put(.6,1.5){\line(5,-3){1.5}}
\put(.6,-.64){\line(5,3){1.5}}
\end{picture}
\vspace{.2cm}
}


\newcommand{\dynkine}[2]{
\setlength{\unitlength}{2.5mm}
\begin{picture}(#1,1.5)(0,0)
\put(0,0){$\circ$}
\multiput(2,0)(2,0){#2}{$\circ$}
\multiput(.63,.4)(2,0){#2}{\line(1,0){1.44}}
\put(4,-2){$\circ$}
\put(4.36,-1.45){\line(0,1){1.56}}
\end{picture}
\vspace{.3cm}}

\newcommand{\dynkinep}[3]{
\setlength{\unitlength}{2.5mm}
\begin{picture}(#1,1.5)(0,0)
\put(0,0){$\circ$}
\multiput(2,0)(2,0){#2}{$\circ$}
\multiput(.63,.4)(2,0){#2}{\line(1,0){1.44}}
\put(#3,0){$\bullet$}
\put(4,-2){$\circ$}
\put(4.36,-1.45){\line(0,1){1.56}}
\end{picture}
\vspace{.3cm}}

\newcommand{\dynkinepp}[4]{
\setlength{\unitlength}{2.5mm}
\begin{picture}(#1,1.5)(0,0)
\put(0,0){$\circ$}
\multiput(2,0)(2,0){#2}{$\circ$}
\multiput(.63,.4)(2,0){#2}{\line(1,0){1.44}}
\put(#3,0){$\bullet$}
\put(4,-2){$\circ$}
\put(4.36,-1.45){\line(0,1){1.56}}
\put(#4,0){$*$}
\end{picture}
\vspace{.2cm}}

\newcommand{\dynkineDeuxMarque}[4]{
\setlength{\unitlength}{2.5mm}
\begin{picture}(#1,1.5)(0,0)
\put(0,0){$\circ$}
\multiput(2,0)(2,0){#2}{$\circ$}
\multiput(.63,.4)(2,0){#2}{\line(1,0){1.44}}
\put(#3,0){$\bullet$}
\put(#4,0){$\bullet$}
\put(4,-2){$\circ$}
\put(4.36,-1.45){\line(0,1){1.56}}
\end{picture}
\vspace{.2cm}}

\newcommand{\dynkineUneMarqueP}[3]{
\setlength{\unitlength}{2.5mm}
\begin{picture}(#1,1.5)(0,0)
\put(0,0){$\circ$}
\multiput(2,0)(2,0){#2}{$\circ$}
\multiput(.63,.4)(2,0){#2}{\line(1,0){1.44}}
\put(#3,0){$\bullet$}
\put(4,-2){$\bullet$}
\put(4.36,-1.45){\line(0,1){1.56}}
\end{picture}
\vspace{.2cm}}

\newcommand{\poidsesix}[6]{
\hspace{-.12cm}
\left [
\begin{array}{ccccc}
{} \hspace{-.2cm} #1 & {} \hspace{-.3cm} #2 & {} \hspace{-.3cm} #3 &
{} \hspace{-.3cm} #4 & {} \hspace{-.3cm} #5 \vspace{-.13cm}\\
\hspace{-.2cm} & \hspace{-.3cm} & {} \hspace{-.3cm} #6 &
{} \hspace{-.3cm} & {} \hspace{-.3cm}
\end{array}
\hspace{-.2cm}
\right ]      }


\newcommand{\spec}{\mathtt{Spec}}
\newcommand{\spe}{\mathtt{Spec}\ }
\newcommand{\proj}{\mathtt{Proj}}
\newcommand{\sz}{{\spec\ \Z}}
\newcommand{\p}{{\mathbb P}}
\renewcommand{\P}{{\mathbb P}}
\newcommand{\A}{{\mathbb A}}
\newcommand{\pz}{\p_\Z}
\newcommand{\cl}{{\cal L}}
\newcommand{\cm}{{\cal M}}
\newcommand{\co}{{\cal O}}
\newcommand{\cf}{{\cal F}}
\newcommand{\cg}{{\cal G}}
\newcommand{\cq}{{\cal Q}}
\newcommand{\Hom}{{\cal H}om}
\newcommand{\End}{{\cal E}nd}


\newcommand{\ssi}{si et seulement si }
\renewcommand{\iff}{if and only if }
\newcommand{\tr}{{}^t}
\newcommand{\trace}{\mathtt{tr}}
\newcommand{\scal}[1]{\langle #1 \rangle}
\newcommand{\im}{\mathtt{Im}}
\renewcommand{\det}{\mathtt{det}}
\newcommand{\Det}{\mathtt{Det}}


\newcommand{\suiteexacte}[3]{#1 \rightarrow #2 \rightarrow #3}
\newcommand{\surmap}{\rightarrow \hspace{-.5cm} \rightarrow}
\newcommand{\limiteinverse}{\lim_\leftarrow}
\newcommand{\liste}{\

\begin{itemize}}
\newcommand{\codim}{\mbox{codim}}
\newcommand{\point}{^{ ^\bullet} \hspace{-.7mm}}
\newcommand{\X}{\mathfrak X}

\newcommand{\res}[2]{\vspace{.15cm} 

\noindent
{\bf #1 :} {\it #2} \vspace{.15cm} 

\noindent}

\newcommand{\fonction}[5]{
\begin{array}[t]{rrcll}
#1 & : & #2 & \rightarrow & #3 \\
   &   & #4 & \mapsto     & #5
\end{array}  }

\newcommand{\fonctionratsansnom}[4]{
\begin{array}[t]{ccl}
#1 & \dasharrow & #2 \\
#3 & \mapsto    & #4
\end{array}  }

\newcommand{\fonc}[3]{
#1 : #2 \mapsto #3  }

\newcommand{\directlim}[1]{
\lim_{\stackrel{\rightarrow}{#1}}    }

\newcommand{\inverselim}[1]{
\lim_{\stackrel{\rightarrow}{#1}}    }

\newcommand{\suitecourte}[3]{
0 \rightarrow #1 \rightarrow #2 \rightarrow #3 \rightarrow 0 }

\newcommand{\matdd}[4]{
\left (
\begin{array}{cc}
#1 & #2  \\
#3 & #4
\end{array}
\right )   }

\newcommand{\matddr}[4]{
\left (
\begin{array}{cc}
\hspace{-.2cm} #1 & \hspace{-.2cm}#2  \\
\hspace{-.2cm} #3 & \hspace{-.2cm}#4
\end{array}
\hspace{-.2cm} \right )   }

\newcommand{\mattt}[9]{
\left (
\begin{array}{ccc}
{} #1 & {} #2 & {} #3 \\
{} #4 & {} #5 & {} #6 \\
{} #7 & {} #8 & {} #9
\end{array}
\right )   }

\newcommand{\matttr}[9]{
\left (
\begin{array}{ccc}
{} \hspace{-.2cm} #1 & {} \hspace{-.2cm} #2 & {} \hspace{-.2cm} #3 \\
{} \hspace{-.2cm} #4 & {} \hspace{-.2cm} #5 & {} \hspace{-.2cm} #6 \\
{} \hspace{-.2cm} #7 & {} \hspace{-.2cm} #8 & {} \hspace{-.2cm} #9
\end{array}
\hspace{-.15cm}
\right )   }

%% file: g37.tex

\psset{unit=6mm}

\begin{pspicture*}(-0.5,-0.5)(5.5,5.5)

\psline(1.88,4.88)(1.12,4.12)
\psline(2.12,4.88)(2.88,4.12)
\psline(0.88,3.88)(0.12,3.12)
\psline(1.12,3.88)(1.88,3.12)
\psline(2.88,3.88)(2.12,3.12)
\psline(3.12,3.88)(3.88,3.12)
\psline(0.12,2.88)(0.88,2.12)
\psline(1.88,2.88)(1.12,2.12)
\psline(2.12,2.88)(2.88,2.12)
\psline(3.88,2.88)(3.12,2.12)
\psline(4.12,2.88)(4.88,2.12)
\psline(1.12,1.88)(1.88,1.12)
\psline(2.88,1.88)(2.12,1.12)
\psline(3.12,1.88)(3.88,1.12)
\psline(4.88,1.88)(4.12,1.12)
\psline(2.12,0.88)(2.88,0.12)
\psline(3.88,0.88)(3.12,0.12)

\pscircle*(0,3){0.16}
\pscircle*(1,4){0.16}
\pscircle*(1,2){0.16}
\pscircle*(2,5){0.16}
\pscircle*(2,3){0.16}
\pscircle*(2,1){0.16}
\pscircle*(3,4){0.16}
\pscircle*(3,2){0.16}
\pscircle*(3,0){0.16}
\pscircle*(4,3){0.16}
\pscircle*(4,1){0.16}
\pscircle*(5,2){0.16}

\end{pspicture*}

%% file: gq5.tex

\psset{unit=6mm}

\begin{pspicture*}(-0.5,-0.5)(3.5,6.5)

\psline(2.88,5.88)(2.12,5.12)
\psline(1.88,4.88)(1.12,4.12)
\psline(2.12,4.88)(2.88,4.12)
\psline(0.88,3.88)(0.12,3.12)
\psline(1.12,3.88)(1.88,3.12)
\psline(2.88,3.88)(2.12,3.12)
\psline(0.12,2.88)(0.88,2.12)
\psline(1.88,2.88)(1.12,2.12)
\psline(2.12,2.88)(2.88,2.12)
\psline(1.12,1.88)(1.88,1.12)
\psline(2.88,1.88)(2.12,1.12)
\psline(2.12,0.88)(2.88,0.12)

\pscircle*(0,3){0.16}
\pscircle*(1,4){0.16}
\pscircle*(1,2){0.16}
\pscircle*(2,5){0.16}
\pscircle*(2,3){0.16}
\pscircle*(2,1){0.16}
\pscircle(3,6){0.16}  \psline(2.88,5.88)(3.12,6.12)  \psline(2.88,6.12)(3.12,5.88)
\pscircle(3,4){0.16}  \psline(2.84,4)(3.16,4)  \psline(3,3.84)(3,4.16)
\pscircle(3,2){0.16}  \psline(2.88,1.88)(3.12,2.12)  \psline(2.88,2.12)(3.12,1.88)
\pscircle(3,0){0.16}  \psline(2.84,0)(3.16,0)  \psline(3,-0.16)(3,0.16)

\end{pspicture*}

%% file: e6.tex

\psset{unit=6mm}

\begin{pspicture*}(-0.5,-0.5)(4.5,10.5)

\psline(3.88,9.88)(3.12,9.12)
\psline(2.88,8.88)(2.12,8.12)
\psline(1.88,7.88)(1.12,7.12)
\psline(2,7.84)(2,7.16)
\psline(0.88,6.88)(0.12,6.12)
\psline(1.12,6.88)(1.88,6.12)
\psline(2,6.84)(2,6.16)
\psline(0.12,5.88)(0.88,5.12)
\psline(1.88,5.88)(1.12,5.12)
\psline(2.12,5.88)(2.88,5.12)
\psline(1.12,4.88)(1.88,4.12)
\psline(2.88,4.88)(2.12,4.12)
\psline(3.12,4.88)(3.88,4.12)
\psline(2,3.84)(2,3.16)
\psline(2.12,3.88)(2.88,3.12)
\psline(3.88,3.88)(3.12,3.12)
\psline(2,2.84)(2,2.16)
\psline(2.88,2.88)(2.12,2.12)
\psline(1.88,1.88)(1.12,1.12)
\psline(0.88,0.88)(0.12,0.12)

\pscircle*(0,6){0.16}
\pscircle*(0,0){0.16}
\pscircle*(1,7){0.16}
\pscircle*(1,5){0.16}
\pscircle*(1,1){0.16}
\pscircle*(2,8){0.16}
\pscircle*(2,7){0.16}
\pscircle*(2,6){0.16}
\pscircle*(2,4){0.16}
\pscircle*(2,3){0.16}
\pscircle*(2,2){0.16}
\pscircle*(3,9){0.16}
\pscircle*(3,5){0.16}
\pscircle*(3,3){0.16}
\pscircle*(4,10){0.16}
\pscircle*(4,4){0.16}

\end{pspicture*}

%% file: hasseq8.tex

\psset{unit=6mm}

\begin{pspicture*}(-0.5,-0.5)(4.5,6.5)

\psline(0.12,5.88)(0.88,5.12)
\psline(1.12,4.88)(1.88,4.12)
\psline(2,3.84)(2,3.16)
\psline(2.12,3.88)(2.88,3.12)
\psline(2,2.84)(2,2.16)
\psline(2.88,2.88)(2.12,2.12)
\psline(1.88,1.88)(1.12,1.12)
\psline(0.88,0.88)(0.12,0.12)

\pscircle*(0,6){0.16}   \put(0.208,5.859){$-\epsilon_1$}
\pscircle*(0,0){0.16}   \put(0.208,-0.141){$\epsilon_1$}
\pscircle*(1,5){0.16}   \put(1.208,4.859){$-\epsilon_2$}
\pscircle*(1,1){0.16}   \put(1.208,0.859){$\epsilon_2$}
\pscircle*(2,4){0.16}   \put(2.208,3.859){$-\epsilon_3$}
\pscircle*(2,3){0.16}   \put(0.792,2.859){$-\epsilon_4$}
\pscircle*(2,2){0.16}   \put(2.208,1.859){$\epsilon_3$}
\pscircle*(3,3){0.16}   \put(3.208,2.859){$\epsilon_4$}

\end{pspicture*}

%% file: racineq10.tex

\psset{unit=6mm}

\begin{pspicture*}(-0.5,-0.5)(6.5,6.5)

\psline(0.12,5.88)(0.88,5.12)
\psline(1.12,4.88)(1.88,4.12)
\psline(2,3.84)(2,3.16)
\psline(2.12,3.88)(2.88,3.12)
\psline(2,2.84)(2,2.16)
\psline(2.88,2.88)(2.12,2.12)
\psline(1.88,1.88)(1.12,1.12)
\psline(0.88,0.88)(0.12,0.12)

\pscircle*(0,6){0.16}   \put(0.208,5.859){$\epsilon_0+\epsilon_1$}
\pscircle*(0,0){0.16}   \put(0.208,-0.141){$\epsilon_0-\epsilon_1$}
\pscircle*(1,5){0.16}   \put(1.208,4.859){$\epsilon_0+\epsilon_2$}
\pscircle*(1,1){0.16}   \put(1.208,0.859){$\epsilon_0-\epsilon_2$}
\pscircle*(2,4){0.16}   \put(2.208,3.859){$\epsilon_0+\epsilon_3$}
\pscircle*(2,3){0.16}   \put(0.042,2.859){$\epsilon_0+\epsilon_4$}
\pscircle*(2,2){0.16}   \put(2.208,1.859){$\epsilon_0-\epsilon_3$}
\pscircle*(3,3){0.16}   \put(3.208,2.859){$\epsilon_0-\epsilon_4$}

\end{pspicture*}

%% file: w_max_1.tex

\psset{unit=6mm}

\begin{pspicture*}(-0.5,-0.5)(4.5,10.5)

\psline(3.88,9.88)(3.12,9.12)
\psline(2.88,8.88)(2.12,8.12)
\psline(1.88,7.88)(1.12,7.12)
\psline(2,7.84)(2,7.16)
\psline(0.88,6.88)(0.12,6.12)
\psline(1.12,6.88)(1.88,6.12)
\psline(2,6.84)(2,6.16)
\psline(0.12,5.88)(0.88,5.12)
\psline(1.88,5.88)(1.12,5.12)
\psline(2.12,5.88)(2.88,5.12)
\psline(1.12,4.88)(1.88,4.12)
\psline(2.88,4.88)(2.12,4.12)
\psline(3.12,4.88)(3.88,4.12)
\psline(2,3.84)(2,3.16)
\psline(2.12,3.88)(2.88,3.12)
\psline(3.88,3.88)(3.12,3.12)
\psline(2,2.84)(2,2.16)
\psline(2.88,2.88)(2.12,2.12)
\psline(1.88,1.88)(1.12,1.12)
\psline(0.88,0.88)(0.12,0.12)

\pscircle*(0,6){0.16}
\pscircle*(0,0){0.16}
\pscircle(1,7){0.16}
\pscircle*(1,5){0.16}
\pscircle*(1,1){0.16}
\pscircle(2,8){0.16}
\pscircle(2,7){0.16}
\pscircle*(2,6){0.16}
\pscircle*(2,4){0.16}
\pscircle*(2,3){0.16}
\pscircle*(2,2){0.16}
\pscircle(3,9){0.16}
\pscircle*(3,5){0.16}
\pscircle*(3,3){0.16}
\pscircle(4,10){0.16}
\pscircle*(4,4){0.16}

\end{pspicture*}

%% file: h.tex

\psset{unit=6mm}

\begin{pspicture*}(-2.7,-0.5)(6.7,10.5)

\psline(3.88,9.88)(3.12,9.12)
\psline(2.88,8.88)(2.12,8.12)
\psline(1.88,7.88)(1.12,7.12)
\psline(2,7.84)(2,7.16)
\psline(0.88,6.88)(0.12,6.12)
\psline(1.12,6.88)(1.88,6.12)
\psline(2,6.84)(2,6.16)
\psline(0.12,5.88)(0.88,5.12)
\psline(1.88,5.88)(1.12,5.12)
\psline(2.12,5.88)(2.88,5.12)
\psline(1.12,4.88)(1.88,4.12)
\psline(2.88,4.88)(2.12,4.12)
\psline(3.12,4.88)(3.88,4.12)
\psline(2,3.84)(2,3.16)
\psline(2.12,3.88)(2.88,3.12)
\psline(3.88,3.88)(3.12,3.12)
\psline(2,2.84)(2,2.16)
\psline(2.88,2.88)(2.12,2.12)
\psline(1.88,1.88)(1.12,1.12)
\psline(0.88,0.88)(0.12,0.12)

\pscircle*(0,6){0.16}   \put(-2.541,5.859){\scriptsize{$+---+$}}
\pscircle*(0,0){0.16}   \put(0.208,-0.141){\scriptsize{$++++-$}}
\pscircle*(1,7){0.16}   \put(-1.541,6.859){\scriptsize{$-+--+$}}
\pscircle*(1,5){0.16}   \put(-1.541,4.859){\scriptsize{$+--+-$}}
\pscircle*(1,1){0.16}   \put(1.208,0.859){\scriptsize{$+++-+$}}
\pscircle*(2,8){0.16}   \put(-0.541,7.859){\scriptsize{$--+-+$}}
\pscircle*(2,7){0.16}   \put(2.208,6.859){\scriptsize{$--++-$}}
\pscircle*(2,6){0.16}   \put(2.208,5.859){\scriptsize{$-+-+-$}}
\pscircle*(2,4){0.16}   \put(-0.541,3.859){\scriptsize{$+-+--$}}
\pscircle*(2,3){0.16}   \put(-0.541,2.859){\scriptsize{$++---$}}
\pscircle*(2,2){0.16}   \put(2.208,1.859){\scriptsize{$++-++$}}
\pscircle*(3,9){0.16}   \put(0.459,8.859){\scriptsize{$---++$}}
\pscircle*(3,5){0.16}   \put(3.208,4.859){\scriptsize{$-++--$}}
\pscircle*(3,3){0.16}   \put(3.208,2.859){\scriptsize{$+-+++$}}
\pscircle*(4,10){0.16}   \put(1.459,9.859){\scriptsize{$-----$}}
\pscircle*(4,4){0.16}   \put(4.208,3.859){\scriptsize{$-++++$}}

\end{pspicture*}

%% file: hd.tex

\psset{unit=6mm}

\begin{pspicture*}(-2.7,-0.5)(6.7,10.5)

\psline(3.88,9.88)(3.12,9.12)
\psline(2.88,8.88)(2.12,8.12)
\psline(1.88,7.88)(1.12,7.12)
\psline(2,7.84)(2,7.16)
\psline(0.88,6.88)(0.12,6.12)
\psline(1.12,6.88)(1.88,6.12)
\psline(2,6.84)(2,6.16)
\psline(0.12,5.88)(0.88,5.12)
\psline(1.88,5.88)(1.12,5.12)
\psline(2.12,5.88)(2.88,5.12)
\psline(1.12,4.88)(1.88,4.12)
\psline(2.88,4.88)(2.12,4.12)
\psline(3.12,4.88)(3.88,4.12)
\psline(2,3.84)(2,3.16)
\psline(2.12,3.88)(2.88,3.12)
\psline(3.88,3.88)(3.12,3.12)
\psline(2,2.84)(2,2.16)
\psline(2.88,2.88)(2.12,2.12)
\psline(1.88,1.88)(1.12,1.12)
\psline(0.88,0.88)(0.12,0.12)

\pscircle*(0,6){0.16}   \put(-2.541,5.859){\scriptsize{$+----$}}
\pscircle*(0,0){0.16}   \put(0.208,-0.141){\scriptsize{$+++++$}}
\pscircle*(1,7){0.16}   \put(-1.541,6.859){\scriptsize{$-+---$}}
\pscircle*(1,5){0.16}   \put(-1.541,4.859){\scriptsize{$+--++$}}
\pscircle*(1,1){0.16}   \put(1.208,0.859){\scriptsize{$+++--$}}
\pscircle*(2,8){0.16}   \put(-0.541,7.859){\scriptsize{$--+--$}}
\pscircle*(2,7){0.16}   \put(2.208,6.859){\scriptsize{$--+++$}}
\pscircle*(2,6){0.16}   \put(2.208,5.859){\scriptsize{$-+-++$}}
\pscircle*(2,4){0.16}   \put(-0.541,3.859){\scriptsize{$+-+-+$}}
\pscircle*(2,3){0.16}   \put(-0.541,2.859){\scriptsize{$++--+$}}
\pscircle*(2,2){0.16}   \put(2.208,1.859){\scriptsize{$++-+-$}}
\pscircle*(3,9){0.16}   \put(0.459,8.859){\scriptsize{$---+-$}}
\pscircle*(3,5){0.16}   \put(3.208,4.859){\scriptsize{$-++-+$}}
\pscircle*(3,3){0.16}   \put(3.208,2.859){\scriptsize{$+-++-$}}
\pscircle*(4,10){0.16}   \put(1.459,9.859){\scriptsize{$----+$}}
\pscircle*(4,4){0.16}   \put(4.208,3.859){\scriptsize{$-+++-$}}

\end{pspicture*}

%% file: cas0.tex

\psset{unit=5mm}

\begin{pspicture*}(-0.5,-0.5)(4.5,10.5)

\psline(3.88,9.88)(3.12,9.12)
\psline(2.88,8.88)(2.12,8.12)
\psline(1.88,7.88)(1.12,7.12)
\psline(2,7.84)(2,7.16)
\psline(0.88,6.88)(0.12,6.12)
\psline(1.12,6.88)(1.88,6.12)
\psline(2,6.84)(2,6.16)
\psline(0.12,5.88)(0.88,5.12)
\psline(1.88,5.88)(1.12,5.12)
\psline(2.12,5.88)(2.88,5.12)
\psline(1.12,4.88)(1.88,4.12)
\psline(2.88,4.88)(2.12,4.12)
\psline(3.12,4.88)(3.88,4.12)
\psline(2,3.84)(2,3.16)
\psline(2.12,3.88)(2.88,3.12)
\psline(3.88,3.88)(3.12,3.12)
\psline(2,2.84)(2,2.16)
\psline(2.88,2.88)(2.12,2.12)
\psline(1.88,1.88)(1.12,1.12)
\psline(0.88,0.88)(0.12,0.12)

\pscircle*(0,6){0.16}
\pscircle(0,0){0.16}
\pscircle*(1,7){0.16}
\pscircle*(1,5){0.16}
\pscircle(1,1){0.16}
\pscircle*(2,8){0.16}
\pscircle*(2,7){0.16}
\pscircle*(2,6){0.16}
\pscircle*(2,4){0.16}
\pscircle(2,3){0.16}
\pscircle(2,2){0.16}
\pscircle*(3,9){0.16}
\pscircle*(3,5){0.16}
\pscircle(3,3){0.16}
\pscircle*(4,10){0.16}
\pscircle*(4,4){0.16}

\end{pspicture*}

%% file: cas1.tex

\psset{unit=5mm}

\begin{pspicture*}(-0.5,-0.5)(4.5,10.5)

\psline(3.88,9.88)(3.12,9.12)
\psline(2.88,8.88)(2.12,8.12)
\psline(1.88,7.88)(1.12,7.12)
\psline(2,7.84)(2,7.16)
\psline(0.88,6.88)(0.12,6.12)
\psline(1.12,6.88)(1.88,6.12)
\psline(2,6.84)(2,6.16)
\psline(0.12,5.88)(0.88,5.12)
\psline(1.88,5.88)(1.12,5.12)
\psline(2.12,5.88)(2.88,5.12)
\psline(1.12,4.88)(1.88,4.12)
\psline(2.88,4.88)(2.12,4.12)
\psline(3.12,4.88)(3.88,4.12)
\psline(2,3.84)(2,3.16)
\psline(2.12,3.88)(2.88,3.12)
\psline(3.88,3.88)(3.12,3.12)
\psline(2,2.84)(2,2.16)
\psline(2.88,2.88)(2.12,2.12)
\psline(1.88,1.88)(1.12,1.12)
\psline(0.88,0.88)(0.12,0.12)

\pscircle*(0,6){0.16}
\pscircle(0,0){0.16}
\pscircle*(1,7){0.16}
\pscircle*(1,5){0.16}
\pscircle(1,1){0.16}
\pscircle*(2,8){0.16}
\pscircle*(2,7){0.16}
\pscircle*(2,6){0.16}
\pscircle*(2,4){0.16}
\pscircle(2,3){0.16}
\pscircle(2,2){0.16}
\pscircle*(3,9){0.16}
\pscircle*(3,5){0.16}
\pscircle*(3,3){0.16}
\pscircle*(4,10){0.16}
\pscircle*(4,4){0.16}

\end{pspicture*}

%% file: cas2.tex

\psset{unit=5mm}

\begin{pspicture*}(-0.5,-0.5)(4.5,10.5)

\psline(3.88,9.88)(3.12,9.12)
\psline(2.88,8.88)(2.12,8.12)
\psline(1.88,7.88)(1.12,7.12)
\psline(2,7.84)(2,7.16)
\psline(0.88,6.88)(0.12,6.12)
\psline(1.12,6.88)(1.88,6.12)
\psline(2,6.84)(2,6.16)
\psline(0.12,5.88)(0.88,5.12)
\psline(1.88,5.88)(1.12,5.12)
\psline(2.12,5.88)(2.88,5.12)
\psline(1.12,4.88)(1.88,4.12)
\psline(2.88,4.88)(2.12,4.12)
\psline(3.12,4.88)(3.88,4.12)
\psline(2,3.84)(2,3.16)
\psline(2.12,3.88)(2.88,3.12)
\psline(3.88,3.88)(3.12,3.12)
\psline(2,2.84)(2,2.16)
\psline(2.88,2.88)(2.12,2.12)
\psline(1.88,1.88)(1.12,1.12)
\psline(0.88,0.88)(0.12,0.12)

\pscircle*(0,6){0.16}
\pscircle(0,0){0.16}
\pscircle*(1,7){0.16}
\pscircle*(1,5){0.16}
\pscircle(1,1){0.16}
\pscircle*(2,8){0.16}
\pscircle*(2,7){0.16}
\pscircle*(2,6){0.16}
\pscircle*(2,4){0.16}
\pscircle*(2,3){0.16}
\pscircle(2,2){0.16}
\pscircle*(3,9){0.16}
\pscircle*(3,5){0.16}
\pscircle(3,3){0.16}
\pscircle*(4,10){0.16}
\pscircle(4,4){0.16}

\end{pspicture*}

%% file: p5.tex

\psset{unit=5mm}

\begin{pspicture*}(-0.5,-0.5)(4.5,10.5)

\psline(3.88,9.88)(3.12,9.12)
\psline(2.88,8.88)(2.12,8.12)
\psline(1.88,7.88)(1.12,7.12)
\psline(2,7.84)(2,7.16)
\psline(0.88,6.88)(0.12,6.12)
\psline(1.12,6.88)(1.88,6.12)
\psline(2,6.84)(2,6.16)
\psline(0.12,5.88)(0.88,5.12)
\psline(1.88,5.88)(1.12,5.12)
\psline(2.12,5.88)(2.88,5.12)
\psline(1.12,4.88)(1.88,4.12)
\psline(2.88,4.88)(2.12,4.12)
\psline(3.12,4.88)(3.88,4.12)
\psline(2,3.84)(2,3.16)
\psline(2.12,3.88)(2.88,3.12)
\psline(3.88,3.88)(3.12,3.12)
\psline(2,2.84)(2,2.16)
\psline(2.88,2.88)(2.12,2.12)
\psline(1.88,1.88)(1.12,1.12)
\psline(0.88,0.88)(0.12,0.12)

\pscircle(0,6){0.16}
\pscircle*(0,0){0.16}
\pscircle(1,7){0.16}
\pscircle(1,5){0.16}
\pscircle*(1,1){0.16}
\pscircle(2,8){0.16}
\pscircle(2,7){0.16}
\pscircle(2,6){0.16}
\pscircle(2,4){0.16}
\pscircle(2,3){0.16}
\pscircle*(2,2){0.16}
\pscircle(3,9){0.16}
\pscircle(3,5){0.16}
\pscircle*(3,3){0.16}
\pscircle(4,10){0.16}
\pscircle*(4,4){0.16}

\end{pspicture*}

%% file: ip5.tex

\psset{unit=5mm}

\begin{pspicture*}(-0.5,-0.5)(4.5,10.5)

\psline(3.88,9.88)(3.12,9.12)
\psline(2.88,8.88)(2.12,8.12)
\psline(1.88,7.88)(1.12,7.12)
\psline(2,7.84)(2,7.16)
\psline(0.88,6.88)(0.12,6.12)
\psline(1.12,6.88)(1.88,6.12)
\psline(2,6.84)(2,6.16)
\psline(0.12,5.88)(0.88,5.12)
\psline(1.88,5.88)(1.12,5.12)
\psline(2.12,5.88)(2.88,5.12)
\psline(1.12,4.88)(1.88,4.12)
\psline(2.88,4.88)(2.12,4.12)
\psline(3.12,4.88)(3.88,4.12)
\psline(2,3.84)(2,3.16)
\psline(2.12,3.88)(2.88,3.12)
\psline(3.88,3.88)(3.12,3.12)
\psline(2,2.84)(2,2.16)
\psline(2.88,2.88)(2.12,2.12)
\psline(1.88,1.88)(1.12,1.12)
\psline(0.88,0.88)(0.12,0.12)

\pscircle*(0,6){0.16}
\pscircle(0,0){0.16}
\pscircle*(1,7){0.16}
\pscircle(1,5){0.16}
\pscircle(1,1){0.16}
\pscircle*(2,8){0.16}
\pscircle(2,7){0.16}
\pscircle(2,6){0.16}
\pscircle(2,4){0.16}
\pscircle(2,3){0.16}
\pscircle(2,2){0.16}
\pscircle*(3,9){0.16}
\pscircle(3,5){0.16}
\pscircle(3,3){0.16}
\pscircle*(4,10){0.16}
\pscircle(4,4){0.16}

\end{pspicture*}

%% file: dp5.tex

\psset{unit=5mm}

\begin{pspicture*}(-0.5,-0.5)(4.5,10.5)

\psline(0.12,9.88)(0.88,9.12)
\psline(1.12,8.88)(1.88,8.12)
\psline(2,7.84)(2,7.16)
\psline(2.12,7.88)(2.88,7.12)
\psline(2,6.84)(2,6.16)
\psline(2.88,6.88)(2.12,6.12)
\psline(3.12,6.88)(3.88,6.12)
\psline(1.88,5.88)(1.12,5.12)
\psline(2.12,5.88)(2.88,5.12)
\psline(3.88,5.88)(3.12,5.12)
\psline(0.88,4.88)(0.12,4.12)
\psline(1.12,4.88)(1.88,4.12)
\psline(2.88,4.88)(2.12,4.12)
\psline(0.12,3.88)(0.88,3.12)
\psline(1.88,3.88)(1.12,3.12)
\psline(2,3.84)(2,3.16)
\psline(1.12,2.88)(1.88,2.12)
\psline(2,2.84)(2,2.16)
\psline(2.12,1.88)(2.88,1.12)
\psline(3.12,0.88)(3.88,0.12)

\pscircle(0,10){0.16}
\pscircle*(0,4){0.16}
\pscircle(1,9){0.16}
\pscircle(1,5){0.16}
\pscircle*(1,3){0.16}
\pscircle(2,8){0.16}
\pscircle(2,7){0.16}
\pscircle(2,6){0.16}
\pscircle(2,4){0.16}
\pscircle(2,3){0.16}
\pscircle*(2,2){0.16}
\pscircle(3,7){0.16}
\pscircle(3,5){0.16}
\pscircle*(3,1){0.16}
\pscircle(4,6){0.16}
\pscircle*(4,0){0.16}

\end{pspicture*}

%% file: hasse_e62.tex

\psset{unit=7.5mm}

\begin{pspicture*}(-0.5,-0.5)(3.5,7.5)

\psline[linecolor=red](2.88,6.88)(2.12,6.12)
\psline(3,6.84)(3,6.16)
\psline[linecolor=red](1.88,5.88)(1.12,5.12)
\psline(2,5.84)(2,5.16)
\psline[linecolor=red](2.12,5.88)(2.88,5.12)
\psline[linecolor=blue](2.88,5.88)(2.12,5.12)
\psline[linecolor=red](0.88,4.88)(0.12,4.12)
\psline(1,4.84)(1,4.16)
\psline[linecolor=red](1.12,4.88)(1.88,4.12)
\psline[linecolor=blue](1.88,4.88)(1.12,4.12)
\psline[linecolor=blue](2.12,4.88)(2.88,4.12)
\psline[linecolor=red](2.88,4.88)(2.12,4.12)
\psline(3,4.84)(3,4.16)
\psline(0,3.84)(0,3.16)
\psline[linecolor=red](0.12,3.88)(0.88,3.12)
\psline[linecolor=blue](0.88,3.88)(0.12,3.12)
\psline[linecolor=blue](1.12,3.88)(1.88,3.12)
\psline[linecolor=red](1.88,3.88)(1.12,3.12)
\psline(2,3.84)(2,3.16)
\psline[linecolor=red](2.12,3.88)(2.88,3.12)
\psline[linecolor=blue](2.88,3.88)(2.12,3.12)
\psline[linecolor=blue](0.12,2.88)(0.88,2.12)
\psline(1,2.84)(1,2.16)
\psline[linecolor=red](1.12,2.88)(1.88,2.12)
\psline[linecolor=blue](1.88,2.88)(1.12,2.12)
\psline[linecolor=blue](2.12,2.88)(2.88,2.12)
\psline[linecolor=red](2.88,2.88)(2.12,2.12)
\psline(3,2.84)(3,2.16)
\psline[linecolor=blue](1.12,1.88)(1.88,1.12)
\psline(2,1.84)(2,1.16)
\psline[linecolor=red](2.12,1.88)(2.88,1.12)
\psline[linecolor=blue](2.88,1.88)(2.12,1.12)
\psline[linecolor=blue](2.12,0.88)(2.88,0.12)
\psline(3,0.84)(3,0.16)

\pscircle*[linecolor=red](0,4){0.16}
\pscircle*[linecolor=blue](0,3){0.16}
\pscircle*[linecolor=red](1,5){0.16}
\pscircle*[linecolor=blue](1,4){0.16}
\pscircle*[linecolor=red](1,3){0.16}
\pscircle*[linecolor=blue](1,2){0.16}
\pscircle*[linecolor=red](2,6){0.16}
\pscircle*[linecolor=blue](2,5){0.16}
\pscircle*[linecolor=red](2,4){0.16}
\pscircle*[linecolor=blue](2,3){0.16}
\pscircle*[linecolor=red](2,2){0.16}
\pscircle*[linecolor=blue](2,1){0.16}
\pscircle*[linecolor=red](3,7){0.16}
\pscircle*[linecolor=blue](3,6){0.16}
\pscircle*[linecolor=red](3,5){0.16}
\pscircle*[linecolor=blue](3,4){0.16}
\pscircle*[linecolor=red](3,3){0.16}
\pscircle*[linecolor=blue](3,2){0.16}
\pscircle*[linecolor=red](3,1){0.16}
\pscircle*[linecolor=blue](3,0){0.16}

\end{pspicture*}

%% file: qmax.tex

\psset{unit=7.5mm}

\begin{pspicture*}(-0.5,-0.5)(3.5,7.5)

\psline[linecolor=red](2.88,6.88)(2.12,6.12)
\psline(3,6.84)(3,6.16)
\psline[linecolor=red](1.88,5.88)(1.12,5.12)
\psline(2,5.84)(2,5.16)
\psline[linecolor=red](2.12,5.88)(2.88,5.12)
\psline[linecolor=blue](2.88,5.88)(2.12,5.12)
\psline[linecolor=red](0.88,4.88)(0.12,4.12)
\psline(1,4.84)(1,4.16)
\psline[linecolor=red](1.12,4.88)(1.88,4.12)
\psline[linecolor=blue](1.88,4.88)(1.12,4.12)
\psline[linecolor=blue](2.12,4.88)(2.88,4.12)
\psline[linecolor=red](2.88,4.88)(2.12,4.12)
\psline(3,4.84)(3,4.16)
\psline(0,3.84)(0,3.16)
\psline[linecolor=red](0.12,3.88)(0.88,3.12)
\psline[linecolor=blue](0.88,3.88)(0.12,3.12)
\psline[linecolor=blue](1.12,3.88)(1.88,3.12)
\psline[linecolor=red](1.88,3.88)(1.12,3.12)
\psline(2,3.84)(2,3.16)
\psline[linecolor=red](2.12,3.88)(2.88,3.12)
\psline[linecolor=blue](2.88,3.88)(2.12,3.12)
\psline[linecolor=blue](0.12,2.88)(0.88,2.12)
\psline(1,2.84)(1,2.16)
\psline[linecolor=red](1.12,2.88)(1.88,2.12)
\psline[linecolor=blue](1.88,2.88)(1.12,2.12)
\psline[linecolor=blue](2.12,2.88)(2.88,2.12)
\psline[linecolor=red](2.88,2.88)(2.12,2.12)
\psline(3,2.84)(3,2.16)
\psline[linecolor=blue](1.12,1.88)(1.88,1.12)
\psline(2,1.84)(2,1.16)
\psline[linecolor=red](2.12,1.88)(2.88,1.12)
\psline[linecolor=blue](2.88,1.88)(2.12,1.12)
\psline[linecolor=blue](2.12,0.88)(2.88,0.12)
\psline(3,0.84)(3,0.16)

\pscircle[linecolor=red](0,4){0.16}
\pscircle*[linecolor=blue](0,3){0.16}
\pscircle[linecolor=red](1,5){0.16}
\pscircle[linecolor=blue](1,4){0.16}
\pscircle*[linecolor=red](1,3){0.16}
\pscircle*[linecolor=blue](1,2){0.16}
\pscircle[linecolor=red](2,6){0.16}
\pscircle[linecolor=blue](2,5){0.16}
\pscircle[linecolor=red](2,4){0.16}
\pscircle*[linecolor=blue](2,3){0.16}
\pscircle*[linecolor=red](2,2){0.16}
\pscircle*[linecolor=blue](2,1){0.16}
\pscircle[linecolor=red](3,7){0.16}
\pscircle[linecolor=blue](3,6){0.16}
\pscircle[linecolor=red](3,5){0.16}
\pscircle[linecolor=blue](3,4){0.16}
\pscircle*[linecolor=red](3,3){0.16}
\pscircle*[linecolor=blue](3,2){0.16}
\pscircle*[linecolor=red](3,1){0.16}
\pscircle*[linecolor=blue](3,0){0.16}

\end{pspicture*}

%% file: e6_2.tex

\psset{unit=7.5mm}

\begin{pspicture*}(-0.5,-0.5)(4,7.5)

\psline[linecolor=red](2.88,6.88)(2.12,6.12)
\psline(3,6.84)(3,6.16)
\psline[linecolor=red](1.88,5.88)(1.12,5.12)
\psline(2,5.84)(2,5.16)
\psline[linecolor=red](2.12,5.88)(2.88,5.12)
\psline[linecolor=blue](2.88,5.88)(2.12,5.12)
\psline[linecolor=red](0.88,4.88)(0.12,4.12)
\psline(1,4.84)(1,4.16)
\psline[linecolor=red](1.12,4.88)(1.88,4.12)
\psline[linecolor=blue](1.88,4.88)(1.12,4.12)
\psline[linecolor=blue](2.12,4.88)(2.88,4.12)
\psline[linecolor=red](2.88,4.88)(2.12,4.12)
\psline(3,4.84)(3,4.16)
\psline(0,3.84)(0,3.16)
\psline[linecolor=red](0.12,3.88)(0.88,3.12)
\psline[linecolor=blue](0.88,3.88)(0.12,3.12)
\psline[linecolor=blue](1.12,3.88)(1.88,3.12)
\psline[linecolor=red](1.88,3.88)(1.12,3.12)
\psline(2,3.84)(2,3.16)
\psline[linecolor=red](2.12,3.88)(2.88,3.12)
\psline[linecolor=blue](2.88,3.88)(2.12,3.12)
\psline[linecolor=blue](0.12,2.88)(0.88,2.12)
\psline(1,2.84)(1,2.16)
\psline[linecolor=red](1.12,2.88)(1.88,2.12)
\psline[linecolor=blue](1.88,2.88)(1.12,2.12)
\psline[linecolor=blue](2.12,2.88)(2.88,2.12)
\psline[linecolor=red](2.88,2.88)(2.12,2.12)
\psline(3,2.84)(3,2.16)
\psline[linecolor=blue](1.12,1.88)(1.88,1.12)
\psline(2,1.84)(2,1.16)
\psline[linecolor=red](2.12,1.88)(2.88,1.12)
\psline[linecolor=blue](2.88,1.88)(2.12,1.12)
\psline[linecolor=blue](2.12,0.88)(2.88,0.12)
\psline(3,0.84)(3,0.16)

\pscircle*[linecolor=red](0,4){0.16}   \put(0.208,3.887){15}
\pscircle*[linecolor=blue](0,3){0.16}   \put(0.208,2.887){\underline{15}}
\pscircle*[linecolor=red](1,5){0.16}   \put(1.208,4.887){14}
\pscircle*[linecolor=blue](1,4){0.16}   \put(1.208,3.887){\underline{14}}
\pscircle*[linecolor=red](1,3){0.16}   \put(1.208,2.887){25}
\pscircle*[linecolor=blue](1,2){0.16}   \put(1.208,1.887){\underline{25}}
\pscircle*[linecolor=red](2,6){0.16}   \put(2.208,5.887){13}
\pscircle*[linecolor=blue](2,5){0.16}   \put(2.208,4.887){\underline{13}}
\pscircle*[linecolor=red](2,4){0.16}   \put(2.208,3.887){24}
\pscircle*[linecolor=blue](2,3){0.16}   \put(2.208,2.887){\underline{24}}
\pscircle*[linecolor=red](2,2){0.16}   \put(2.208,1.887){35}
\pscircle*[linecolor=blue](2,1){0.16}   \put(2.208,0.887){\underline{35}}
\pscircle*[linecolor=red](3,7){0.16}   \put(3.208,6.887){12}
\pscircle*[linecolor=blue](3,6){0.16}   \put(3.208,5.887){\underline{12}}
\pscircle*[linecolor=red](3,5){0.16}   \put(3.208,4.887){23}
\pscircle*[linecolor=blue](3,4){0.16}   \put(3.208,3.887){\underline{23}}
\pscircle*[linecolor=red](3,3){0.16}   \put(3.208,2.887){34}
\pscircle*[linecolor=blue](3,2){0.16}   \put(3.208,1.887){\underline{34}}
\pscircle*[linecolor=red](3,1){0.16}   \put(3.208,0.887){45}
\pscircle*[linecolor=blue](3,0){0.16}   \put(3.208,-0.113){\underline{45}}

\end{pspicture*}

%% file: Cas0.tex

\psset{unit=7mm}

\begin{pspicture*}(-0.5,-0.5)(3.5,7.5)

\psline[linecolor=red](2.88,6.88)(2.12,6.12)
\psline(3,6.84)(3,6.16)
\psline[linecolor=red](1.88,5.88)(1.12,5.12)
\psline(2,5.84)(2,5.16)
\psline[linecolor=red](2.12,5.88)(2.88,5.12)
\psline[linecolor=blue](2.88,5.88)(2.12,5.12)
\psline[linecolor=red](0.88,4.88)(0.12,4.12)
\psline(1,4.84)(1,4.16)
\psline[linecolor=red](1.12,4.88)(1.88,4.12)
\psline[linecolor=blue](1.88,4.88)(1.12,4.12)
\psline[linecolor=blue](2.12,4.88)(2.88,4.12)
\psline[linecolor=red](2.88,4.88)(2.12,4.12)
\psline(3,4.84)(3,4.16)
\psline(0,3.84)(0,3.16)
\psline[linecolor=red](0.12,3.88)(0.88,3.12)
\psline[linecolor=blue](0.88,3.88)(0.12,3.12)
\psline[linecolor=blue](1.12,3.88)(1.88,3.12)
\psline[linecolor=red](1.88,3.88)(1.12,3.12)
\psline(2,3.84)(2,3.16)
\psline[linecolor=red](2.12,3.88)(2.88,3.12)
\psline[linecolor=blue](2.88,3.88)(2.12,3.12)
\psline[linecolor=blue](0.12,2.88)(0.88,2.12)
\psline(1,2.84)(1,2.16)
\psline[linecolor=red](1.12,2.88)(1.88,2.12)
\psline[linecolor=blue](1.88,2.88)(1.12,2.12)
\psline[linecolor=blue](2.12,2.88)(2.88,2.12)
\psline[linecolor=red](2.88,2.88)(2.12,2.12)
\psline(3,2.84)(3,2.16)
\psline[linecolor=blue](1.12,1.88)(1.88,1.12)
\psline(2,1.84)(2,1.16)
\psline[linecolor=red](2.12,1.88)(2.88,1.12)
\psline[linecolor=blue](2.88,1.88)(2.12,1.12)
\psline[linecolor=blue](2.12,0.88)(2.88,0.12)
\psline(3,0.84)(3,0.16)

\pscircle*[linecolor=red](0,4){0.16}
\pscircle[linecolor=blue](0,3){0.16}
\pscircle*[linecolor=red](1,5){0.16}
\pscircle*[linecolor=blue](1,4){0.16}
\pscircle[linecolor=red](1,3){0.16}
\pscircle[linecolor=blue](1,2){0.16}
\pscircle*[linecolor=red](2,6){0.16}
\pscircle*[linecolor=blue](2,5){0.16}
\pscircle*[linecolor=red](2,4){0.16}
\pscircle[linecolor=blue](2,3){0.16}
\pscircle[linecolor=red](2,2){0.16}
\pscircle[linecolor=blue](2,1){0.16}
\pscircle*[linecolor=red](3,7){0.16}
\pscircle*[linecolor=blue](3,6){0.16}
\pscircle*[linecolor=red](3,5){0.16}
\pscircle*[linecolor=blue](3,4){0.16}
\pscircle[linecolor=red](3,3){0.16}
\pscircle[linecolor=blue](3,2){0.16}
\pscircle[linecolor=red](3,1){0.16}
\pscircle[linecolor=blue](3,0){0.16}

\end{pspicture*}

%% file: Cas1.tex

\psset{unit=7mm}

\begin{pspicture*}(-0.5,-0.5)(3.5,7.5)

\psline[linecolor=red](2.88,6.88)(2.12,6.12)
\psline(3,6.84)(3,6.16)
\psline[linecolor=red](1.88,5.88)(1.12,5.12)
\psline(2,5.84)(2,5.16)
\psline[linecolor=red](2.12,5.88)(2.88,5.12)
\psline[linecolor=blue](2.88,5.88)(2.12,5.12)
\psline[linecolor=red](0.88,4.88)(0.12,4.12)
\psline(1,4.84)(1,4.16)
\psline[linecolor=red](1.12,4.88)(1.88,4.12)
\psline[linecolor=blue](1.88,4.88)(1.12,4.12)
\psline[linecolor=blue](2.12,4.88)(2.88,4.12)
\psline[linecolor=red](2.88,4.88)(2.12,4.12)
\psline(3,4.84)(3,4.16)
\psline(0,3.84)(0,3.16)
\psline[linecolor=red](0.12,3.88)(0.88,3.12)
\psline[linecolor=blue](0.88,3.88)(0.12,3.12)
\psline[linecolor=blue](1.12,3.88)(1.88,3.12)
\psline[linecolor=red](1.88,3.88)(1.12,3.12)
\psline(2,3.84)(2,3.16)
\psline[linecolor=red](2.12,3.88)(2.88,3.12)
\psline[linecolor=blue](2.88,3.88)(2.12,3.12)
\psline[linecolor=blue](0.12,2.88)(0.88,2.12)
\psline(1,2.84)(1,2.16)
\psline[linecolor=red](1.12,2.88)(1.88,2.12)
\psline[linecolor=blue](1.88,2.88)(1.12,2.12)
\psline[linecolor=blue](2.12,2.88)(2.88,2.12)
\psline[linecolor=red](2.88,2.88)(2.12,2.12)
\psline(3,2.84)(3,2.16)
\psline[linecolor=blue](1.12,1.88)(1.88,1.12)
\psline(2,1.84)(2,1.16)
\psline[linecolor=red](2.12,1.88)(2.88,1.12)
\psline[linecolor=blue](2.88,1.88)(2.12,1.12)
\psline[linecolor=blue](2.12,0.88)(2.88,0.12)
\psline(3,0.84)(3,0.16)

\pscircle*[linecolor=red](0,4){0.16}
\pscircle*[linecolor=blue](0,3){0.16}
\pscircle*[linecolor=red](1,5){0.16}
\pscircle*[linecolor=blue](1,4){0.16}
\pscircle[linecolor=red](1,3){0.16}
\pscircle[linecolor=blue](1,2){0.16}
\pscircle*[linecolor=red](2,6){0.16}
\pscircle*[linecolor=blue](2,5){0.16}
\pscircle[linecolor=red](2,4){0.16}
\pscircle[linecolor=blue](2,3){0.16}
\pscircle[linecolor=red](2,2){0.16}
\pscircle[linecolor=blue](2,1){0.16}
\pscircle*[linecolor=red](3,7){0.16}
\pscircle*[linecolor=blue](3,6){0.16}
\pscircle[linecolor=red](3,5){0.16}
\pscircle[linecolor=blue](3,4){0.16}
\pscircle[linecolor=red](3,3){0.16}
\pscircle[linecolor=blue](3,2){0.16}
\pscircle[linecolor=red](3,1){0.16}
\pscircle[linecolor=blue](3,0){0.16}

\end{pspicture*}

%% file: Cas2.tex

\psset{unit=7mm}

\begin{pspicture*}(-0.5,-0.5)(3.5,7.5)

\psline[linecolor=red](2.88,6.88)(2.12,6.12)
\psline(3,6.84)(3,6.16)
\psline[linecolor=red](1.88,5.88)(1.12,5.12)
\psline(2,5.84)(2,5.16)
\psline[linecolor=red](2.12,5.88)(2.88,5.12)
\psline[linecolor=blue](2.88,5.88)(2.12,5.12)
\psline[linecolor=red](0.88,4.88)(0.12,4.12)
\psline(1,4.84)(1,4.16)
\psline[linecolor=red](1.12,4.88)(1.88,4.12)
\psline[linecolor=blue](1.88,4.88)(1.12,4.12)
\psline[linecolor=blue](2.12,4.88)(2.88,4.12)
\psline[linecolor=red](2.88,4.88)(2.12,4.12)
\psline(3,4.84)(3,4.16)
\psline(0,3.84)(0,3.16)
\psline[linecolor=red](0.12,3.88)(0.88,3.12)
\psline[linecolor=blue](0.88,3.88)(0.12,3.12)
\psline[linecolor=blue](1.12,3.88)(1.88,3.12)
\psline[linecolor=red](1.88,3.88)(1.12,3.12)
\psline(2,3.84)(2,3.16)
\psline[linecolor=red](2.12,3.88)(2.88,3.12)
\psline[linecolor=blue](2.88,3.88)(2.12,3.12)
\psline[linecolor=blue](0.12,2.88)(0.88,2.12)
\psline(1,2.84)(1,2.16)
\psline[linecolor=red](1.12,2.88)(1.88,2.12)
\psline[linecolor=blue](1.88,2.88)(1.12,2.12)
\psline[linecolor=blue](2.12,2.88)(2.88,2.12)
\psline[linecolor=red](2.88,2.88)(2.12,2.12)
\psline(3,2.84)(3,2.16)
\psline[linecolor=blue](1.12,1.88)(1.88,1.12)
\psline(2,1.84)(2,1.16)
\psline[linecolor=red](2.12,1.88)(2.88,1.12)
\psline[linecolor=blue](2.88,1.88)(2.12,1.12)
\psline[linecolor=blue](2.12,0.88)(2.88,0.12)
\psline(3,0.84)(3,0.16)

\pscircle*[linecolor=red](0,4){0.16}
\pscircle[linecolor=blue](0,3){0.16}
\pscircle*[linecolor=red](1,5){0.16}
\pscircle*[linecolor=blue](1,4){0.16}
\pscircle*[linecolor=red](1,3){0.16}
\pscircle[linecolor=blue](1,2){0.16}
\pscircle*[linecolor=red](2,6){0.16}
\pscircle*[linecolor=blue](2,5){0.16}
\pscircle*[linecolor=red](2,4){0.16}
\pscircle[linecolor=blue](2,3){0.16}
\pscircle[linecolor=red](2,2){0.16}
\pscircle[linecolor=blue](2,1){0.16}
\pscircle*[linecolor=red](3,7){0.16}
\pscircle*[linecolor=blue](3,6){0.16}
\pscircle*[linecolor=red](3,5){0.16}
\pscircle[linecolor=blue](3,4){0.16}
\pscircle[linecolor=red](3,3){0.16}
\pscircle[linecolor=blue](3,2){0.16}
\pscircle[linecolor=red](3,1){0.16}
\pscircle[linecolor=blue](3,0){0.16}

\end{pspicture*}

%% file: Cas3.tex

\psset{unit=7mm}

\begin{pspicture*}(-0.5,-0.5)(3.5,7.5)

\psline[linecolor=red](2.88,6.88)(2.12,6.12)
\psline(3,6.84)(3,6.16)
\psline[linecolor=red](1.88,5.88)(1.12,5.12)
\psline(2,5.84)(2,5.16)
\psline[linecolor=red](2.12,5.88)(2.88,5.12)
\psline[linecolor=blue](2.88,5.88)(2.12,5.12)
\psline[linecolor=red](0.88,4.88)(0.12,4.12)
\psline(1,4.84)(1,4.16)
\psline[linecolor=red](1.12,4.88)(1.88,4.12)
\psline[linecolor=blue](1.88,4.88)(1.12,4.12)
\psline[linecolor=blue](2.12,4.88)(2.88,4.12)
\psline[linecolor=red](2.88,4.88)(2.12,4.12)
\psline(3,4.84)(3,4.16)
\psline(0,3.84)(0,3.16)
\psline[linecolor=red](0.12,3.88)(0.88,3.12)
\psline[linecolor=blue](0.88,3.88)(0.12,3.12)
\psline[linecolor=blue](1.12,3.88)(1.88,3.12)
\psline[linecolor=red](1.88,3.88)(1.12,3.12)
\psline(2,3.84)(2,3.16)
\psline[linecolor=red](2.12,3.88)(2.88,3.12)
\psline[linecolor=blue](2.88,3.88)(2.12,3.12)
\psline[linecolor=blue](0.12,2.88)(0.88,2.12)
\psline(1,2.84)(1,2.16)
\psline[linecolor=red](1.12,2.88)(1.88,2.12)
\psline[linecolor=blue](1.88,2.88)(1.12,2.12)
\psline[linecolor=blue](2.12,2.88)(2.88,2.12)
\psline[linecolor=red](2.88,2.88)(2.12,2.12)
\psline(3,2.84)(3,2.16)
\psline[linecolor=blue](1.12,1.88)(1.88,1.12)
\psline(2,1.84)(2,1.16)
\psline[linecolor=red](2.12,1.88)(2.88,1.12)
\psline[linecolor=blue](2.88,1.88)(2.12,1.12)
\psline[linecolor=blue](2.12,0.88)(2.88,0.12)
\psline(3,0.84)(3,0.16)

\pscircle[linecolor=red](0,4){0.16}
\pscircle[linecolor=blue](0,3){0.16}
\pscircle*[linecolor=red](1,5){0.16}
\pscircle*[linecolor=blue](1,4){0.16}
\pscircle[linecolor=red](1,3){0.16}
\pscircle[linecolor=blue](1,2){0.16}
\pscircle*[linecolor=red](2,6){0.16}
\pscircle*[linecolor=blue](2,5){0.16}
\pscircle*[linecolor=red](2,4){0.16}
\pscircle*[linecolor=blue](2,3){0.16}
\pscircle[linecolor=red](2,2){0.16}
\pscircle[linecolor=blue](2,1){0.16}
\pscircle*[linecolor=red](3,7){0.16}
\pscircle*[linecolor=blue](3,6){0.16}
\pscircle*[linecolor=red](3,5){0.16}
\pscircle*[linecolor=blue](3,4){0.16}
\pscircle[linecolor=red](3,3){0.16}
\pscircle[linecolor=blue](3,2){0.16}
\pscircle[linecolor=red](3,1){0.16}
\pscircle[linecolor=blue](3,0){0.16}

\end{pspicture*}

%% file: Cas4.tex

\psset{unit=7mm}

\begin{pspicture*}(-0.5,-0.5)(3.5,7.5)

\psline[linecolor=red](2.88,6.88)(2.12,6.12)
\psline(3,6.84)(3,6.16)
\psline[linecolor=red](1.88,5.88)(1.12,5.12)
\psline(2,5.84)(2,5.16)
\psline[linecolor=red](2.12,5.88)(2.88,5.12)
\psline[linecolor=blue](2.88,5.88)(2.12,5.12)
\psline[linecolor=red](0.88,4.88)(0.12,4.12)
\psline(1,4.84)(1,4.16)
\psline[linecolor=red](1.12,4.88)(1.88,4.12)
\psline[linecolor=blue](1.88,4.88)(1.12,4.12)
\psline[linecolor=blue](2.12,4.88)(2.88,4.12)
\psline[linecolor=red](2.88,4.88)(2.12,4.12)
\psline(3,4.84)(3,4.16)
\psline(0,3.84)(0,3.16)
\psline[linecolor=red](0.12,3.88)(0.88,3.12)
\psline[linecolor=blue](0.88,3.88)(0.12,3.12)
\psline[linecolor=blue](1.12,3.88)(1.88,3.12)
\psline[linecolor=red](1.88,3.88)(1.12,3.12)
\psline(2,3.84)(2,3.16)
\psline[linecolor=red](2.12,3.88)(2.88,3.12)
\psline[linecolor=blue](2.88,3.88)(2.12,3.12)
\psline[linecolor=blue](0.12,2.88)(0.88,2.12)
\psline(1,2.84)(1,2.16)
\psline[linecolor=red](1.12,2.88)(1.88,2.12)
\psline[linecolor=blue](1.88,2.88)(1.12,2.12)
\psline[linecolor=blue](2.12,2.88)(2.88,2.12)
\psline[linecolor=red](2.88,2.88)(2.12,2.12)
\psline(3,2.84)(3,2.16)
\psline[linecolor=blue](1.12,1.88)(1.88,1.12)
\psline(2,1.84)(2,1.16)
\psline[linecolor=red](2.12,1.88)(2.88,1.12)
\psline[linecolor=blue](2.88,1.88)(2.12,1.12)
\psline[linecolor=blue](2.12,0.88)(2.88,0.12)
\psline(3,0.84)(3,0.16)

\pscircle[linecolor=red](0,4){0.16}
\pscircle[linecolor=blue](0,3){0.16}
\pscircle*[linecolor=red](1,5){0.16}
\pscircle[linecolor=blue](1,4){0.16}
\pscircle[linecolor=red](1,3){0.16}
\pscircle[linecolor=blue](1,2){0.16}
\pscircle*[linecolor=red](2,6){0.16}
\pscircle*[linecolor=blue](2,5){0.16}
\pscircle*[linecolor=red](2,4){0.16}
\pscircle[linecolor=blue](2,3){0.16}
\pscircle[linecolor=red](2,2){0.16}
\pscircle[linecolor=blue](2,1){0.16}
\pscircle*[linecolor=red](3,7){0.16}
\pscircle*[linecolor=blue](3,6){0.16}
\pscircle*[linecolor=red](3,5){0.16}
\pscircle*[linecolor=blue](3,4){0.16}
\pscircle*[linecolor=red](3,3){0.16}
\pscircle[linecolor=blue](3,2){0.16}
\pscircle[linecolor=red](3,1){0.16}
\pscircle[linecolor=blue](3,0){0.16}

\end{pspicture*}